\documentclass[12pt]{amsart}
\usepackage{amssymb, amscd }
\usepackage{tikz-cd}

\usepackage{color}
\definecolor{darkgreen}{cmyk}{1,0,1,.2}
\definecolor{m}{rgb}{1,0.1,1}
\definecolor{green}{cmyk}{1,0,1,0}

\definecolor{test}{rgb}{1,0,0}
\definecolor{cmyk}{cmyk}{0,1,1,0}

\newtheorem{theorem}{Theorem}[section]

\newtheorem{proposition}[theorem]{Proposition}
\newtheorem{lemma}[theorem]{Lemma}

\newtheorem{remark}[theorem]{Remark}

\theoremstyle{remark}

\theoremstyle{definition}
\newtheorem{definition}[theorem]{Definition}

\numberwithin{equation}{subsection}

\def\CC{\operatorname{CC}}
\def\HC{\operatorname{HC}}

\def\ch{\operatorname{ch}}
\def\Ch{\operatorname{Ch}}
\def\Vect{\operatorname{Vect}}

\def\ind{\operatorname{ind}}
\def\End{\operatorname{End}}

\def\Hom{\operatorname{Hom}}

\def\Ind{\operatorname{Ind}}

\def\Ker{\operatorname{Ker}}

\def\Tr{\operatorname{Tr}}
\def\STr{\operatorname{Str}}
\def\tr{\operatorname{tr}}

\DeclareMathOperator{\str}{str}
\DeclareMathOperator{\id}{id}
\DeclareMathOperator{\pr}{pr}

\def\A{\mathbb A}

\def\C{\mathbb C}

\def\R{\mathbb R}

\def\Z{\mathbb Z}
\def\N{\mathbb N}

\def\cA{{\mathcal A}}
\def\maA{{\mathcal A}}
\def\maB{{\mathcal B}}
\def\maC{{\mathcal C}}
\def\maD{{\mathcal D}}
\def\maE{{\mathcal E}}
\def\maF{{\mathcal F}}

\def\maL{{\mathcal L}}

\def\maG{{\mathcal G}}

\def\tP{{\mathcal P}}

\def\cG{{\mathcal G}}
\def\maH{{\mathcal H}}
\def\cH{{\mathcal H}}

\def\cP{{\mathcal P}}

\def\maS{{\mathcal S}}

\def\what{\widehat}

\def\pa{\partial}

\DeclareMathOperator{\del}{\partial}
\DeclareMathOperator{\ord}{ord}
\def\cE{{\mathcal E}}
\def\cH{{\mathcal H}}
\def\cP{{\mathcal P}}
\def\cL{{\mathcal L}}
\def\u{u^*}

\renewcommand{\t}{\sigma}

\begin{document}

\title
{The higher twisted index theorem for foliations }


\author[M.-T. Benameur]{Moulay-Tahar Benameur}
\address{UMR  5149 du CNRS, Institut Montpelli\'erain Alexander Grothendieck, Montpellier, France}
\email{moulay.benameur@umontpellier.fr}
\author[A. Gorokhovsky ]{Alexander Gorokhovsky}
\address{Department of Mathematics, University of Colorado, Boulder, CO 80305-0395, USA}
\email{Alexander.Gorokhovsky@colorado.edu}
\author[E. Leichtnam]{Eric Leichtnam}
\address{UMR 7586 du CNRS, Institut Math\'ematique de Jussieu-PRG,  Paris, France}
\email{eric.leichtnam@imj-prg.fr}

\thanks{The second author was partially supported
by the NSF grant DMS-0900968.
}


\begin{abstract}\ Given a gerbe $L$, on the holonomy groupoid $\maG$ of the foliation $(M, \maF)$, whose pull-back to $M$ is torsion, we construct a Connes $\Phi$-map from the twisted Dupont-Sullivan bicomplex of $\maG$ to the cyclic complex of the $L$-projective leafwise smoothing operators on $(M, \maF)$. Our construction allows to couple the $K$-theory analytic indices of $L$-projective leafwise elliptic operators with the twisted cohomology of $B\maG$ producing scalar higher  invariants. Finally by adapting the Bismut-Quillen superconnection approach, we compute these higher twisted indices as integrals over the ambiant manifold of the expected twisted characteristic classes.

\end{abstract}

\maketitle
%
%
\section*{Introduction}

The Atiyah-Singer families index theorem has been  extended by Mathai, Melrose and Singer in \cite{MMS1, MMS2} to the class of projective families of operators. Given a smooth fibration $\pi: M\to B$ of closed manifolds, the authors of \cite{MMS1} consider local families of pseudodifferential operators over  open sets of a trivializing cover of $B$ such that the usual compatibility condition on triple overlaps fails by a scalar factor. This yields an integral  $3$-cohomology class of $B$ which is known to correspond to a Dixmier-Douady (DD) class of a  gerbe. A prototype example of such family arises from the attempt to define a Dirac operator along the fibers of some non $K$-oriented fibrations. Main geometric examples are though given by usual Dirac-type operators along the fibers with coefficients in  auxiliary (horizontally) twisted vector bundles over the total manifold $M$. See also \cite{NistorTroitsky}, where a closely related index theorem was  proved as well. The proof given in \cite{MMS1} of the projective index theorem  is the generalization to the twisted $K$-theory of the classical constructions of Atiyah-Singer for the families index theorem and they assumed the gerbe on the base manifold to be torsion. Later on, in \cite{MMS2}, the projective index theorem was extended to deal with more general gerbes, called decomposable. In \cite{BenaGoro}, the Bismut superconnection approach to the families index theorem \cite{Bismut, Quillen} allowed to give a new proof of the index formula, for Dirac-type operators with coefficients in twisted  vector bundles, under the assumption that only the pull-back of the gerbe to $M$ is torsion.  \\

%
%
%

In this paper we extend the previous results and in particular those of \cite{BenaGoro} to projective operators along the leaves of foliations. In the non-twisted case, the index problem for elliptic operators along the leaves of a smooth foliation on a closed manifold, was solved in the seminal works of Connes \cite{ConnesIntegration, ConnesSurvey} and then of Connes-Skandalis \cite{ConnesSkandalis}. In \cite{ConnesSurvey} and then in \cite{ConnesSkandalis}, Connes and Connes-Skandalis proved the   index theorem for leafwise elliptic operators, an equality  in the $K$-theory of the holonomy groupoid $C^*$-algebra. When the foliation is a  fibration,  the Connes-Skandalis theorem reduces to the Atiyah-Singer theorem for families  \cite{AtiyahSinger4}, and the formula then lives in the $K$-group of the base manifold,  so its Chern character  is totally computable in terms of characteristic classes. For general foliations though, the  $K$-group of the holonomy groupoid $C^*$-algebra is  harder to compute, and many fundamental noncommutative geometry constructions were developed in order to provide  computable pairings with such $K$-groups.
The de Rham homology of the base manifold must now be replaced by the cyclic cohomology of some dense subalgebra of the groupoid $C^*$-algebra. An interesting and rich enough class of cyclic cocycles arises from the homology of the classifying space of the holonomy  groupoid  \cite{ConnesBook, Gorokhovsky, GorokhovskyLott}. Recall for instance the breakthrough results of \cite{ConnesTransverse} which insure the extension of the pairing with the so-called fundamental class to the $C^*$-algebra $K$-theory, with deep topological consequences (see again \cite{ConnesTransverse}). The  index formulae obtained using the full homology of the classifying space $B\cG$ cannot in general  be deduced from the Connes-Skandalis theorem since the pairings are  only defined at the level of smooth dense subalgebras \cite{Gorokhovsky, GorokhovskyLott, Carrillo}, and most of the higher index formulae are then tackled independently by using the superconnection formalism of Quillen \cite{Quillen} and by extending to foliations the classical  approach to the families index problem \cite{Bismut, BenameurHeitsch2}.\\

In the present paper, we  use this superconnection approach in the context of foliations where the Dirac-type operator along the leaves   acts on a  vector bundle, which is twisted by  a gerbe on the (reduced) holonomy groupoid.  This groupoid plays in noncommutative geometry the role of the singular space of leaves. We also introduce many constructions which are valid for all \'etale Hausdorff Lie groupoids $\cG$ and for all $\cG$-equivariant smooth submersions $\pi: P\rightarrow B$.
Roughly speaking, a gerbe $L$ on  $\cG$ or a $\cG$-gerbe, is provided by a triple $(\maL, \varsigma, \nu)$ where $\maL$ is a gerbe  on the unit manifold $\cG^{(0)}$, $\varsigma$ is a $\cG$-action, that is a line bundle over the arrows $\cG^{(1)}$, with the natural compatibility condition in terms of a $2$-morphisms $\nu$ over the $2$-arrows, see Definition \ref{def:gerbe}. Our definition then agrees with Morita equivalence and with the approach in  \cite{TuXuGe}.
The space of compactly supported smooth sections of $\varsigma$ over $\cG^{(1)}$ is then endowed with the structure of an involutive convolution algebra by using the $2$-morphism $\nu$. If we complete it with respect  to the regular representations, we  obtain the regular $L$-twisted $C^*$-algebra $C^*_r (\cG, L)$ associated with the groupoid $\cG$ and the $\cG$-gerbe $L$. In the case of foliations, this latter $C^*$-algebra is the $L$-twisted $C^*$-algebra of the foliation $(M, \maF)$ and it is the natural receptacle for indices of $L$-projective leafwise elliptic pseudodifferential operators, exactly as in the untwisted case \cite{ConnesIntegration}. On the other hand, the DD-class associated with any $\cG$-gerbe $L$ lives in the integral $3$-cohomology $H^3(\cG, \Z)$ associated with the Lie groupoid $\cG$, or equivalently in the $3$-cohomology  of its classifying space $B\cG$.  We give a  construction of the Chern-Weil representative of this DD-class  in the Dupont-Sullivan complex of compatible forms. The complex of the compatible forms was first introduced by D. Sullivan in \cite{SullivanInfinitesimal} and applied to the study of classifying spaces and characteristic classes by J. Dupont, cf.  \cite{dupont}.  When pulled-back to the lifted groupoid $P\rtimes \cG$ with units $P$, we get a $P\rtimes \cG$-gerbe,  whose DD-class is the pull-back $\pi^*[L]$ of $[L]$. In the case of foliations, the $\maG$-submersion is free and proper and the pull-back DD-class $\pi^*[L]$ in $H^3 (P\rtimes \cG, \Z)\simeq H^3 (M, \Z)$ is  also the   pull-back of the DD-class on $B\maG$ under the classifying map of our foliation. So we naturally ask that this class is torsion. This corresponds for fibrations to asking that the DD-class pulled-back to the total manifold of the fibration is  torsion and we recover the usual assumption, see  \cite{MMS1, MMS2, BenaGoro}.  Now for the free and proper $\maG$-equivariant submersion $\pi:P\to B=\maG^{(0)}$, $L$-twisted vector bundles on $M=P/\maG$,  are the $\pi^*\maL$-twisted  bundles over $P$ which are $\cG$-equivariant, see Definition \ref{TwistedBundles}.  Given such bundle $\maE\to P$, the algebra $\Psi_{\pi^*\maL} (P\vert \cG^{(0)}, \maE)^\cG$ of $L$-projective leafwise pseudodifferential operators  is  defined as the algebra of (properly supported) fiberwise $\pi^*\maL$-projective pseudodifferential operators which are $\cG$-equivariant. These are  the $L$-projective leafwise pseudodifferential operators on the foliation $\maF$ of $M$ obtained by moding out the fibers of $\pi$ and it can suggestively also be denoted $\Psi_{L} (M, \maF; \maE)$.
 Its two-sided  ideal  $\Psi^{-\infty}_{L} (M, \maF; \maE)$ of $L$-projective leafwise smoothing operators  is then  $\Psi_{\pi^*\maL}^{-\infty} (P\vert \cG^{(0)}, \maE)^\cG$.


Leafwise ellipticity of an element $A$  then insures the existence of a $K$-theory analytic index class $\ind_a (A)$ which lives in the $K$-theory group $K \left(\Psi_{L}^{-\infty} (M, \maF;  \maE)\right)$ and hence a $C^*$-index $\Ind_a (A)$ which lives in the $K$-theory of its $C^*$-completion , i.e.
$$
\ind_a (A) \; \in \; K \left(\Psi_{L}^{-\infty} (M, \maF; \maE)\right) \text{ and } \Ind_a(A) \; \in  \; K \left({\overline{\Psi_{L}^{-\infty} (M, \maF; \maE))}}\right).
$$
As for the non twisted case, the index class yields an analytic class in the $K$-theory of the gerbe $C^*$-algebra $C^*_r (\cG, L)$ which is compatible with the Morita equivalence
$$
C^*_r (\cG, L) \sim {\overline{\Psi_{L}^{-\infty} (M, \maF; \maE)}}
$$
 If we denote by $\maA$ the $\cG$-equivariant Azumaya $C^*$-algebra bundle  associated with $\maE$ then  the principal symbol of any $L$-twisted leafwise elliptic operator defines a class
$$
[\sigma (A)] \; \in \;
 K \left(C_0(T^*P\vert \cG^{(0)}, \maA)^\cG\right),
$$
and we obtain in this way a well defined $L$-projective leafwise index map
$$
\Ind_a: K \left(C_0(T^*P\vert \cG^{(0)}, \maA)^\cG\right)
\longrightarrow K\left( C^*_r (\cG, L)\right).
$$

The goal here  is to define and compute  the higher  indices of $L$-projective leafwise elliptic operators associated with the pairings against the (periodized) twisted  homology classes of the classifying space $B\cG$. Our approach uses superconnections and yields directly to a topological formula in terms of twisted characteristic classes. More precisely, we define the pairing of the smooth index class $\ind_a (A)$ with cyclic cocycles arising from the twisted Dupont-Sullivan complex of smooth currents $(\maC^\bullet_\Delta(\cG)[(\u)^{-1}, \u], d^t_{\varOmega(L)})$, where $\u$ is  a formal variable of degree $+2$. These currents correspond to differential forms which incorporate an orientation so that they can be integrated even on nonoriented manifolds.  We achieve this pairing by
constructing  a morphism $\Lambda^\cE$ from this Dupont-Sullivan complex to the (periodic) cyclic cohomology of the ideal $\Psi_{\pi^*\maL}^{-\infty} (P\vert \maG^{(0)}, \maE)^\cG$, representing the morphism
\begin{equation*}
\Lambda^{\cE} \colon H_{\varOmega(L)} (B\cG) \longrightarrow \HC^\bullet \left(\Psi_{L}^{-\infty} (M, \maF; \maE)\right).
\end{equation*}

When the $L$-projective bundle $\cE$ is $\Z_2$-graded by the involution $\gamma$, we obtain in the same way the corresponding morphism $\Lambda^{\cE, \gamma}$.
We thus obtain the  geometric construction of cyclic cocycles we use, and of higher indices of any $L$-projective leafwise elliptic operator $D$ acting from $\cE^+$ to $\cE^-$ associated with a given even degree cocycle $c$ in the Dupont-Sullivan complex. The precise formula is:
\begin{equation*}
\ind_c (D) := \langle (u^*)^{-\frac{\deg c}{2}}\Lambda^{\cE, \gamma} (c), \ind_a D \rangle.
\end{equation*}
The pairing in the RHS is the usual normalized pairing between $K$-theory and cyclic cohomology.

\medskip

Our morphism $\Lambda^\cE$ relies on a twisted and simplicial extension of the so-called   Connes $\Phi$-map  (see Section 2.$\delta$ in \cite{ConnesBook}).

\medskip

{\bf{Theorem \ref{PhiE}}} [Twisted Connes' $\Phi$-map]
{\em {There exists
a morphism of complexes
\begin{equation*}
\Phi^\cE \colon (\maC^\bullet_\Delta(\maG)[(\u)^{-1}, \u], d^t_{\varOmega(L)}) \longrightarrow C^\bullet(\maG, \CC^\bullet (\Psi^{-\infty}_{\pi^*\maL}(P|\maG^{(0)}, \maE))),
\end{equation*}
which extends to the twisted case, the classical  Connes $\Phi$-map.}}

\medskip

Recall that the classical map $\Phi$ has been introduced in \cite{ConnesTransverse, ConnesBook} as a map from the cohomology of $B\maG$ with the coefficients in the orientation bundle to the cyclic cohomology of the smooth convolution algebra $C^\infty_c(\maG)$. This map is an injection, cf. \cite{ConnesBook}, Section 2.$\delta$. In \cite{bryni} Brylinski-Nistor computed the cohomology of $C^\infty_c(\maG)$ for a general Hausdorff \'etale $\maG$. This computation has been extended to the non-Hausdorff case by Crainic in \cite{crainic}.  Our theorem above extends Connes' construction to the twisted case.


Then using  constructions again from \cite{bryni, crainic} and the isomorphism
\[
 \HC^\bullet(\Psi^{-\infty}_{\pi^*\maL}(P|\maG^{(0)}, \maE) \rtimes \maG )
 \longrightarrow  \HC^\bullet \left(\Psi^{-\infty}_{\pi^*\maL}(P|\maG^{(0)}, \maE)^\maG\right)\,
\]
the allowed morphism $\Lambda^{\cE}$  is deduced (see Equation \ref{eq:Lamb}).

\medskip

As explained above, we also include the computation of the higher index $\ind_c (D)$ by a topological formula {\em{\`a la Atiyah-Singer}}. We concentrate on  Dirac-type operators along the leaves twisted by an $L$-projective hermitian bundle.
The topological formula that we obtain is, as expected, some pairing of the class of $c$ with the twisted Atiyah-Singer  characteristic class. In order to pair the class of $c$ with cohomology classes on $M=P/\cG$, we pull them back to the groupoid $\cP:=P\rtimes \cG$ and use the Dupont-Sullivan complex of forms on $\cP$ and the well-defined cap product
$$
(\Omega_\bullet^\Delta(\cP), d_{\varOmega(L)}) \otimes (\maC^\bullet_\Delta(\cP)[(\u)^{-1}, \u], d^t_{\varOmega(L)}) \longrightarrow (\maC^\bullet_\Delta(\cP)[(\u)^{-1}, \u], d).
$$
Recall that  the pull-back of $L$ to $\cP$ is torsion so that the Chern-Weil representative of its DD-class in the Dupont-Sullivan complex of $\cP$ is exact. Since the Lie groupoid $\cP$ is proper, free and cocompact, the twisted (and non-twisted) cohomology of $\cP$ is isomorphic to the cohomology of the closed manifold $M$. More precisely, the $\hat{A}$-form of the tangent bundle to the foliation $\maF$  is easily represented by a closed form ${\hat{A}} (P\vert \cG^{(0)})$ in the non-twisted Dupont-Sullivan complex, while the Chern character of the $L$-projective bundle $\maE$, with respect to the spinor bundle $\maS$ as in the notations of \cite{bgv},  is naturally represented by a closed form $\Ch_L(\maE/S)$ in the $\pi^*\varOmega (L)$-twisted Dupont-Sullivan complex of forms.  Hence, we end up with a $d_{\pi^*\varOmega (L)}$-closed form $ {\hat{A}} (P\vert \cG^{(0)})\Ch_L(\maE/S)$ which can be paired with the $d^t_{\pi^*\varOmega (L)}$-closed smooth current $\pi^*c$ to yield a  smooth closed current (for the non-twisted differential)
$$
(\pi^*c)\cap \left(\widehat{A}(P|\cG^{(0)}) \Ch_L(\maE/S)\right) \; \in \; \maC^\bullet_\Delta(\cP)[(\u)^{-1}, \u].
$$
Simplicial integration and integration over $M=P/\cG$ allows to integrate the closed smooth current $(\pi^*c)\cap \left(\widehat{A}(P|\cG^{(0)}) \Ch_L(\maE/S)\right)$ and extract eventually a well defined complex number:
$$
\left\langle \pi^*c\; , \; \widehat{A}(P|\cG^{(0)}) \Ch_L(\maE/S)\right\rangle.
$$
Our main result reads:

\medskip

{\bf{Theorem  \ref{thm:main}}. }  {\em{Let $c$ be a cocycle of even degree in $(\maC^\bullet_\Delta(\cG)[(\u)^{-1}, \u], d^t_{\varOmega(L)})$.
\begin{equation*}
\ind_c (D) = \left\langle  \pi^*c\; , \;  \widehat{A}(P|\cG^{(0)}) \Ch_L(\maE/S)\right\rangle.
\end{equation*}}}

Notice that  the cocycle $c$ represents a class $[c]$ in $H_{\varOmega(L)} (B\maG)$ and if $f: M\to B\maG$ is a classifying map for our foliation $\maF$ on $M$, then $f^*[c]=[\pi^*c]$. An explicit example illustrating the geometric set up underlying Theorem \ref{thm:main} is given in Subsection \ref{sub:ex}.

Our theorem  unifies several previous results.
First of all in the non-twisted case, we recover the results of \cite{GorokhovskyLott, ConnesBook} as well as the results of \cite{BenameurHeitsch1}. In the case of a fibration, the  gerbe is defined on the base of the fibration and we recover the results of \cite{MMS1, MMS2, BenaGoro}. In the case of a covering manifold, this result has been proved by Mathai-Marcolli (\cite{MM1}, \cite{MM2}), extending Connes-Moscovici higher index theorem \cite{ConnesMoscovici90} to the case of projective actions.

Let us recall  that the twisted version of the $K$-theoretic Connes-Skandalis index theorem has been recently investigated by Carrillo Rouse and Wang (\cite{CaW}) by using the Connes tangent groupoid approach. More precisely, \cite{CaW} used a geometric set up \`a la Hilsum-Skandalis for the definition of twisted groupoid which is well adapted to the proof of twisted $K-$theoretic index theorems.
Their geometric set up is equivalent to ours. Their focus however is on the $K$-theoretic results.  The focus of  this paper, on the other hand,  is on cohomological index theorems in suitable cyclic cohomological groups and on explicit cohomological index formula.

Finally, we mention that Carey and Wang (\cite{CW}) have proven a Riemann-Roch theorem in twisted $K-$theory and used it to prove a cohomological formula for the index pairing in twisted $K-$theory, so as to get applications to D-branes and D-charges. See  also the paper
\cite{CMW} for  definition and properties of differential twisted $K-$theory.
 \\

{\em{Acknowledgements.}} We benefited from discussions with many colleagues, in particular during the NCG workshop at MFO and we thank them all. We are in particular indebted to the following people: A. Carey, A. Connes, T. Fack, J. Heitsch, V. Mathai, V. Nistor, G. Skandalis and J.-L. Tu. MB and EL wish to thank the french National Research Agency for support via the project ANR-14-CE25-0012-01 (SINGSTAR). AG wishes to thank the National Science Foundation (NSF)  for support via the  grant DMS-0900968 \\

\bigskip

\tableofcontents

\section{Homologies of \'Etale  Lie groupoids}\label{AppendixHomology}



Let $\maG$ be an \'etale Hausdorff locally compact groupoid with units $\maG^{(0)}$ and denote for $k\geq 1$ by $\maG^{(k)}$ the set of composable $k$-tuples, i.e.
$$
\maG^{(k)}= \{(\gamma_1, \cdots, \gamma_k)\in \maG^k\text{ such that } r(\gamma_{i+1})=s(\gamma_i), i=1, \cdots, k-1\}
$$
 We define   for $l \ge1$ and $0\leq i_1<\cdots <i_{l+1}\leq k$, the maps $\pr^k_{i_1i_2\ldots i_{l+1}} : \maG^{(k)} \rightarrow \maG^{(l)}$ by
\begin{equation}\label{prk}
\pr^k_{i_1i_2\ldots i_{l+1}} (\gamma_1, \gamma_2, \ldots, \gamma_k): =
(\gamma_{i_1+1}\ldots \gamma_{i_2}, \gamma_{i_2+1}\ldots \gamma_{i_3}, \ldots , \gamma_{i_l+1}\ldots \gamma_{i_{l+1}}).
\end{equation}
For $l=1$ for instance, we have
\begin{equation*}
\pr_{ij}^k (\gamma_1, \cdots, \gamma_k) := \gamma_{i+1} \cdots \gamma_j, \quad  \text{ for } 0\leq i <j\leq k.
\end{equation*}
For $l=0$ we define the maps $\pr^k_{i} \colon \maG^{(k)} \to \maG^{(0)}$ by
\begin{equation*}
\pr^k_{i} ( \gamma_1, \gamma_2, \ldots, \gamma_k) : =
s(\gamma_i) = r(\gamma_{i+1}).
\end{equation*}

 By $\maG^{(\bullet)}$ we denote the simplicial space defined as follows. The face maps are given for $k\geq 2$ by
\begin{equation}\label{face}
\del_i(\gamma_1, \ldots , \gamma_k)=
\begin{cases}
(\gamma_2, \ldots , \gamma_k) &\text{ for } i=0\\
(\gamma_1, \ldots ,\gamma_i \gamma_{i+1}, \ldots, \gamma_k) &\text{ for } 0<i<k\\
(\gamma_1, \ldots , \gamma_{k-1}) &\text{ for } i=k
\end{cases}
\end{equation}
and for $k=1$, we define $\del_0(\gamma_1)= r(\gamma_1)$ and $\del_1 (\gamma_1)= s(\gamma_1)$.

We denote as well by $r$ the map
$$
r=\pr^k_0:\cG^{(k)} \longrightarrow \cG^{(0)} \text{  given by } \pr^k_0 (\gamma_1, \cdots, \gamma_k) = r(\gamma_1).
$$
We shall only be interested in (Hausdorff \'etale) Lie groupoids which  we assume from now on, hence all the groupoid operations are smooth.

\subsection{The de Rham bicomplex}

We recall now the de Rham bi-complex $\Omega^\bullet (\maG^{(\bullet)}, \R)$ associated with the nerve of the groupoid $\maG$ and which computes the cohomology (with real coefficients)  of its classifying space $B\maG$.  So,
$$
\Omega^m (\maG^{(k)}, \R) := C^\infty_c (\maG^{(k)}, r^*\Lambda^m(T^*\cG^{(0)})),
$$
and if $d$ is the de Rham differential on the manifold $\cG^{(0)}$ then it lifts to the differential $r^*d$ which will  still be denoted $d$ for simplicity. So, the differentials are
$$
d: \Omega^{\bullet} (\maG^{(\bullet)}, \R) \rightarrow \Omega^{\bullet+1} (\maG^{(\bullet)}, \R) \text{ and } \del: \Omega^{\bullet} (\maG^{(\bullet)}, \R) \rightarrow \Omega^{\bullet} (\maG^{(\bullet+1)}, \R)
$$
where  $\del$ is induced by the simplicial maps $\del_i \colon \maG^{(k)}\to \maG^{(k-1)}$ given  in \eqref{face}. More precisely,
$$
\del (\omega) = (-1)^m  \sum (-1)^i \del_i^*\omega \text{ for }\omega\in \Omega^m(\cG^{(\bullet)}).
$$
Notice that $\del_1^*\omega$ uses the bundle transformations $\gamma:\Lambda^\bullet (T^*_{s(\gamma)} \cG^{(0)}) \to \Lambda^\bullet (T^*_{r(\gamma)}\cG^{(0)} )$ which are induced by the corresponding usual bundle maps $T_{r(\gamma)} \cG^{(0)} \to T_{s(\gamma)} \cG^{(0)}$, these latter being well defined for any $\gamma$ since $\maG$ is \'etale.

\begin{lemma}\label{proper}
Given a Hausdorff \'etale groupoid ${\what{\cG}}$ which is proper, the  cohomology of the complex  $(\Omega^{m} (\what{\maG}^{(\bullet)}), \partial)$   is trivial in positive degrees.
\end{lemma}

\begin{proof}
Since ${\what{\cG}}$ is proper and Hausdorff there exists $h \in C^\infty_c(\what{\cG}^{(0)})$ such that
\begin{equation*}\forall z \in \what{\cG}^{(0)} ,\; \sum \limits_{\{g\in \what{\cG}\ |\ s(g)=z\}} h(zg)=1.
\end{equation*}
Recall that we have for $m\geq 1$ the isomorphisms
$$
\gamma:\Lambda^m T_{s(\gamma)}^*\what{\cG} \longrightarrow \Lambda^m T_{r(\gamma)}^*\what{\cG},
$$
with the usual functorial properties.
Define then  the operator $K \colon \Omega^{m} (\what{\maG}^{(\bullet)}) \to \Omega^{m} (\what{\maG}^{(\bullet-1)})$ by
\begin{equation*}
K\omega (\gamma_1, \cdots, \gamma_k) := \sum_{s(\gamma)=r(\gamma_1)} h(r(\gamma)) \times \gamma^{-1} \omega (\gamma, \gamma_1,\cdots, \gamma_k).
\end{equation*}
A straightforward computation gives
\begin{equation*}
K\circ \del +\del \circ K = \id \colon \Omega^{m} (\what{\maG}^{(i)}) \to \Omega^{m} (\what{\maG}^{(i)}), \text{ for } i\ge 1.
\end{equation*}
\end{proof}

We shall be interested in the following differentiable situation of a  $\maG$-action on a smooth manifold $P$.
A $\maG$ equivariant submersion is given by a submersion
$\pi \colon P \to \maG^{(0)}$ and a   diffeomorphism
$\lambda_\gamma \colon P_{r(\gamma)} \to P_{s(\gamma)}$ for every  $\gamma \in \maG$ such that $\lambda_{\gamma_1 \gamma_2} =\lambda_{\gamma_2 }\circ \lambda_{ \gamma_1}$ and $\lambda_{\gamma^{-1}} = \lambda_{\gamma}^{-1}$. We set $\lambda_\gamma (p) := p\gamma$ for $p\in P_{r(\gamma)}$.
Given such an action we can form a new \'etale groupoid $\cP = P \rtimes \maG$ with  the space of objects $P$ and the morphisms given by pairs
$(p, \gamma) \in P \times \cG$ such that $\pi(p) = r(\gamma)$ with $r(p, \gamma) = p$, $ s(p, \gamma) = p\gamma$. The map $\pi$ defines a morphism of groupoids
$(p, \gamma) \mapsto \gamma$. The groupoid $\cP$ is proper if and only if the $\maG$-action on $P$ is proper.

The associated simplicial space is described explicitly as follows:
$$
\cP^{(k)}= (P \rtimes \cG)^{(k)}=\{(p, \gamma_1, \ldots , \gamma_k) \in P \times \cG^{(k)}
\ |\  \pi(p)= r(\gamma_{1})\}.
$$
The face maps $\del_i\colon \cP^{(k)} \to \cP^{(k-1)}$ are given explicitly by
\begin{equation*}
\del_i(p, \gamma_1, \ldots , \gamma_k)=
\begin{cases}
(p\gamma_1, \gamma_2, \ldots , \gamma_k) &\text{ for } i=0\\
(p, \gamma_1, \ldots ,\gamma_i \gamma_{i+1}, \ldots, \gamma_k) &\text{ for } 0<i<k\\
(p, \gamma_1, \ldots , \gamma_{k-1}) &\text{ for } i=k
\end{cases}
\end{equation*}

$\pi$ induces  submersions $ \pi^{(n)} \colon \cP^{(n)} \to \cG^{(n)}$ by
\begin{equation*}
\pi^{(n)}(p, \gamma_1, \ldots , \gamma_k)=(\gamma_1, \ldots , \gamma_k)
\end{equation*}
These submersions are compatible  with the face maps.

\subsection{The Dupont-Sullivan bicomplex}\label{defcurrent}
We shall also use the Dupont-Sullivan bicomplex associated with $\maG$. Denote by $\Delta^k$ the standard $k$-simplex
\begin{equation*}
\Delta^k := \{(t_0, \cdots , t_{k})\in \R^{k+1}\ | \  0\leq t_i,  \sum_{i=0}^k t_i=1\}.
\end{equation*}
and by $\delta_i \colon \Delta^{k-1} \to \Delta^k$ the face maps
\begin{equation*}
\delta_i\ (t_0, \cdots , t_{k-1}) = (t_0, \cdots,t_{i-1}, 0, t_i, \cdots , t_{k-1})
\end{equation*}

 The algebra of Dupont-Sullivan forms was introduced in \cite{dupont}  and will be  denoted here by $\Omega^*_{DS}(\cG)$. So $\Omega^*_{DS}(\cG) \subset \prod \limits_k \Omega^*(\cG^{(k)}) \otimes \Omega^*(\Delta^k)$ and an element of $\Omega^*_{DS}(\cG)$
is given by a collection $\omega=(\omega_k)_k$ of forms
\begin{equation*}\omega_k \in  \Omega^*(\maG^{(k)})  \otimes \Omega^*(\Delta^k) =\Omega^*(\cG^{(k)} \times \Delta^k)
\end{equation*}
satisfying for every $k\geq 1$ the following relation
\begin{equation*}
(\id \times \delta_i)^* \omega_{k} = (\del_i \times \id)^* \omega_{k-1} \in \Omega^*(\cG^{(k)} \times \Delta^{k-1}).
\end{equation*}
There is a natural bigrading on $\Omega^*_{DS}(\cG)$ given by
\[\Omega^{i,j}_{DS}(\cG) =  \Omega^*_{DS}(\cG)\bigcap \left(\prod \limits_k \Omega^i(\cG^{(k)}) \otimes \Omega^j(\Delta^k)\right),
\]
and we set as usual $\Omega^m_{DS} (\cG) \,  := \, \bigoplus_{i+j=m} \Omega^{i,j}_{DS}(\cG).$
The differential is given by $d_\maG+d_\Delta$ where $d_\maG$ is de Rham differential on $\maG$
and $d_\Delta$ is de Rham differential on $\Delta$.
We will also consider the complex $\Omega_\bullet^\Delta(\cG)= \Omega^*_{DS}(\cG)[u]$
where the differential has now degree $-1$ and is given  by   $u(d_\maG+d_\Delta)$.

The dual version of currents will also be needed.
We refer the reader  to \ref{Currents} for the definition of the spaces  $\maC^\bullet(B):=\maC_*(B)[\u]$ and $\maC_*(B)[(\u)^{-1}, \u]:=\maC^\bullet(B)[(\u)^{-1}]$, where $\u$ is a formal variable of degree $+2$. For an \'etale   map $f\colon B \to B'$ there is a well-defined pull-back
\begin{equation*}
f^* \colon \maC^\bullet(B') \to \maC^\bullet(B),
\end{equation*}
dual to the ``integration'' along the fibers. We can therefore define the space of simplicial smooth  currents
$\maC_\Delta^\bullet(\cG) \subset \prod \limits_k \maC_*(\cG^{(k)})[(\u)^{-1}, \u] \otimes \Omega^*(\Delta^k)$
as given by a collection $\{\xi_k\}$ where
\begin{equation*}\xi_k \in  \maC_*(\maG^{(k)})[(\u)^{-1}, \u] \otimes \Omega^*(\Delta^k)
\end{equation*}
satisfying for every $k$
\begin{equation*}
(\id \otimes \delta_i^* )\xi_{k} = (\del_i^* \otimes \id) \xi_{k-1} \in \maC_*(\cG^{(k)})[(\u)^{-1}, \u] \otimes \Omega(\Delta^{k-1}).
\end{equation*}
The differential is given by $d^t=\u d^t_\maG+d_\Delta$.



Set for  $\alpha \otimes \beta \in \Omega^i(\cG^{(k)})[u] \otimes \Omega^j(\Delta^k)$, $\xi \otimes \omega \in \maC_*(\cG^{(k)})[ \u] \otimes \Omega^*(\Delta^k)$ (cf. \ref{Currents}).
\begin{equation*}
\langle (\xi \otimes \omega) ,(\alpha \otimes \beta)\rangle = (-1)^{\deg \alpha \deg \omega} \langle (\u)^{-j}\xi,  \alpha \rangle  \int_{\Delta^k} \omega \wedge \beta,
\end{equation*}
We then have, with $\phi \in \Omega^*(\cG^{(k)})[u] \otimes \Omega^*(\Delta^k)$,  $c \in \maC_*(\cG^{(k)})[ \u] \otimes \Omega^*(\Delta^k)$,
\begin{equation*}
\langle d_\cG^t c,  \phi  \rangle + (-1)^{\deg c} \langle  c,  d_\cG\phi  \rangle =0
\end{equation*}
and
\begin{equation*}
\langle d_\Delta  c,  \phi  \rangle + (-1)^{\deg c} \langle  c,  d_\Delta\phi  \rangle =\sum_{i=0}^k (-1)^i \langle (\id \otimes \delta_i^* )c, (\id \otimes \delta_i^* ) \phi\rangle
\end{equation*}
 Define also a cap product $\maC^\bullet_\Delta(\cG) \otimes \Omega_\bullet^\Delta(\cG)   \to  \maC^\bullet_\Delta(\cG) $ by setting
\begin{equation*}
(\xi \otimes \omega) \cap (\alpha \otimes \beta) = (-1)^{\deg \alpha \deg \omega} (\u)^{-j}(\xi \cap \alpha) \otimes (\omega \wedge \beta)
\end{equation*}

If $\varOmega \in \Omega^3_{DS}(\cG)$ is a closed form we define the associated twisted complexes which compute the twisted cohomology/homology of the classifying space $B\cG$. So,  $(\Omega_\bullet^\Delta(\cG), d_\varOmega)$ is obtained by considering the differential $d_\varOmega= u(d_\maG+d_\Delta) +u^2 \varOmega \wedge \cdot$, where the wedge is induced from the wedge products on the manifolds $\cG^{(k)} \times \Delta^k$. Similarly, we define the twisted complex $(\maC^\bullet_\Delta(\cG), d^t_\varOmega )$ by setting for $c \in \maC^\bullet_\Delta(\cG)$
\begin{equation} \label{defdtsimp}
d^t_\varOmega c  = d^t c + (-1)^{\deg c} c \cap  \varOmega \, = \, \u d^t_\maG c+d_\Delta c +(-1)^{\deg c} c \cap  \varOmega.
\end{equation}
The homology of this complex will be denoted $H^\bullet_\Omega (B\cG)$.   We point out that only $3$ forms $\varOmega$ with no $(3,0)$ component arise in the present paper.

\medskip

\section{Gerbes and twisted bundles}\label{preliminaries}

We review in this section some background material about connections on gerbes and their morphisms. We also recall many results that will be used in the next sections. Our companion paper is \cite{BenaGoro} as we shall use many constructions from there. We give here only the general overview, referring the reader to \cite{Brylinski}, \cite{Hitchin},  \cite{Murray} and \cite{MS}
for the details. The differential geometry of not necessarily abelian gerbes is described in \cite{BM}.

\subsection{Gerbes and connections} \label{subs:curve} $\;$

We shall describe the gerbes in terms of their descent data.

Let $M$ be a smooth manifold. Given an open cover $(U_\alpha)_{\alpha\in \Lambda}$ of $M$, we set as usual
$$
U_{\alpha_1\cdots \alpha_k} = \bigcap_{1\leq j \leq k} U_{\alpha_j}.
$$

\begin{definition}\label{def:1-m}
A descent datum for a gerbe $\maL$ on $M$   is  the collection $(U_\alpha, \maL_{\alpha\beta}, \mu_{\alpha\beta\gamma})$ where $(U_\alpha)_{\alpha\in \Lambda}$ is an open cover of $M$, $(\maL_{\alpha\beta} \to U_{\alpha\beta})_{\alpha, \beta\in \Lambda}$ is a collection of line bundles and $\mu_{\alpha\beta\gamma}: \maL_{\alpha\beta} \otimes\maL_{\beta\gamma} \to \maL_{\alpha\gamma}$ is a collection of  line bundle isomorphisms over each triple intersection $U_{\alpha\beta\gamma}$ such that over each quadruple intersection $U_{\alpha\beta\gamma\delta}$, the following diagram commutes

\begin{equation*}
\begin{CD}
\maL_{\alpha\beta}\otimes \maL_{\beta\gamma}\otimes\maL_{\gamma\delta}
@>\mu_{\alpha\beta\gamma}\otimes \id>>
\maL_{\alpha\gamma}\otimes\maL_{\gamma\delta}\\
@V \id\otimes\mu_{\beta\gamma\delta}VV  @VV\mu_{\alpha\gamma\delta}V\\
\maL_{\alpha\beta}\otimes \maL_{\beta\delta} @>\ \ \mu_{\alpha\beta\delta}\ \ >>\maL_{\alpha\delta}
\end{CD}
\end{equation*}
\end{definition}

Notice that we don't need to assume in this definition  that the open sets $(U_\alpha)_{\alpha\in \Lambda}$
are contractible, however, our covers will always be good covers \cite{BottTu}.

Given two  descent data $ (U_\alpha, \maL_{\alpha\beta}, \mu_{\alpha\beta\gamma})$ and
$ (U_\alpha, \maL_{\alpha\beta}', \mu_{\alpha\beta\gamma}')$ on the same open cover $\{U_{\alpha}\}$,
an isomorphism between them is given by line bundles $S_{\alpha}$ on $U_{\alpha}$ and isomorphisms of line
bundles $\lambda_{\alpha \beta} \colon S_{\alpha}^{-1} \otimes \maL_{\alpha \beta} \otimes S_{\beta} \to \maL_{\alpha \beta}'$ over $U_{\alpha \beta}$ so that the following diagram commutes

\begin{equation*}
\begin{CD}
S_{\alpha}^{-1}\otimes \maL_{\alpha \beta}
\otimes S_{\beta}\otimes S_{\beta}^{-1}\otimes \maL_{\beta \gamma} \otimes S_{\gamma}
@>\id \otimes \mu_{\alpha \beta \gamma} \otimes \id>>
S_{\alpha}^{-1}\otimes\maL_{\alpha \gamma} \otimes S_{\gamma} \\
@V\lambda_{\alpha \beta}\otimes\lambda_{\beta \gamma}VV @VV\lambda_{\alpha \gamma}V \\
\maL_{\alpha \beta}' \otimes \maL_{\beta \gamma}' @>\ \ \ \mu_{\alpha \beta \gamma}'\ \ \ >>\maL_{\alpha \gamma}'
\end{CD}
\end{equation*}

It is clear how to compose such (iso)morphisms.

\begin{definition} \label{def:2-m}
  Given two
 isomorphisms     $(S_{\alpha}, \lambda_{\alpha \beta})$ and $(S_{\alpha}', \lambda_{\alpha \beta}')$
between  $ (U_\alpha, \maL_{\alpha\beta}, \mu_{\alpha\beta\gamma})$ and
$ (U_\alpha, \maL_{\alpha\beta}',  \mu_{\alpha\beta\gamma}')  $, a $2$-morphism between them is a collection
of line bundle isomorphisms $\nu_{\alpha} \colon S_{\alpha} \to S_{\alpha}'$ such that
\begin{equation}\label{2morphism}
\lambda_{\alpha \beta}' \circ (\nu_{\alpha}^{-1} \otimes \id \otimes \nu_{\beta}) =\lambda_{\alpha \beta}
\end{equation}
where we denote by
$\nu_{\alpha}^{-1}$ the isomorphism $S_{\alpha}^{-1} \to (S_{\alpha}')^{-1}$ induced by $\nu_{\alpha}$.
\end{definition}

Notice that Equation \ref{2morphism} allows to define the collection $(\lambda_{\alpha\beta})$ which then satisfies the commutativity of the previous diagram.
More generally, a $2$-morphism can be defined between any isomorphisms of gerbes as follows.

\begin{definition} \label{def:2-m-general}\
  Consider an
 isomorphism     $\varsigma_1=(S_{1,\alpha}, \lambda_{1, \alpha \beta})$ between  the descent data $ (U_\alpha, \maL_{1,\alpha\beta}, \mu_{1,\alpha\beta\gamma})$ and
$ (U_\alpha, \maL_{1, \alpha\beta}',  \mu_{1, \alpha\beta\gamma}')  $ and, an isomorphism $\varsigma_2=(S_{2, \alpha}', \lambda_{2, \alpha \beta}')$
between  the descent data $ (U_\alpha, \maL_{2,\alpha\beta}, \mu_{2,\alpha\beta\gamma})$ and
$ (U_\alpha, \maL_{2, \alpha\beta}',  \mu_{2, \alpha\beta\gamma}')  $. Then,
 a $2$-morphism between them is given by the following data:
 \begin{itemize}
 \item An isomorphism $\varsigma=(S_{\alpha}, \lambda_{ \alpha \beta})$ between  the descent data $ (U_\alpha, \maL_{1,\alpha\beta}, \mu_{1,\alpha\beta\gamma})$ and
$ (U_\alpha, \maL_{2, \alpha\beta},  \mu_{2, \alpha\beta\gamma})  $ and an isomorphism $\varsigma'=(S_{ \alpha}', \lambda_{ \alpha \beta}')$
between  the descent data $ (U_\alpha, \maL_{1,\alpha\beta}', \mu_{1,\alpha\beta\gamma}')$ and
$ (U_\alpha, \maL_{2, \alpha\beta}',  \mu_{2, \alpha\beta\gamma}')  $,  in the sense of Definition \ref{def:1-m} above;
\item A $2$-morphism $\nu$ between the  corresponding descent data for the composite isomorphism $\varsigma'\circ \varsigma_1$ and for the composite isomorphism $\varsigma_2\circ \varsigma$, in the sense of Definition \ref{def:2-m}.
\end{itemize}
\end{definition}

We shall refer to such $2$-morphism by $(\nu, \varsigma,\varsigma')$. Definition \ref{def:2-m} corresponds to the identity gerbes isomorphisms and in this case we refer to the $2$-morphism simply by $\nu$.

Let $(V_i, \varrho)_{i\in I}$ be a refinement of the open cover   $(U_\alpha)_{\alpha\in \Lambda}$ of $M$. So $(V_i)_{i\in I}$ is an open cover of $M$ with
$\varrho: I \rightarrow \Lambda$ such that $V_i \subset U_{\varrho(i)}.$
Then restriction to the refinement $(V_i, \varrho)_{i\in I}$ of $ (U_\alpha, \maL_{\alpha\beta}, \mu_{\alpha\beta\gamma})$ is the descent datum $\maL'=(V_i, \maL'_{ij}, \mu'_{ijk})$ given by:
$$
\maL'_{ij} := \maL_{\varrho(i)\varrho(j)}|_{V_{ij}} \text{ and } \mu'_{ijk} := \mu_{\varrho(i)\varrho(j)\varrho(k)}|_{V_{ijk}}.
$$
Similarly one defines restriction of the isomorphisms and $2$-morphisms to a refinement. We do not
distinguish between a descent datum and  its restriction to a refinement.
Thus  for instance, the isomorphism between two  descent data $ (U_\alpha, \maL_{\alpha\beta}, \mu_{\alpha\beta\gamma})$ and
$ (U_\alpha', \maL_{\alpha\beta}', \mu_{\alpha\beta\gamma}')$ is an isomorphism between their restrictions to some
common refinement of $\{ U_\alpha\}$ and $\{ U_\alpha'\}$. Equivalence classes of descent data are thus well defined, and an equivalence class of  (Dixmier-Douady) gerbes  on $M$ is an equivalence class of such descent data on $M$.

Recall more, a  gerbe can be defined as a maximal collection of descent data $\maD_i$, $i \in A$  together with the
isomorphisms $s_{ij}\colon \maD_j \to \maD_i$ for each $i, j \in A$ and $2$-morphisms $\nu_{ijk} \colon
s_{ij} s_{jk} \to s_{ik}$ satisfying the natural associativity condition. We refer the reader
to the book \cite{Brylinski} for the details.

If the cover $\{U_{\alpha}\}$ is good \cite{BottTu},  all the bundles $\maL_{\alpha \beta}$ are trivializable.
After choice of such a trivialization,
the collection $(\mu_{\alpha\beta\gamma})$ can be viewed  as a \v{C}ech $2$-cochain with coefficients in the sheaf $\underline{\C^*}$ of smooth functions with values in the nonzero complex numbers $\C^*$. We can and will, unless otherwise specified,  always work with good covers.
The compatibility  condition over $U_{\alpha\beta\gamma\delta}$ tells us that $\mu$ is a $2$-cocycle and hence defines a cohomology class  $[\mu]\in H^2(M; \underline{\C^*}) \cong H^3(M, \Z)$. This class is a well defined invariant of the gerbe called the Dixmier-Douady class. We denote this class by $[\maL]$.  Every class in $H^3(M, \Z)$ is a class of a gerbe defined by this class uniquely up to an isomorphism  (see again \cite{Brylinski}).

Given a smooth map $f: M'\to M$ between smooth manifolds $M'$ and $M$, we can pull-back any descent datum for a gerbe on $M$ to a descent datum on $M'$. The pull backs of isomorphic descent data are isomorphic and thus we obtain a well-defined pull-back of a gerbe. Clearly the Dixmier-Douady class of the pull-back is the pull-back of the Dixmier-Douady class.

An unitary descent datum is  $(U_\alpha, \maL_{\alpha\beta},  \mu_{\alpha\beta\gamma})$ together with a choice of metric on each $\maL_{\alpha\beta}$  such that each
$\mu_{\alpha\beta\gamma}$ is an isometry.
A notion of unitary equivalence of two unitary descent data on the same
open cover $U_{\alpha}$ is obtained from the notion of equivalence above by requiring that each
line bundle $S_{\alpha}$ is Hermitian and each $\lambda_{\alpha \beta}$ is an isometry.  The definition
of $2$-morphisms is modified by requiring each $\nu_{\alpha}$ to be an isometry.
It is clear that the restriction of a unitary descent datum to a refinement is again unitary.
Then a unitary gerbe is an equivalence class of unitary descent data in the sense described above.

\begin{lemma}\label{descconnection}
Let  $(U_\alpha, \maL_{\alpha\beta}, \mu_{\alpha\beta\gamma})$  be a descent datum on $M$.  There exists  a collection $(\nabla_{\alpha\beta})$ of connections on $(\maL_{\alpha\beta})$ such that for any $(\alpha, \beta, \gamma)\in \Lambda^3$ with $U_{\alpha \beta \gamma} \ne \emptyset$:
$$
\mu_{\alpha\beta\gamma}^* \nabla_{\alpha\gamma} = \nabla_{\alpha\beta}\otimes \id + \id \otimes \nabla_{\beta\gamma}.
$$
If the descent datum is unitary  each $\nabla_{\alpha\beta}$ can be chosen Hermitian.
\end{lemma}


\begin{lemma} \label{lem:curv}
Let $(\nabla_{\alpha\beta})$ be as above, and denote by $\omega_{\alpha\beta}=\nabla^2_{\alpha\beta}$ the curvatures
of the connections $\nabla_{\alpha \beta}$. Then
there exists a collection   of differential $2$-forms $\omega_\alpha\in \Omega^2(U_\alpha)$ such that
$$
\omega_{\alpha\beta}=\omega_{\alpha} - \omega_{\beta}, \quad \text{ for } U_{\alpha\beta} \not = \emptyset.
$$
\end{lemma}

%

\begin{definition} The collection  $\nabla^\maL= (\nabla_{\alpha \beta})$ is a connective structure on $\maL$,
while $\omega^\maL= (\omega_{\alpha})$ is a curving (compatible with $\nabla^\maL$). We  say that  the collection $ (\nabla_{\alpha \beta}, \omega_{\alpha})$ is a connection on the descent
datum $(U_\alpha, \maL_{\alpha\beta}, \mu_{\alpha\beta\gamma})$.

\end{definition}

If $\xi \in \Omega^2(M)$ then $\omega^\maL+\xi:=(\omega_{\alpha}+\xi)$
is another curving compatible with $\nabla^\maL$.
A connection on a descent datum yields a connection on its restriction to a refinement
in an obvious manner. A connection will be  identified  with its restriction.
%
%



\begin{definition}
Let  $\varsigma= (S_{\alpha}, \lambda_{\alpha \beta})$ be an isomorphism of descent datum $(U_\alpha, \maL_{\alpha\beta}, \mu_{\alpha\beta\gamma})$
with connective structure $\{\nabla_{\alpha \beta}\}$
and descent datum  $(U_\alpha, \maL_{\alpha\beta}', \mu_{\alpha\beta\gamma}')$
with connective structure $\{\nabla_{\alpha \beta}' \}$. A connection on $\varsigma$  is   a collection of connections $(\nabla_{\alpha} )$ on $(S_\alpha)$, satisfying the equality
\begin{equation}\label{nabla'}
(\lambda_{\alpha \beta})^*\nabla_{\alpha \beta}' =  \nabla_{\alpha}^* \otimes \id \otimes \id +\id \otimes \nabla_{\alpha \beta} \otimes \id +\id \otimes \id \otimes \nabla_{\beta}.
\end{equation}
Here  $\nabla_{\alpha}^*$ denotes the dual connection on $S_{\alpha}^{-1}$.
\end{definition}
If $ \{\nabla_\alpha\}$ is a connection on $\varsigma$ and $\xi \in \Omega^1(M)$ then $\{\nabla_\alpha+\xi\}$ is another such connection.

Let $\varsigma= (S_{\alpha}, \lambda_{\alpha \beta}, \nabla_{\alpha} )$
be a morphism of descent data with a connection, and
let $\varpi_\alpha =\nabla_{\alpha}^2$ be the associated collection of curvatures.

\begin{lemma}\label{curvature1morphism}
There exists a global form $c(\varsigma) \in \Omega^2(M)$ such that
\begin{equation*}
c(\varsigma)|_{U_\alpha} =  \omega_{\alpha}'+\varpi_{\alpha}-\omega_{\alpha}.
\end{equation*}
\end{lemma}

We say that $\varsigma$ is a flat morphism if $c(\varsigma)=0$, i.e. if
\begin{equation}\label{omega'}\
\omega_{\alpha} =\omega_{\alpha}'+\varpi_{\alpha}.
\end{equation}

\begin{lemma}
Let $\varsigma$ be an isomorphism between the  two gerbes $\maL$ and $\maL'$ with connective structures.
Then there exists a connection   on $\varsigma$.
\end{lemma}

\begin{proof} We can assume that the gerbes and the isomorphism are represented by data over the same good open cover $(U_\alpha)_\alpha$. So $(U_\alpha, S_\alpha, \lambda_{\alpha\beta})$ is an isomorphism between $(U_\alpha, \maL_{\alpha\beta}, \mu_{\alpha\beta\gamma})$ and $(U_\alpha, \maL'_{\alpha\beta}, \mu'_{\alpha\beta\gamma})$. We are given connective structures $(\nabla_{\alpha\beta})$ and $(\nabla'_{\alpha\beta})$ on these descent data for $\maL$ and $\maL'$ respectively, and we fix  a collection of connections $(\nabla_\alpha)_\alpha$ on the line bundles $(S_\alpha)_\alpha$. Denote by $\bar\nabla_{\alpha\beta}$ the connection defined out of the collections $(\nabla'_{\alpha\beta})$ and $(\nabla_\alpha)$ by Equation \ref{nabla'}:
\[
\bar{\nabla}_{\alpha \beta}' := (\lambda_{\alpha \beta}^{-1})^* \left(\nabla_{\alpha}^* \otimes \id \otimes \id +\id \otimes \nabla_{\alpha \beta} \otimes \id +\id \otimes \id \otimes \nabla_{\beta}\right).
\]
Then the difference $\nabla_{\alpha\beta} - \bar\nabla_{\alpha\beta}$ is a Cech $2$-cocyle with values in the soft sheaf of $1$-forms, hence there exists a collection $(\varphi_\alpha)_\alpha$ of $1$-forms on $(U_\alpha)_\alpha$ such that over $U_{\alpha\beta}\neq \emptyset$, we have
$$
\nabla_{\alpha\beta} - \bar\nabla_{\alpha\beta} = \varphi_\alpha - \varphi_\beta.
$$
It is then easy to check that the collection of connections $(\nabla_\alpha + \varphi_\alpha)_\alpha$ allows to endow $\varsigma$ with the structure of an isomorphism of our gerbes with connective structures.

\end{proof}

If $(S_{\alpha}, \lambda_{\alpha \beta}, \nabla_{\alpha} )$ and  $
(S_{\alpha}', \lambda_{\alpha \beta}', \nabla_{\alpha}' )$
are  two  morphisms between the descent data with connections
as above, the $2$-morphisms between them are thus the same as the $2$-morphisms
between $(S_{\alpha}, \lambda_{\alpha \beta})$ and  $
(S_{\alpha}', \lambda_{\alpha \beta}')$.

Assume that we are given two gerbes $\maL$ and $\maL'$ and two  morphisms $\varsigma$, $\varsigma'$ between them.
On an open cover $\{U_\alpha \}_{\alpha\in \Lambda}$ represent $\varsigma$ by $(S_{\alpha}, \lambda_{\alpha \beta}, \nabla_{\alpha} )$ and  $\varsigma' $ by
$(S_{\alpha}', \lambda_{\alpha \beta}', \nabla_{\alpha}' )$. For a $2$-morphism $\nu \colon  \varsigma \to \varsigma'$, we  define the $1$-form
$c(\nu)|_{U_\alpha} $ by
\begin{equation*}
c(\nu)|_{U_\alpha} = \nu_\alpha^* \nabla'_{\alpha}-\nabla_\alpha
\end{equation*}

\begin{lemma}
The collection of local $1$-forms $(c(\nu)|_{U_\alpha})_\alpha$ defines a global form $c(\nu) \in \Omega^1(M)$.
\end{lemma}

\begin{proof}
Indeed, notice that the isomorphism
$$
\nu_\alpha\otimes \id \otimes \nu_\beta^{-1} : S_\alpha\otimes \maL_{\alpha\beta}'\otimes S_\beta^{-1} \longrightarrow S'_\alpha\otimes \maL_{\alpha\beta}'\otimes {S'}_\beta^{-1}
$$
pulls back the connection $\nabla'_\alpha\otimes\id + \id\otimes \nabla'_{\alpha\beta}\otimes \id + \id\otimes {\nabla'}_\beta^{-1}$ to the connection $\nabla_\alpha\otimes\id + \id\otimes \nabla'_{\alpha\beta}\otimes \id + \id\otimes {\nabla}_\beta^{-1}$. From this the following relation holds
$$
(\nu_\alpha^*\nabla'_\alpha-\nabla_\alpha )\otimes \id + \id\otimes ((\nu_\beta^{-1})^*{\nabla'}_\beta^{-1}-\nabla_\beta^{-1} ) = 0,
$$
and this ends  the proof.
\end{proof}

We now state the following lemma which is proved in \cite{BenaGoro}[page 7]:

\begin{lemma}\label{defomega} Let $\maL$ be a gerbe represented by a descent
datum $(U_\alpha, \maL_{\alpha\beta}, \mu_{\alpha\beta\gamma})$.
Choose a connection on $\maL$ represented by the collection of connections $(\nabla_{\alpha \beta}, \omega_{\alpha})$
on the descent datum.
\begin{enumerate}
\item There exists a well-defined closed form $\varOmega\in \Omega^3(M)$ --
curvature $3$-form of the connection --
such that $\varOmega|_{U_{\alpha}}=\frac{d\omega_{\alpha}}{2\pi i}$.

\item  Let $(\nabla_{\alpha \beta}', \omega_{\alpha}')$ be  another
connection on $\maL$ and let $\varOmega'$ be the corresponding curvature $3$-form.
Then there exists a canonical  $\eta \in \Omega^2(M)/d\Omega^1(M)$ such that
$\varOmega'= \varOmega+ d\eta$.

\item Let $(\nabla_{\alpha \beta}, \omega_{\alpha})$, $(\nabla_{\alpha \beta}', \omega_{\alpha}')$ and
$(\nabla_{\alpha \beta}'', \omega_{\alpha}'')$ be $3$  connections on $\maL$ with the corresponding
curvature $3$-forms $\varOmega$, $\varOmega'$ and $\varOmega''$. Let $\eta, \eta', \eta'' \in
\Omega^2(M)/d\Omega^1(M)$ be the canonical elements constructed above such that $
\varOmega'-\varOmega=d \eta$, $\varOmega''-\varOmega'=d\eta'$, $\varOmega''-\varOmega=d\eta''$.
Then $\eta''=\eta+\eta'$.
\end{enumerate}
\end{lemma}

Choosing a unitary descent datum for $\maL$ and a hermitian connection $(\nabla_{\alpha \beta}, \omega_{\alpha})$, a classical argument shows that the $3$-form $\varOmega$ is a  de Rham representative of the Dixmier-Douady class of the gerbe. An argument similar to the proof of Lemma \ref{descconnection} shows the following proposition which gathers the main properties that will be used in the sequel. The proofs are  straightforward consequences of the definitions and are  omitted:

\begin{proposition}\label{Relations}\
\begin{enumerate}
\item Let $\varsigma \colon \maL \to \maL'$ be a morphism of gerbes with connections. Then
 \begin{equation*} dc(\varsigma) =\varOmega(\maL')-\varOmega(\maL).\end{equation*}
\item Let $\varsigma \colon \maL \to \maL'$, $\varsigma' \colon \maL' \to \maL''$ be morphisms of gerbes with connections. Then
\begin{equation*}
c(\varsigma' \circ \varsigma) = c(\varsigma')+c(\varsigma).
\end{equation*}

\item   Let $\varsigma \colon \maL \to \maL'$, $\varsigma' \colon \maL \to \maL'$ be morphisms of gerbes with connections, and $\nu \colon \varsigma \to \varsigma'$ a $2$-morphism. Then
\begin{equation*} dc(\nu) =c(\varsigma)-c(\varsigma').
\end{equation*}
\item The correspondence $\nu \mapsto c(\nu)$ is additive under the horizontal and vertical compositions of $2$-morphisms.
\end{enumerate}
\end{proposition}

\subsection{Twisted and transversally twisted bundles}\label{hotwbu} $\;$

\begin{definition}
A descent datum for a twisted vector bundle $\maE$ consists of a descent datum
$(U_{\alpha}, \maL_{\alpha\beta}, \mu_{\alpha\beta\gamma})$ for a gerbe $\maL$ together with a collection
$(\maE_\alpha\to U_\alpha)_{\alpha\in \Lambda}$ of vector bundles  and  a collection of vector bundle isomorphisms $\varphi_{\alpha\beta}: \maE_\alpha\otimes \maL_{\alpha\beta} \stackrel{\cong}{\rightarrow} \maE_\beta$ such that for every $\alpha$, $\beta$, $\gamma$ the following diagram commutes

\begin{equation*}
\begin{CD}
\maE_\alpha\otimes \maL_{\alpha\beta}\otimes \maL_{\beta\gamma} @>\id \otimes \mu_{\alpha\beta\gamma}>>
\maE_\alpha \otimes \maL_{\alpha\gamma}\\
@V\varphi_{\alpha\beta}\otimes \id VV    @VV\varphi_{\alpha\gamma}V \\
\maE_\beta\otimes \maL_{\beta\gamma} @>\ \ \varphi_{\beta\gamma}\ \ >> \maE_{\gamma}
\end{CD}
\end{equation*}

\end{definition}

Restriction of the descent datum for $\maE$ to a refinement is given by the restriction of the descent datum
for $\maL$ together with restriction of the vector bundles $\maE_{\alpha}$ and the isomorphisms $\varphi_{\alpha\beta}$.

\begin{definition}
An (iso)morphism between two descent data $( U_{\alpha}, \maL_{\alpha\beta}, \mu_{\alpha\beta\gamma}, \maE_{\alpha}, \varphi_{\alpha \beta})$ and
$(  U_{\alpha}', \maL_{\alpha\beta}', \mu_{\alpha\beta\gamma}', \maE_{\alpha}', \varphi_{\alpha \beta}')$ is given by the collection  $(\rho_{\alpha}, S_{\alpha}, \lambda_{\alpha \beta})$ where $( S_{\alpha}, \lambda_{\alpha \beta})$ is an isomorphism between
$(  U_{\alpha}, \maL_{\alpha\beta}, \mu_{\alpha\beta\gamma})$ and
$( U_{\alpha}', \maL_{\alpha\beta}', \mu_{\alpha\beta\gamma}')$ and $\rho_{\alpha} \colon \maE_{\alpha} \otimes S_{\alpha} \to \maE_{\alpha}'$ is a collection of (iso)morphisms such that the following diagram commutes
\begin{equation*}
\begin{CD}
\maE_{\alpha} \otimes S_{\alpha}\otimes S_{\alpha}^{-1}\otimes \maL_{\alpha \beta} \otimes S_{\beta}
@>\varphi_{\alpha \beta} \otimes \id>>
\maE_{\beta} \otimes S_{\beta}\\
@V\rho_{\alpha}\otimes\lambda_{\alpha \beta }VV @VV\rho_{\beta}V\\
\maE_{\alpha }' \otimes \maL_{\alpha \beta }' @>\ \varphi_{\alpha \beta \gamma}'\ >> \maE_{\beta}'
\end{CD}
\end{equation*}

An isomorphism  between two descent data on two different covers is defined as an isomorphism
between their restrictions on a common refinement.

\end{definition}
\begin{definition}\label{morphisms for bundles}
A $2$-morphism between two isomorphisms $(\rho_{\alpha}, S_{\alpha}, \lambda_{\alpha \beta})$ and
$(\rho'_{\alpha}, S'_{\alpha}, \lambda'_{\alpha \beta})$ is
the $2$-morphism  between the corresponding isomorphisms of the gerbe descent data, i.e. collection of
morphisms $\nu_\alpha \colon S_\alpha \to S'_\alpha$  subject to \eqref{2morphism}, such that we have in addition
\begin{equation*}
\rho'_\alpha \circ (\id \otimes \nu_\alpha) = \rho_\alpha.
\end{equation*}
\end{definition}

\bigskip
A twisted bundle is then defined as an equivalence  class
of descent data of twisted vector bundles. ``Forgetting'' the bundle data we obtain from the descent datum for a twisted vector
bundle a descent datum for a gerbe, and the same applies to morphisms and $2$-morphisms.
We say that
$\maE$ is an $\maL$-twisted vector bundle if ``forgetting'' the bundle data one obtains the equivalence class
of the gerbe descent data defining $\maL$.


Assume now that the gerbe $\maL$ is unitary. An Hermitian descent datum for $\maE$  consists of a unitary descent datum $(U_{\alpha}, \maL_{\alpha \beta}, \mu_{\alpha \beta \gamma})$ for $\maL$   and a collection $h_{\alpha}$ of metrics on $\maE_{\alpha}$ such that the maps $\varphi_{\alpha \beta}$ are isometries. One obtains a notion of isomorphism of Hermitian descent data by requiring the $\rho_{\alpha}$'s to be isometries. An Hermitian twisted bundle is then an equivalence class of Hermitian descent data.

\begin{remark}
Given a gerbe $(U_{\alpha}, \maL_{\alpha\beta}, \mu_{\alpha\beta\gamma})$ on $M$, it is well known that
 (finite dimensional) twisted vector bundles exist if and only if the gerbe is torsion
(see e.g. \cite{DK, Karoubi} and references therein).
\end{remark}


The following lemma from  \cite{BenaGoro} will be used later.

\begin{lemma} \label{extend}
Let $\maL$ be a gerbe on $M$ and $\maE$ an $\maL$-twisted bundle. Let $(U_{\alpha}, \maL_{\alpha \beta},
\mu_{\alpha \beta \gamma})$ be any descent datum for $\maL$. Then there exists a descent datum for $\maE$ isomorphic
to the one of the form
$(U_{\alpha}, \maL_{\alpha \beta}, \mu_{\alpha \beta \gamma}, \maE_{\alpha}, \varphi_{\alpha \beta})$.
\end{lemma}

We now fix a gerbe $\maL$ with a descent datum $(U_\alpha, \maL_{\alpha\beta}, \mu_{\alpha\beta\gamma})$ on the smooth manifold $M$ together with a twisted vector bundle $\maE$ represented by $(\maE_\alpha, \varphi_{\alpha\beta})$. We denote by $\maA_\alpha$ the collection of bundles of algebras
$$
\maA_\alpha := \End( \maE_\alpha), \quad \alpha\in \Lambda.
$$
For any $U_{\alpha\beta}\not = \emptyset$, we have a canonical vector bundle isomorphism over $U_{\alpha\beta}$
$$
\rho_{\alpha\beta} : \End( \maE_\alpha \otimes \maL_{\alpha\beta}) \longrightarrow \maA_\alpha,
$$
extending the canonical isomorphism $\End(\maL_{\alpha\beta}) \simeq U_{\alpha\beta}\times \C$. Therefore, the bundle isomorphism $\varphi_{\alpha\beta}$ together with the identification $\rho_{\alpha\beta}$ induce the isomorphism of algebra bundles over $U_{\alpha\beta}$ given by
$$
\maA_\beta \stackrel{\varphi_{\alpha\beta}^*}{\longrightarrow} \End( \maE_\alpha \otimes \maL_{\alpha\beta})\stackrel{\rho_{\alpha\beta}}{\longrightarrow} \maA_\alpha.
$$
We denote  this isomorphism simply by $\varphi_{\alpha\beta}^*$. It is then easy to check that
$$
\varphi_{\alpha\beta}^* \circ \varphi_{\beta\gamma}^* = \varphi_{\alpha\gamma}^*, \quad \text{ over } U_{\alpha\beta\gamma}.
$$
Therefore, the collection $\{\maA_\alpha\}$   defines a bundle $\maA $ of algebras over $M$, which is denoted  $\End(\maE)$. It is an easy exercise to  check that the isomorphism class of $\End(\maE)$ depends only on the isomorphism class of the $\maL$-twisted vector bundle $\maE$.

\begin{definition}
The bundle $\maA = \End(\maE) $ is called the Azumaya bundle associated with the $\maL$-twisted vector bundle $\maE$.
\end{definition}

If  $\maE$ is a Hermitian twisted bundle, each of the bundles $\End(\maE_{\alpha})$ is a bundle of *-algebras,
with the *-operation given by taking the adjoint endomorphism. This induces a structure of a bundle of *-algebras on $\End(\maE)$.

Note that for every $\alpha$ we have the trace $ \tr_{\alpha} \colon \End(\maE_{\alpha}) \to C^{\infty}(U_{\alpha})$.
For $a \in \Gamma(U_{\alpha \beta}; \maA_{\beta})$ $\tr_{\alpha}(\varphi_{\alpha \beta}^*(a))= \tr_{\beta}(a)$.
Therefore we obtain the trace $\tr \colon \End(\maE) \to C^{\infty}(M)$ defined by
\begin{equation*} \tr(a)|_{U_{\alpha}}= \tr_{\alpha}(a|_{U_{\alpha}}) \text{ for } a \in \End(\maE)
\end{equation*}
If the bundle $\maE$ is $\Z_2$-graded then $\End(\maE)$ is a bundle
of $\Z_2$-graded algebras with a supertrace
$$
\str \colon \End(\maE) \longrightarrow C^{\infty}(M).
$$
Notice that, more generally, if $\maE$ and $\maE'$ are $\maL$-twisted bundles, we have a well defined
vector bundle $\Hom (\maE, \maE')$, defined similarly by $\left.\Hom (\maE, \maE') \right|_{U_{\alpha}}=\Hom (\maE|_{U_{\alpha}}, \maE'|_{U_{\alpha}})$.

\begin{definition}
Let as before $(U_\alpha, \maL_{\alpha\beta}, \mu_{\alpha\beta\gamma}, \maE_\alpha, \varphi_{\alpha\beta}) $ be a descent datum for the $\maL$-twisted vector bundle $\maE$ on $M$. A connection on $(U_\alpha, \maL_{\alpha\beta}, \mu_{\alpha\beta\gamma}, \maE_\alpha, \varphi_{\alpha\beta}) $  is a collection $(\nabla_\alpha, \nabla_{\alpha \beta}, \omega_{\alpha})$ where $(\nabla_{\alpha \beta}, \omega_{\alpha})$ is a connection on $(U_\alpha, \maL_{\alpha\beta}, \mu_{\alpha\beta\gamma})$ and each $\nabla_{\alpha}$ is a connection on $\maE_{\alpha}$ such that the identities

\begin{equation}\label{connectionon}
\varphi_{\alpha\beta}^* \nabla_\beta = \nabla_\alpha\otimes \id + \id \otimes \nabla_{\alpha\beta}
\end{equation}
$\text{hold whenever } U_{\alpha\beta}\not = \emptyset.$
\end{definition}

Again the following lemma is proved in \cite{BenaGoro}[page 11]:

\begin{lemma}\label{connection existence}
Let $(U_\alpha, \maL_{\alpha\beta}, \mu_{\alpha\beta\gamma}, \maE_\alpha, \varphi_{\alpha\beta}) $ be a descent datum for an $\maL$-twisted vector bundle $\maE$.
Then every connection $(\nabla_{\alpha\beta}, \omega_{\alpha})$ on the descent datum  $(U_\alpha, \maL_{\alpha\beta}, \mu_{\alpha\beta\gamma})$ for $\maL$ can be extended to a connection for the descent datum for $\maE$.
\end{lemma}

\begin{definition}\label{bundle morphism}
An  isomorphism  between two descent
data   $( U_{\alpha}, \maL_{\alpha\beta}, \mu_{\alpha\beta\gamma}, \maE_{\alpha}, \varphi_{\alpha \beta})$ and
$(  U_{\alpha}', \maL_{\alpha\beta}', \mu_{\alpha\beta\gamma}', \maE_{\alpha}', \varphi_{\alpha \beta}')$
for $\maE$, with connections $(\nabla_\alpha, \nabla_{\alpha \beta}, \omega_{\alpha}) $ and
$(\nabla_\alpha', \nabla_{\alpha \beta}', \omega_{\alpha}') $, is given by the collections
 $S=(\rho_{\alpha}, S_{\alpha}, \nabla_\alpha^S, \lambda_{\alpha \beta})$ where
$(\rho_{\alpha}, S_{\alpha}, \lambda_{\alpha \beta})$ is an isomorphism between the descent data
without connections,  $\nabla_\alpha^S$ are connections on $S_\alpha$ such that
$\varsigma= \left(S_{\alpha}, \lambda_{\alpha \beta}, \nabla_\alpha^S\right)$ is a morphism of the corresponding
gerbe descent data with connective structures (i.e. the equation
\eqref{nabla'} is satisfied) and the equality
\begin{equation} \label{eq:s}
\rho_{\alpha}^* \nabla_{\alpha}' =\nabla_{\alpha} \otimes \id + \id \otimes \nabla_{\alpha}^S
 \end{equation}
holds (in the notations of \eqref{nabla'}).\\
An isomorphism is called flat if the corresponding isomorphism of gerbes $\varsigma$ is flat.
\end{definition}

A connection on a twisted vector bundle is then a choice  of connections on each
descent datum of this twisted bundle and lifting of  isomorphisms of descent data to
flat isomorphisms of descent data with connections. So, every connection
on a gerbe $\maL$ can be extended to a connection on any $\maL$-twisted vector bundle.

Recall also the following propositions \ref{global curvature} and \ref{twchern}, cf. e.g. \cite{BenaGoro}[pages 12-13].

\begin{proposition}\label{global curvature}
Let $\maE$ be an $\maL$-twisted bundle with connection. Choose a descent datum $(U_\alpha, \maL_{\alpha\beta}, \mu_{\alpha\beta\gamma}, \maE_\alpha, \varphi_{\alpha\beta}) $   with a connection $(\nabla_\alpha, \nabla_{\alpha \beta}, \omega_{\alpha})$ representing $\maE$. Then the collection $(\theta_\alpha + \omega_\alpha)$, where $\theta_\alpha=\nabla_\alpha^2$ is the curvature of $\nabla_{\alpha}$, defines a global differential $2$-form $\theta$ on $M$ with coefficients in the Azumaya bundle $\maA=\End(\maE)$. This form is independent of the choice of the representing descent datum.
\end{proposition}


In the notations above let $\nabla \colon \maA \to \Omega^1(M, \maA)$ be the connection defined for a fixed descent datum by
\begin{equation}
(\nabla \xi)|_{U_{\alpha}}=[\nabla_{\alpha}, \xi] .
\end{equation}

It is easy to see that $\nabla$ is well defined and is a derivation with respect to the product on $\maA$.
Note that
\begin{equation}
\nabla^2=[\theta, \cdot] \text{ and }\nabla \theta = \varOmega
\end{equation}
where $\varOmega$ is the $3$-curvature form of the connection on $\maL$, see Lemma \ref{defomega}. In the following proposition, we use the notations of Appendix \ref{TwistedCohom}, in particular the twisted differential $d_\varOmega$.

\begin{proposition}\label{twchern}\
Let $\maE$ be an $\maL$-twisted vector bundle and $\nabla$ a connection on $\maE$.
Set $\Ch_{\maL}(\nabla)= \tr e^{-u\theta} \in \Omega^*(M)[u].$
Then
\begin{enumerate}
\item $ d_{\varOmega}\Ch_{\maL}(\nabla)=0$.
\item The class of $ \Ch_{\maL}(\nabla)$ in $H_{\varOmega}^*(M)$ is independent of choice of connections.\\
Namely, assume that we are given a different connection $\nabla'$ on $\maE$ (and hence on $\maL$)
 and let $\varOmega'$ be the associated $3$-curvature form. Then $I([\Ch_{\maL}(\nabla)])=[\Ch_{\maL}(\nabla')]$, where $I$ is the canonical isomorphism of $H^*_\varOmega (M)$ with $H^*_{\varOmega '} (M)$.
\end{enumerate}
\end{proposition}
Identifying as before the canonically isomorphic twisted cohomology spaces $H_{\varOmega}^*(M)$, we  denote the class of $\Ch_{\maL}(\nabla)$ by $\Ch_{\maL}(\maE)$.

Consider now an isomorphism $S$ of twisted bundles $\maE$, $\maE'$ represented by the descent data as in Definition \ref{bundle morphism}.
Let $\theta$ and $\theta'$ be the curvature forms for $\maE$ and $\maE'$ respectively (cf. Proposition \ref{global curvature}). Let $\cA$, $\cA'$ be the
Azumaya bundles of $\maE$ and $\maE'$ and $S^* \colon \cA' \to \cA$ the induced isomorphism. Then an immediate consequence of  \eqref{eq:s} gives
$$  S^*\theta' -\theta  = c(\varsigma).
$$

We will also need  to use superconnections on twisted bundles and we now briefly indicate the modifications
which need to be made to the notion of connection to obtain that of superconnection.
Assume that we are given a gerbe $\maL$  and  a $\Z_2$-graded $\maL$-twisted vector bundle $\maE=\maE^+\oplus \maE^-$.
Let  $(U_\alpha, \maL_{\alpha\beta}, \mu_{\alpha\beta\gamma}, \maE_\alpha, \varphi_{\alpha\beta}) $ be a descent datum for $\maE$.
\begin{definition}
 A superconnection $\A$ on the descent datum   is $\A= (\A_{\alpha}, \nabla_{\alpha \beta}, \omega_{\alpha})$ where
$(\nabla_{\alpha \beta}, \omega_{\alpha})$ is a connection on the descent datum for $\maL$ and
each $\A_{\alpha}$ is a superconnection on  $\maE_{\alpha}$ satisfying the relations
$$
\varphi_{\alpha\beta}^* \A_\beta = \A_\alpha\otimes \id + \id \otimes \nabla_{\alpha\beta},\quad \text{ for } U_{\alpha\beta}\not = \emptyset.
$$
\end{definition}

Each superconnection $\A_{\alpha}$ can be written as $ \A_{\alpha}= \sum_{k\ge 0} \A_{\alpha}^{[k]}$ where $\A_{\alpha}^{[k]} \in \Omega^k(U_{\alpha}; \End(\maE_{\alpha})^-)$ for $k$-even, $\A_{\alpha}^{[k]} \in \Omega^k(U_{\alpha}; \End(\maE_{\alpha})^+)$ for $k$-odd, $k\ne 1$, and $\A_{\alpha}^{[1]}$ is a grading preserving connection on $\maE_{\alpha}$. It is easy to see that for each $k \ne 1$ there exists a form $\A^{[k]} \in \Omega^k(M; \End(\maE))$ such that $\A^{[k]}|_{U_{\alpha}}=\A_{\alpha}^{[k]}$. For $k=1$, $(\A_{\alpha}^{[1]},  \nabla_{\alpha \beta}, \omega_{\alpha})$ defines a  connection on (the descent datum of) $\maE$.

Let $\varOmega$ be as before the curvature $3$-form  of the connection on $\maL$. Define the curvature $\theta^{\A}$  of the superconnection $\A$ by
\begin{equation*}
\theta^{\A}|_{U_{\alpha}}=(\A_{\alpha})^2+\omega_{\alpha}.
\end{equation*}

We introduce  a formal variable $u^{1/2}$  of degree $-1$ such that $(u^{1/2})^2=u$ and define the new superconnection
\begin{equation*}
\A_{u^{-1}} := \sum u^{(k-1)/2}\A^{[k]}.
\end{equation*}
Then the curvature $\theta^{{\A}_{u^{-1}}}$ of the  superconnection $\A_{u^{-1}}$ is given by
\begin{equation*}
\theta^{{\A}_{u^{-1}}}|_{U_{\alpha}}=(\A_{u^{-1}, \alpha})^2+\omega_{\alpha}.
\end{equation*}

Notice then that $u\theta^{{\A}_{u^{-1}}} \in \Omega^{even}(M, \End (\maE)^+)[u] +u^{1/2}\Omega^{odd}(M, \End (\maE)^-)[u]$. We therefore deduce that
$$
\exp\left(-u\theta^{{\A}_{u^{-1}}}\right) \; \in \; \Omega^{even}(M, \End (\maE)^+)[u] +u^{1/2}\Omega^{odd}(M, \End (\maE)^-)[u],
$$
and the differential form $\str  \exp \left(-u\theta^{{\A}_{u^{-1}}}\right)$ belongs to $\Omega^{even}(M)[u]$.
The following is an analogue of  Proposition \ref{twchern} for the superconnections with the essentially identical proof.

\begin{proposition}\label{chernsuperconnection}\

\begin{enumerate}
\item Set $\Ch_{\maL}(\A)= \str \exp \left(-u\theta^{{\A}_{u^{-1}}}\right)$. Then $ d_{\varOmega}\Ch_{\maL}(\A)=0$
\item
The class of $ \Ch_{\maL}(\A)$ in $H_{\varOmega}(M)$ is independent of choice of superconnection. \\Specifically, assume we are given a different superconnection $\A'$ on $\maE$ (and therefore a different connection on $\maL$) and let $\varOmega'$ be the associated $3$-curvature form. Then $I([\Ch_{\maL}(\A)])=[\Ch_{\maL}(\A')]$, where $I$ is the canonical isomorphism of $H^*_{\varOmega} (M)$ with $H^*_{\varOmega '} (M)$.
\end{enumerate}
\end{proposition}
\medskip

Let $\pi: M \to B$ be a smooth fibration, then we need to extend the previous definitions  to gerbes on $M$ which are pulled back from $B$ as follows. Let $\maL$ be a gerbe on $B$.

\begin{definition}\
\begin{itemize}
\item
 A descent datum
for a transversally $\maL$-twisted bundle $\maE$ on $M$ consists of the descent datum
$(U_{\alpha}, \maL_{\alpha\beta}, \mu_{\alpha\beta\gamma})$ for $\maL$ together with a collection
$(\maE_\alpha\to \pi^{-1}U_\alpha)_{\alpha\in \Lambda}$ of vector bundles and a collection of vector bundle isomorphisms $\varphi_{\alpha\beta}\colon \maE_\alpha\otimes \pi^*\maL_{\alpha\beta} \cong \maE_\beta$  so that
$$
(\pi^{-1}U_{\alpha}, \pi^*\maL_{\alpha\beta}, \pi^*\mu_{\alpha\beta\gamma}, \maE_\alpha, \varphi_{\alpha\beta} )
$$
is a descent datum for a twisted vector bundle on $M$.
\item An (iso)morphism between two such descent data
 is given by the collection  $(\rho_{\alpha}, S_{\alpha}, \lambda_{\alpha \beta})$ where $( S_{\alpha}, \lambda_{\alpha \beta})$ is an isomorphism between
$(  U_{\alpha}, \maL_{\alpha\beta}, \mu_{\alpha\beta\gamma})$ and
$( U_{\alpha}', \maL_{\alpha\beta}', \mu_{\alpha\beta\gamma}')$  and $\rho_{\alpha} \colon \maE_{\alpha} \otimes \pi^*S_{\alpha} \to \maE_{\alpha}'$ is such that $(\rho_{\alpha}, \pi^*S_{\alpha}, \pi^*\lambda_{\alpha \beta})$ is
an (iso)morphism in the usual sense between  $(\pi^{-1} U_{\alpha}, \pi^*\maL_{\alpha\beta}, \pi^*\mu_{\alpha\beta\gamma}, \maE_{\alpha}, \varphi_{\alpha \beta})$ and
$( \pi^{-1} U_{\alpha}', \pi^*\maL_{\alpha\beta}', \pi^*\mu_{\alpha\beta\gamma}', \maE_{\alpha}', \varphi_{\alpha \beta}')$.
\end{itemize}
\end{definition}
The notion of $2$-morphism is defined similarly. A transversally twisted bundle is then an equivalence class of descent data. \\

Let $\Vect_{\pi^*\maL}(M)$ denote the set of isomorphism classes of
all $\pi^*\maL$-twisted bundles on $M$ and $\Vect_{\pi^*\maL}^\nmid (M)$ denote the set of isomorphism classes of
all transversally $\maL$-twisted bundles. Then according to Lemma \ref{extend}, every class in $\Vect_{\pi^*\maL}(M)$ can be represented by a descent datum of a transversally $\maL$-twisted bundle and the obvious forgetful map
$$
\Vect_{\pi^*\maL}^\nmid (M) \longrightarrow \Vect_{\pi^*\maL}(M),
$$
is thus surjective. Recall that (see \cite{MMS1},  \cite{MMS2}) $\Vect_{\pi^*\maL}(M)$  is not empty if and only if the class $\pi^* [\maL]$ is torsion in $H^3(M, \Z)$.
We then deduce immediately the following Proposition:
\begin{proposition}
Transversely $\maL$-twisted vector bundles exist on $M$ if and only if the class $\pi^* [\maL]$ is torsion in $H^3(M, \Z)$.
\end{proposition}
Notice that the above forgetful map is not injective. Indeed, if $(U_{\alpha}, \maL_{\alpha\beta}, \mu_{\alpha\beta\gamma}, \maE_\alpha, \varphi_{\alpha \beta})$ is a descent datum for a transversally twisted bundle and
$S$ is a line bundle on $M$ then $(U_{\alpha}, \maL_{\alpha\beta}, \mu_{\alpha\beta\gamma}, \maE_\alpha \otimes S\vert_{\pi^{-1}U_\alpha}, \varphi_{\alpha \beta}\otimes \id)$ is another such descent datum. These data define the same
element of $\Vect_{\pi^*\maL}(M)$ but different elements of $\Vect_{\pi^*\maL}^\nmid (M)$, unless $S$ is a pull-back of a line bundle from $B$.

A connection on the descent datum $( U_{\alpha}, \maL_{\alpha\beta}, \mu_{\alpha\beta\gamma}, \maE_{\alpha}, \varphi_{\alpha \beta})$ for a transversely $\maL$-twisted bundle is a collection $(\nabla_{\alpha}, \nabla_{\alpha \beta}, \omega_{\alpha})$ where
$(\nabla_{\alpha \beta}, \omega_{\alpha})$ is a connection on the descent datum for $\maL$ and $\nabla_{\alpha}$ is a
connection on $\maE_{\alpha}$ such that $(\nabla_{\alpha}, \pi^*\nabla_{\alpha \beta}, \pi^*\omega_{\alpha})$
is a connection on the descent datum  $(\pi^{-1} U_{\alpha}, \pi^*\maL_{\alpha\beta}, \pi^*\mu_{\alpha\beta\gamma}, \maE_{\alpha}, \varphi_{\alpha \beta})$.
With these
definitions one can now define a notion of connection on the transversally $\maL$-twisted bundle in complete
analogy with the definition for the twisted bundles. If $\nabla$  is such a connection and $\varOmega$ is the curvature $3$-form of the gerbe $\maL$, one defines $\Ch_{\maL}(\nabla)$,  a closed form in $\left(\Omega^*(M), d_{\pi^*\varOmega}\right)$.
The analogues of Propositions \ref{global curvature} and \ref{twchern} hold in this context with the same proofs.

\subsection{Projective families of pseudodifferential operators}\label{Projective families}
We collect in this subsection several facts about   families of pseudodifferential operators, associated with our submersion $M\to B$ and with transversely $\maL$-twisted bundles. More details were given in \cite{BenaGoro} and we shall be brief. Let us  start with  the untwisted calculus.

%
%
%
%

We denote by $\Psi^m (M|B ; E)$ the space of classical  fiberwise pseudodifferential operators of order $\leq m$ on   $\pi$, acting on the sections of  the vector bundle $E$, which are {\underline{fiberwise properly supported}} over $M$
in the following sense. For $A \in \Psi^m (M|B ; E)$ let $S_A \subset M \times M \times B$ be the support of the Schwartz kernel of A. We require that the maps $p_i \colon S_A \to X  \times Y$, $i=1$, $2$ induced by the projections
on the $i$th factor in $X \times X \times Y$ are proper.
 As usual, we set
$$
\Psi (M|B ; E) := \bigcup_{m\in \Z} \Psi^m (M|B ; E) \text{ and }
 \Psi^{-\infty} (M|B; E) := \bigcap_{m\in \Z} \Psi^m (M|B ; E).
$$
Recall that composition endows each $\Psi(M|B ; E)$ with the structure of a filtered algebra and a module over $C^{\infty}(B)$ such that the composition is
$C^{\infty}(B)$-linear and $\Psi^{-\infty} (M|B; E)$ is an ideal in this algebra.
We have the following elementary general result which will be used in the sequel for different  submersions:
\begin{lemma} Assume that $\pi\colon X\to Y$ is a smooth submersion, $E$ is a vector bundle on $X$ and $L$ is a line bundle on $Y$.
Define a map $\chi_L\colon \Psi  (X|Y ; E) \to \Psi (X|Y ; E\otimes \pi^* L )$ by
\begin{equation*}
\chi_L(D) (e\otimes \pi^* (l)) = D(e)\otimes \pi^*l
\end{equation*}
for $D\in \Psi (X|Y ; \maE)$, $e \in \Gamma_c(E)$, $l \in \Gamma(L)$. Then $\chi_L$ is a well-defined  isomorphism of  algebras and $C^{\infty}(Y)$-modules. Moreover, we have for line bundles $L_1, L_2 \rightarrow Y$:
\begin{equation*}
\chi_{L_1 \otimes L_2} = \chi_{L_2} \circ \chi_{L_1}.
\end{equation*}
Assume that $A \in \Psi  (X|Y ; E)$ is positive elliptic selfadjoint, and $f$ is a Schwartz function.
Then
$f(\chi_L(A)) =\chi_L(f(A))$.
\end{lemma}

If we have two vector bundles $E$, $E'$ on $M$ we denote by $\Psi(M|B; E, E')$ the set of fiberwise pseudodifferential
operators $\Gamma_c(E) \to \Gamma(E')$. For a line bundle $L$ on $B$,  we again  have the isomorphism of $C^{\infty}(B)$-modules $\chi_L \colon   \Psi(M|B; E, E') \to \Psi(M|B; E\otimes \pi^*L, E'\otimes \pi^*L)$ defined
by the same formula.

We have the vertical cotangent bundle $T^*(M|B)= T^*M/(\Ker \pi_*)^{\perp}$. $\mathring{T}^*(M|B)$ denotes
(the total space of) this bundle with the zero section removed, and $p \colon \mathring{T}^*(M|B) \to M$ is the natural projection.
Recall that for $P\in \Psi^m (M|B ; E, E')$    the principal symbol $\sigma_m(P)$ is a positively  $m$-homogeneous smooth section over  $\mathring{T}^*(M|B) $  of the vector bundle $p^*\Hom( E, E')$. Then identifying the canonically isomorphic bundles $\Hom( E, E')$ and  $\Hom( E\otimes \pi^*L, E'\otimes \pi^*L)$ we have $\sigma_m(P) =\sigma_m(\chi_L(P))$.\\

\medskip

We now introduce the twist by a gerbe and assume that  $\pi^*\maL$-twisted vector bundles exist on $M$, then we know that transversely $\maL$-twisted bundles do exist as well, cf. Section \ref{hotwbu}.  Let $\maE$ be such a transversally $\maL$-twisted vector bundle  on $M$.
We fix a descent datum $(U_{\alpha}, \maL_{\alpha\beta}, \mu_{\alpha\beta\gamma}, \maE_\alpha, \varphi_{\alpha \beta})$ for $\maE$.
  For any $(\alpha, \beta)\in \Lambda^2$ with $U_{\alpha\beta}\not = \emptyset$
we have an isomorphism of filtered  algebras, respecting the $C^{\infty}(U_{\alpha \beta})$-module structure:
\begin{equation}\label{phi}
 \phi_{\alpha\beta}\colon \Psi (\pi^{-1}U_{\alpha\beta}|U_{\alpha\beta} ; \maE_\beta)\rightarrow \Psi (\pi^{-1}U_{\alpha\beta}|U_{\alpha\beta} ; \maE_\alpha),
\end{equation}

It is defined as the composition
\begin{equation*}
\Psi(\pi^{-1}U_{\alpha\beta}|U_{\alpha\beta} ; \maE_\beta) \overset{\psi_{\alpha \beta}}\to
\Psi (\pi^{-1}U_{\alpha\beta}|U_{\alpha\beta} ; \maE_\beta \otimes \pi^* \maL_{\alpha \beta}) \overset{\varphi_{\alpha \beta}}\to \Psi (\pi^{-1}U_{\alpha\beta}|U_{\alpha\beta} ; \maE_\alpha)
\end{equation*}
where $\psi_{\alpha \beta} =\chi_{\maL_{\alpha \beta}}^{-1}=\chi_{\maL_{\beta \alpha }}$.

Recall (cf \cite{bgv}) that  for every $\alpha \in \Lambda$ we have an infinite dimensional bundle $\pi_* \maE_{\alpha}$
on $U_{\alpha}$ defined by $\Gamma(V, \pi_*\maE_{\alpha}) = \Gamma(\pi^{-1}V, \maE_{\alpha})$, $V \subset U_{\alpha}$.
Over $U_{\alpha \beta}$ we have isomorphisms $\pi_*\varphi_{\alpha \beta} \colon \pi_*\maE_{\alpha} \otimes \maL_{\alpha \beta} \to \pi_* \maE_{\beta}$ defined by
\begin{equation*}
\pi_*\varphi_{\alpha \beta} (\xi \otimes l) = \varphi (\xi \otimes \pi^*l).
\end{equation*}
Here $\xi \in \Gamma(U_{\alpha \beta}, \pi_*\maE_{\alpha})= \Gamma(\pi^{-1}U_{\alpha \beta}, \maE_{\alpha})$, $l \in \Gamma(U_{\alpha \beta}, \maL_{\alpha \beta})$.

Note that the isomorphisms $\pi_*\varphi_{\alpha \beta} \colon \pi_*\maE_{\alpha} \otimes \maL_{\alpha \beta} \to \pi_* \maE_{\beta}$ induce the isomorphisms
\begin{equation*}(\pi_*\varphi_{\alpha \beta})^* \colon \End(\pi_* \maE_{\beta})\to
\End(\pi_*\maE_{\alpha} \otimes \maL_{\alpha \beta})\cong \End(\pi_*\maE_{\alpha})\end{equation*}
over $U_{\alpha \beta}$. These isomorphisms restrict and yield  the isomorphisms
$$
\phi_{\alpha \beta} \, :\, \Psi (\pi^{-1}U_{\alpha \beta}|U_{\alpha \beta}, \maE_\beta) \longrightarrow \Psi (\pi^{-1}U_{\alpha \beta}|U_{\alpha \beta}, \maE_\alpha).
$$
 Since the isomorphisms $(\pi_*\varphi_{\alpha \beta})^*$ satisfy the natural cocycle identity
we have the following:

\begin{lemma}\label{Comp.Local}
The isomorphisms $\phi_{\alpha \beta}$ satisfy
\begin{equation*}\phi_{\alpha\beta}\circ \phi_{\beta\gamma} = \phi_{\alpha\gamma}\quad \text{ whenever }U_{\alpha\beta\gamma}\not = \emptyset.
\end{equation*}
\end{lemma}

\begin{definition}
A fiberwise pseudodifferential operator $P$ of order $\leq m$ with coefficients in the transversally $\maL$-twisted vector bundle $\maE$ is a collection $\{P_\alpha\}_{\alpha\in \Lambda}$,  $P_\alpha \in \Psi^m (\pi^{-1}U_\alpha|U_\alpha ; \maE_\alpha)$ such that
$$
P_\alpha = \phi_{\alpha\beta} (P_\beta).
$$
where $\phi_{\alpha\beta}$ is defined in Equation \ref{phi}.
The space of fiberwise pseudodifferential operators of order $\leq m$, with coefficients in the transversely $\maL$-twisted vector bundle $\maE$, is denoted by $\Psi_\maL^m (M|B ; \maE)$.
\end{definition}
Recall that the isomorphisms $\varphi_{\alpha \beta}$ induce the natural isomorphisms $\varphi_{\alpha \beta}^* \colon \End(\maE_{\beta}) \to \End(\maE_{\alpha})$. Then
\begin{equation}\label{symbolinv}
\sigma_m \circ \phi_{\alpha \beta} = p^*\left(\varphi_{\alpha \beta}^*\right) \circ \sigma_m
\end{equation}

Equation \eqref{symbolinv} implies that if $P=\{P_{\alpha}\} \in \Psi_\maL^m (M|B ; \maE)$ then
the collection $\sigma_m(P_{\alpha})$ defines a section of the (untwisted) bundle $p^*\End(\maE)$. We will call this
section the principal symbol of $P= \{P_{\alpha}\}$ and denote it $\sigma_m(P)$.

We define in the same way the space $\Psi_\maL^m (M|B ; \maE, \maE')$ of fiberwise pseudodifferential operators of order $\leq m$, from the transversally $\maL$-twisted vector bundle $\maE$ to the transversally $\maL$-twisted vector bundle $\maE'$. In particular $\Psi_\maL^m (M|B ; \maE, \maE) = \Psi_\maL^m (M|B ; \maE)$. We also have a principal symbol map
$$
\sigma_m \colon \Psi_\maL^m (M|B ; \maE, \maE ') \longrightarrow C^\infty (\mathring{T}^*(M|B),  p^* \Hom( \maE, \maE')).
$$
We set
$$
\Psi_\maL (M|B ; \maE):= \bigcup_{m\in \Z} \Psi_\maL^m (M|B ; \maE) \text{ and } \Psi_\maL^{-\infty} (M|B ; \maE):= \bigcap_{m\in \Z} \Psi_\maL^m (M|B ; \maE),
$$
and
introduce a composition in $\Psi_\maL (M|B ; \maE)$ by

\begin{equation*}
\{P_\alpha\} \circ \{Q_{\alpha}\} = \{P_{\alpha} Q_{\alpha}\}
\end{equation*}

Since the maps $\phi_{\alpha \beta}$ are algebra isomorphisms the right hand side of this equality defines an element
in $\Psi_\maL (M|B ; \maE)$. Recall that our operators are properly supported.
\begin{proposition}\label{Comp}
The composition of operators is $C^{\infty}(B)$-linear and endows $\Psi_\maL (M|B ; \maE)$  with the structure of an  algebra. Moreover, $\Psi_\maL^{-\infty} (M|B ; \maE)$ is a two-sided  ideal in $\Psi_\maL (M|B ; \maE)$.
\end{proposition}

\begin{remark}\label{Sheafify}
The construction can be sheafified, providing  a sheaf of algebras on $B$, denoted by $\Psi_\maL (\pi; \maE)$, and  given by
$$
 U \longmapsto \Psi_\maL (\pi^{-1}U|U ; \maE), \text{ for open sets }U\subset B.
$$
\end{remark}

We  define as usual the algebra of forms on $B$ with values in $\Psi_\maL(M|B; \maE)$ by
\begin{equation*}
\Omega^*\left(B, \Psi_\maL(M|B; \maE)\right) := \Omega^*(B) \otimes_{C^{\infty}(B)} \Psi_\maL(M|B; \maE).
\end{equation*}
Recall that for every $\alpha$ and $V\subset U_{\alpha}$ we have a fiberwise trace $\Tr_{\alpha} \colon \Psi^{-\infty}(\pi^{-1}V|V; \maE_{\alpha}) \to C^{\infty}(V)$. It is easy to see that for $ P \in
\Psi^{-\infty}(\pi^{-1}U_{\alpha \beta}|U_{\alpha \beta}; \maE_{\beta})$, we have $\Tr_{\alpha} \phi_{\alpha \beta} (P) = \Tr_{\beta} (P)$.
We therefore obtain a well defined map
$$
\Tr \colon \Psi_\maL^{-\infty}(M|B, \maE) \rightarrow C^{\infty}(B) \text{ by setting }\left.\Tr \{P_{\alpha}\}\right|_{U_{\alpha}}= \Tr_{\alpha}(P_{\alpha}).
$$
This trace is a $C^{\infty}(B)$-module map satisfying $\Tr[A, B]=0$.
It extends naturally to define a map
$$
\Tr \colon \Omega^*\left(B, \Psi_\maL^{-\infty}(M|B; \maE)\right) \longrightarrow
\Omega^*\left(B\right).
$$

If the bundle $\maE$ is $\Z_2$ graded we have a similarly defined supertrace
\begin{equation*}\STr \colon \Omega^*\left(B, \Psi_\maL^{-\infty}(M|B; \maE)\right) \longrightarrow
\Omega^*\left(B\right).
\end{equation*}
Note that our definition of $\Psi_\maL (M|B ; \maE)$ depends on the descent datum for $\maE$. It is straightforward however to see that an isomorphism of descent data defines canonically an isomorphism of the corresponding
bundles of algebras. Moreover, every isomorphism $S \colon \maE \to \maE'$ of transversally twisted bundles with the underlying isomorphism $\varsigma$ between the corresponding gerbes $\maL$ and $\maL'$ and with the isomorphisms $\rho_\alpha:\maE_\alpha\otimes S_\alpha\to \maE'_\alpha$, induces an isomorphism of algebras
\begin{equation*}
\Psi(S) \colon \Omega^*\left(B, \Psi_\maL (M | B ; \maE)\right) \longrightarrow \Omega^*\left(B, \Psi_{\maL'} (M | B ; \maE')\right).
\end{equation*}
We also get in the same way an isomorphism of the corresponding sheaves. Indeed, over any intersection $U_{\alpha\beta}$, the maps induced by the line bundles $S_\alpha$ and $S_\beta$ and the isomorphisms $\rho_\alpha$ and $\rho_\beta$, obviously satisfy the relation
$$
\phi'_{\alpha\beta} \circ \Psi (\rho_\beta) = \Psi (\rho_\alpha) \circ \phi_{\alpha\beta}.
$$
They thus define the isomorphism $\Psi(S)$. The following is then an easy consequence of the definitions.

\begin{lemma}\label{Psi(s)}
\item 1] Let $S \colon \maE \to \maE'$, $S' \colon \maE \to \maE'$ be two isomorphisms of transversally twisted bundles. If there exists a $2$-morphism $S \Rightarrow S'$ then $\Psi(S) =\Psi(S')$.

\item 2] Let $S_1 \colon \maE \to \maE'$, $S_2 \colon \maE' \to \maE''$ be two isomorphisms of transversally twisted bundles. Then $\Psi(S_2 \circ S_1) = \Psi(S_2) \circ \Psi(S_1)$.

\item 3] For any $P \in \Omega^*\left(B, \Psi_\maL^{-\infty}(M|B; \maE)\right)$,  $\STr \Psi( S)(P) = \STr P $.
\end{lemma}

\medskip

\subsection{Connections and curvatures}

We continue in the notations of the previous paragraph. We assume that we are given a transversally $\maL$-twisted
bundle $\maE$ with connection represented by a descent datum $(U_{\alpha}, \maL_{\alpha \beta}, \mu_{\alpha \beta \gamma}, \maE_{\alpha}, \varphi_{\alpha \beta})$ with connections $(\nabla_\alpha, \nabla_{\alpha \beta}, \omega_{\alpha})$,   see Subsection \ref{hotwbu}.

Recall the bundles $\pi_*\maE_{\alpha}$ and isomorphisms $\pi_*\varphi_{\alpha \beta}$ defined in Subsection \ref{Projective families}. It is easy to see that
$(U_{\alpha}, \maL_{\alpha \beta}, \mu_{\alpha \beta \gamma}, \pi_*\maE_{\alpha}, \pi_* \varphi_{\alpha \beta})$
is a descent datum for an infinite dimensional twisted bundle.
We now proceed to define a connection on this descent datum.

Choose a horizontal distribution i.e. a subbundle $\maH \subset TM$ such that
$
TM = \maH \oplus T (M|B).
$
This choice together with connections $\nabla^{\maE}_{\alpha}$ defines for each $\alpha$ a connection $\nabla_{\alpha}^{\maH}$ on $\pi_*\maE_\alpha$ as follows:
\begin{equation*}
(\nabla^{\maH}_{\alpha})_X\xi = (\nabla^{\maE}_{\alpha})_{X^{\maH}} \xi
\end{equation*}
where $X^{\maH}$ is the horizontal lift of $X\in \Gamma(B, TB)$.

\begin{lemma}\label{pi*}
$(\pi_*\varphi_{\alpha \beta})^* \nabla^{\maH}_{\beta} = \nabla^{\maH}_{\alpha} \otimes \id+ \id \otimes \nabla_{\alpha \beta}$
\end{lemma}

The curvature of the connection $\nabla^{\maH}_{\alpha}$ is a $2$-form $\theta^{\maH}_{\alpha}$ on $U_{\alpha}$
with values in fiberwise first order differential operators  given by
\begin{equation*}
\theta^{\maH}_{\alpha}(X, Y)=\theta^{\maE}_{\alpha}(X^{\maH}, Y^{\maH})+ (\nabla^{\maE}_{\alpha})_{T^{\maH}(X, Y)}.
\end{equation*}
where
\begin{equation}\label{th}
T^{\maH}(X, Y) = [X^{\maH}, Y^{\maH}]-[X, Y]^{\maH}, \ X, Y \in \Gamma(B, TB).
\end{equation}
Each $\nabla^{\maH}_{\alpha}$ defines a filtration-preserving derivation $\partial^{\maH}_{\alpha}$  of the algebra of fiberwise pseudodifferential operators in the following way:
\begin{equation*}
\partial_{\alpha}^{\maH} \colon \Psi(\pi^{-1}U_{\alpha}|U_{\alpha}, \maE_{\alpha}) \longrightarrow
\Omega^1(U_{\alpha}, \Psi(\pi^{-1}U_{\alpha}|U_{\alpha}, \maE_{\alpha})) \end{equation*}
\begin{equation*}
\hspace{2,2cm} D \longmapsto  \partial^{\maH}_{\alpha} (D) = [\nabla_{\alpha}^{\maH}, D]\,. \end{equation*}

If $D \in \Psi(\pi^{-1}U_{\alpha \beta}|U_{\alpha \beta}, \maE_{\beta})$ then the result of Lemma \ref{pi*}
implies that  \begin{equation*}\partial^\maH_{\alpha} (\phi_{\alpha \beta} (D))= \phi_{\alpha \beta} (\partial^\maH_{\beta}(D)).\end{equation*}
Therefore if $\{D_{\alpha}\}$, $D_{\alpha} \in \Psi_\maL(\pi^{-1}U_{\alpha}|U_{\alpha}, \maE_{\alpha})$,
defines an element in $\Psi_\maL^m(M|B, \maE)$ then
$\{\partial_{\alpha}^{\maH}(D_{\alpha})\}$ defines an element of $\Omega^1(B, \Psi_\maL^m(M|B, \maE))$.
We therefore obtain a  derivation
\begin{equation}\label{defpartial}
\partial^{\maH} \colon \Psi_\maL(M|B, \maE) \longrightarrow  \Omega^1(B, \Psi_\maL(M|B, \maE)),
\end{equation}
which extends to a derivation of the algebra $\Omega^*(B, \Psi_\maL(M|B, \maE))$.

\begin{lemma}\label{pseudotheta}  There exists $\theta^{\maH} \in \Omega^2(B, \Psi_\maL^1(M|B, \maE))$ such that
\begin{equation}\label{deftheta}
\theta^{\maH}|_{U_{\alpha}}= \theta^{\maH}_{\alpha} + \omega_{\alpha}
\end{equation}
\end{lemma}
\begin{proof}
By Lemma \ref{pi*}, $\phi_{\alpha \beta}^* \theta^{\maH}_{\beta} = \theta^{\maH}_{\alpha} +\pi^*(\omega_{\alpha} -\omega_\beta)$, and the statement follows as in Proposition \ref{global curvature}.
\end{proof}
We have
\begin{equation*}
(\partial^{\maH})^2 (D) = [\theta^{\maH}, D]\text{ and }
\partial^{\maH}(\theta^{\maH})= \pi^* \varOmega.
\end{equation*}
where $\varOmega$ is the $3$-curvature form of the connection on $\maL$.

Assume that $S$ is an isomorphism of transversally  twisted bundles $\maE$ and $\maE'$ on $M$ and {{$\varsigma$}} is the underlying isomorphism of the
 underlying gerbes $\maL$, $\maL'$
respectively on $B$. Recall that it induces an isomorphism $\Psi(S) \colon \Psi_{\maL} (M|B, \maE) \to \Psi_{\maL '} (M|B, \maE ')$.

\begin{proposition} The following identities hold:
\begin{enumerate}
\item $(\partial^{\maH})' \circ \Psi(S)= \Psi(S) \circ \partial^{\maH}$.
\item $ \Psi(S) \theta^{\maH} - (\theta^{\maH})' = c(\varsigma).$
\end{enumerate}
\end{proposition}

 \medskip

\section{Gerbes  on groupoids}

\subsection{$\cG$-Gerbes and their morphisms}

Given a gerbe $\maL$ on $\maG^{(0)}$, we get two gerbes $s^*\maL$ and $r^*\maL$ over the manifold $\maG^{(1)}$ by pulling back $\maL$ under the source and range maps $s$ and $r$.

\begin{definition} \label{def:gerbe}
A gerbe $( \maL, \varsigma, \nu)$  on $\maG$, is given by the following data:
\begin{itemize}
\item A gerbe $\maL$ on the unit space $\maG^{(0)}$.
\item An isomorphism $ \varsigma \colon s^*\maL  \rightarrow r^*\maL$ of gerbes on $\maG^{(1)}$.
\item A $2$-morphism  $\nu \colon  \varsigma_{01}^2 \circ  \varsigma_{12}^2 \to  \varsigma_{02}^2$
such that the following  diagram commutes:
\begin{equation*}
\begin{CD}
 \varsigma^3_{01} \circ \varsigma^3_{12} \circ  \varsigma^3_{23} @>\ \ \id\otimes \nu_{012} \ \ >>
\varsigma^3_{02} \circ  \varsigma^3_{23}\\
@V \nu_{123}\otimes\id  VV    @VV\nu_{023}V \\
\varsigma^3_{01} \circ \varsigma^3_{13} @>\ \ \  \nu_{013} \ \ \ >>  \varsigma^3_{03}
\end{CD}
\end{equation*}
\end{itemize}
\end{definition}

In the previous definition, we have used  the notion of  $2$-morphism
of definition \ref{def:2-m-general} and we have denoted $ \varsigma_{ij}^k\colon = (\pr^k_{ij})^*  \varsigma$ and  $\nu_{ij\ell} =\nu^3_{ij\ell} =(\pr^3_{ij\ell})^* \nu$. We shall also refer to a gerbe $(\maL, \varsigma, \nu)$ on $\maG$ as a $\maG$-gerbe.

\begin{definition} A morphism between two gerbes $( \maL, \varsigma, \nu)$,  $( \maL', \varsigma', \nu')$ on $\maG$ is given by the  data
$(u,\varrho)$ where:
\begin{itemize}
\item $u: \maL \rightarrow \maL'$ is a morphism of gerbes over $\maG^{(0)}$.
\item
 $( \varrho, s^*u, r^*u) $ is a  $2$-morphism between
$\varsigma$ and $\varsigma'$ in the sense of Definition \ref{def:2-m-general}, which is compatible with the $2$-morphisms $\nu$ and $\nu'$ over $\cG^{(2)}$.
\end{itemize}

\end{definition}

The compatibility condition with the $2$-morphisms $\nu$ and $\nu'$ over $\cG^{(2)}$ means that the following relation holds (with the obvious notations)
\begin{equation*}
 \varrho^2_{02} \circ (\id \otimes \nu) = (\nu' \otimes \id) \circ (\id \otimes  \varrho^2_{01}) \circ ( \varrho^2_{12} \otimes \id).
\end{equation*}

Notice that  $s^*u: s^*\maL \rightarrow s^*\maL'$ and $r^*u: r^*\maL \rightarrow r^*\maL'$ are the pull-back morphisms and
 that $(s^*u)_{01}= (r^*u)_{12}$, $(r^*u)_{01}=(r^*u)_{02}$ and $(s^*u)_{12}= (s^*u)_{02}$. The compatibility is hence equivalent to  the commutativity of the following diagram

\begin{equation*}
\begin{CD}
 \varsigma'_{01} \circ (r^*u)_{12}\circ  \varsigma_{12} @>\ \ \varrho_{12}\otimes \id \ \ >>
\varsigma'_{01} \circ  \varsigma'_{12}\circ (s^*u)_{12}\\
@V {\id\otimes \varrho_{01}^{-1}} VV    @VV{\id \otimes \nu'}V \\
(r^*u)_{01}\circ \varsigma_{01} \circ \varsigma_{12} @>\ \ \  \varrho_{02}\circ (\nu\otimes\id) \ \ \ >>  \varsigma'_{02}\circ (s^*u)_{02}
\end{CD}
\end{equation*}

\begin{definition} Consider two gerbes $(\maL, \varsigma, \nu)$, $(\maL', \varsigma', \nu')$ on $\maG$ and two isomorphisms $(u, \varrho)$, $(u', \varrho')$ between them. A  $2$-morphism between them is given by a   2-morphism
$
\beta\colon u \rightarrow u',\;$
 such that
\begin{equation*}
\varrho' \circ (s^*\beta \otimes \id) = (\id \otimes r^*\beta) \circ \varrho\,.
\end{equation*}
\end{definition}


We now comment briefly on the case when the gerbe $\maL$ on $\maG^{(0)}$ is trivialized. The morphism $\varsigma$ in this case is a line bundle on $\maG^{(1)}$. The gerbe $L$ therefore is given by the following data (cf. Definition \ref{def:gerbe}):
        \begin{itemize}
        \item A line bundle on $\maG^{(1)}$.
        \item an isomorphism of line bundles  $\nu \colon  \varsigma_{01}^2 \circ  \varsigma_{12}^2 \to  \varsigma_{02}^2$
        such that the following  diagram commutes:
        \begin{equation*}
        \begin{CD}
         \varsigma^3_{01} \circ \varsigma^3_{12} \circ  \varsigma^3_{23} @>\ \ \id\otimes \nu_{012} \ \ >>
        \varsigma^3_{02} \circ  \varsigma^3_{23}\\
        @V \nu_{123}\otimes\id  VV    @VV\nu_{023}V \\
        \varsigma^3_{01} \circ \varsigma^3_{13} @>\ \ \  \nu_{013} \ \ \ >>  \varsigma^3_{03}
        \end{CD}
        \end{equation*}
        \end{itemize}
        One can construct now a convolution algebra of  $\maG$ twisted by $L$.
        As a vector space, it coincides with the compactly supported sections of $\varsigma$, $C^\infty_c(\maG^{(1)}, \varsigma)$.
        The convolution product is given, for $s_1$, $s_2 \in C^\infty_c(\maG^{(1)}, \varsigma)$ by
        \[
        s_1 s_2 = (\pr^2_{02})_* (\nu( (\pr^2_{02})^* s_1 \otimes (\pr^2_{02})^*s_2)),
        \]
        i.e.
        \[
        s_1 s_2 (\gamma) = \sum \limits_{\gamma_1\gamma_2=\gamma} \nu(  s_1(\gamma_1) \otimes s_2(\gamma_2)).
        \]

\subsection{Connections on $\maG$-gerbes} $\;$

A connection on a $\maG$-gerbe $(\maL, \varsigma, \nu)$ is given by
\begin{itemize}
\item A choice of connection  $(\nabla^\maL, \omega^\maL)$ on $\maL$
\item A choice of connection  $\nabla^\varsigma$ on $ \varsigma \colon s^*\maL  \rightarrow r^*\maL$ so that it becomes a  morphism of gerbes with connective structures.
\end{itemize}

Then the constructions given in Subsection \ref{subs:curve}  can be applied to the gerbe $\maL$ on $\maG^{(0)}$,  the morphism $\varsigma$ of gerbes
with connection $r^*\maL$,   $s^*\maL$ on $\maG^{(1)}$  and the $2$-morphism $\nu \colon  \varsigma_{01} \circ  \varsigma_{12} \to  \varsigma_{02}$
on $\maG^{(2)}$ to show the existence of  connections on $\maG$-gerbes and to
give the following differential forms:
\begin{equation*}
\varOmega(\maL) \in \Omega^3(\maG^{(0)}), \quad
c(\varsigma) \in \Omega^2(\maG^{(1)}), \quad
c(\nu) \in \Omega^1(\maG^{(2)})
\end{equation*}

We introduce now the following filtration on $\Omega^{\bullet} (\maG^{(\bullet)})$
\begin{equation*}
F^k \Omega^{\bullet} (\maG^{(\bullet)}) = \bigoplus\limits_{k+i \ge 0} \Omega^{i} (\maG^{(\bullet)}).
\end{equation*}

\begin{lemma}\
The degree $3$ element $(\varOmega(\maL),
c(\varsigma),
c(\nu))$   in $F^{-1}\Omega^{\bullet} (\maG^{(\bullet)})$ is closed.
Its class in the cohomology  of $F^{-1} \Omega^{\bullet} (\maG^{(\bullet)})$ is independent of choices made
(connection on $\maL$ and lift of $\varsigma$).
\end{lemma}

\begin{proof} A straightforward computation gives:
$$
d\varOmega(\maL) = 0,\;
\del_0 \varOmega(\maL) - d c(\varsigma) = 0,\;
 \del_0 c (\varsigma) -\del_1 c (\varsigma) + d c(\nu) = 0\; \text{ and } \; \del c(\nu)=0.
$$
Hence we see that
$
(d+\del) \left(\varOmega(\maL) +
c(\varsigma)+ c(\nu) \right) = 0.$
The second part of the lemma  is a consequence of Proposition \ref{Relations}.
\end{proof}

\begin{definition}\label{defclassgerbe}
Associated with any $\cG$ gerbe, there is a well defined DD class in the cohomology of the subcomplex $F^{-1}\Omega^{\bullet} (\maG^{(\bullet)})$, namely the cohomology class of the $3$-cocycle
\[
\frac{1}{2\pi i}(\varOmega(\maL), c(\varsigma), c(\nu)).
\]
\end{definition}

\begin{remark}
Recall that the complexes $\Omega^{\bullet} (\maG^{(\bullet)})$ and $\Omega^\bullet_{DS}(\cG)$ are quasi-isomorphic.
An explicit formula for the quasi-isomorphism $I_\Delta \colon \Omega^\bullet_{DS}(\cG) \to \Omega^{\bullet} (\maG^{(\bullet)})$    can be found in \cite{dupont}, Theorem 6.4.
\end{remark}

To end this paragraph, let us recall that a  $\maG$-sheaf  is a sheaf on $\maG^{(0)}$ such that for every $\gamma \in \maG$ there is morphism $\cA_{s(\gamma)} \to \cA_{r(\gamma)}$, satisfying the standard action axioms.  Here and as usual $\cA_x$ is the stalk of $\cA$ at $x$. Examples are for instance the $\maG$-sheaves $C^\infty_\maG)$ and $(C^\infty_\maG)^*$ corresponding respectively to smooth and smooth nonvanishing functions.

\begin{remark} \label{Intgroupoidclass}
Similar to the case of gerbes on  manifolds, one associates with any $\maG$  gerbe as above,  its class in $H^2(\maG, (C^\infty_\maG)^*)$. The exponential sequence
\[
0\to \mathbb{Z}\longrightarrow C^\infty_\maG \longrightarrow (C^\infty_\maG)^* \to 0
\]
then gives rise to an exact sequence
\[\ldots \to H^2(\maG, C^\infty_\maG) \to H^2(\maG, (C^\infty_\maG)^*) \to H^3(\maG, \mathbb{Z}) \to H^3(\maG, C^\infty_\maG)\to \ldots
\]
However  the cohomology groups $H^\bullet(\maG, C^\infty_\maG)$ in positive degree do not vanish  in general,
and   the  map $H^2(\maG, (C^\infty_\maG)^*) \to H^3(\maG, \mathbb{Z})$  is therefore not an isomorphism.
It is an isomorphism, however, if $\maG$ is proper, see e.g. Lemma \ref{proper}.
\end{remark}

\begin{remark}
          The morphism of sheaves
         $d \log \colon (C^\infty_\maG)^* \to (\Omega^1_\maG)^{cl}$, where $(\Omega^1_\maG)^{cl}$ is the sheaf of closed $1$-forms, induces a map $H^\bullet( \maG, (C^\infty_\maG)^* ) \to H^\bullet(\maG, (\Omega^1_\maG)^{cl})$.
         Using the resolution $(\Omega^{\ge 1}_\maG, d_{dR})$ of $(\Omega^1_\maG)^{cl}$, we obtain that the latter cohomology  is computed  by the complex $F^{-1}\Omega^{\bullet} (\maG^{(\bullet)})[1]$. We therefore obtain a map
         $H^2( \maG, (C^\infty_\maG)^* ) \to H^3(F^{-1}\Omega^{\bullet} (\maG^{(\bullet)}))$. The class given in Definition
         \ref{defclassgerbe}  is  the image of the  class introduced in  the Remark \ref{Intgroupoidclass} under this map.
         However we will not use this fact, since we only use the class from Corollary \ref{defclassgerbe}.
         \end{remark}

\medskip

\subsection{DD class in the Dupont-Sullivan bicomplex}\label{DS-DDclass}

We now  describe explicitely a representative (a  $3$-cocycle $\varOmega (L)$  in the Dupont-Sullivan bicomplex) of the canonical DD class associated with the $\cG$ gerbe $L$ satisfying
$$
I_\Delta\left(\left[\frac{\varOmega(L)}{2i\pi}\right]\right)\; = \;  \left[\frac{1}{2i\pi} \, (\varOmega(\maL), c(\varsigma), c(\nu))\right].
$$
The notation should be clear and no confusion should occur between this $3$-cocycle $\varOmega(L)$ and the differential $3$-form $\varOmega (\maL)$ associated with the gerbe $\maL$ on $\cG^{(0)}$. The components of $\varOmega(L)$ are given in Lemma \ref{DS-Curvature} below.

Assume then that $L=(\maL, \varsigma, \nu)$ is a gerbe with connection on $\maG$.
Let $\cL^{(k)}$  be the pullback of $\cL$  to   $\maG^{(k)}$ via the map $\pr^k_0 \colon (\gamma_1, \ldots, \gamma_k)\mapsto r(\gamma_1)$.
There are isomorphisms
\begin{equation*}(\pr^k_{0i})^* \varsigma  \colon (\pr^k_i)^* \cL \to \cL^{(k)}
\end{equation*}
of gerbes with connective structures, obtained by using the pulled-back connections.
Since $\pr^{k-1}_0 \circ \del_i = \pr^k_0$, for $i\ge 1$, and $\pr^{k-1}_0 \circ \del_0 = \pr^k_1$,
we have canonical isomorphisms of gerbes with connections
\begin{equation}\label{gsimplicial1}
\del_i^*\cL^{(k-1)} \cong  \cL^{(k)}, \quad \text{for } i \ge 1,
\end{equation}
 and  a morphism of gerbes
\begin{equation}\label{gsimplicial2}
(\pr_{01}^k)^* \varsigma \colon \del_0^*\cL^{(k-1)} \to  \cL^{(k)}, \quad \text{for } i =0.
\end{equation}
Both of the gerbes in the equation \eqref{gsimplicial2} are equipped with connections, and
$(\pr_{01}^k)^* \varsigma$ becomes a morphism of gerbes with connective structures when equipped with any of the
connections $\nabla_i$ determined by the condition that the $2$-morphism
\begin{equation*}(\pr_{01i}^k)^* \nu \colon
(\pr_{01}^k)^* \varsigma \circ (\pr_{1 i}^k)^* \varsigma \to (\pr_{0i}^k)^* \varsigma
\end{equation*}
is flat. Here the left hand side is equipped with the connection $\nabla_i \otimes (\pr_{1 i}^k)^*\nabla^\varsigma$ and the right hand side is equipped with the connection
$(\pr_{0 i}^k)^*\nabla^\varsigma$. For simplicity, we omit from the notations the pullback to the simplex $\Delta^k$.

 We have
\begin{equation*}
\nabla_1 = (\pr_{0 1}^k)^*\nabla^\varsigma \text{ and } \nabla_i- (\pr_{0 1}^k)^*\nabla^\varsigma = c((\pr_{01i}^k)^* \nu).
\end{equation*}
In addition to its curving $\omega_0:=(\pr^k_0)^* \omega^\cL$,  the gerbe $\cL^{(k)}$ has curvings
\begin{equation*}
\omega_i:=((\pr^k_{0i})^* \varsigma^{-1})^*(\pr^k_i)^* \omega^\cL\; \text{ for any } i=1, \ldots, k\, \end{equation*}
which are also compatible with the connective structure, and we have
\begin{equation*}
\omega_i-\omega_j= c((\pr_{ij}^k)^*\varsigma).
\end{equation*}

Let $\cL^{(k)}_\Delta$ denote the pull-back of $\cL^{(k)}$ to $\cG^{(k)} \times \Delta^k$ via the projection on the first factor.
A direct calculation shows that we have for any $i\ge 1$ flat isomorphisms of gerbes with connections
\begin{equation}\label{gsimp1}
(\del_i \times \id)^*  \cL^{(k-1)}_\Delta \cong (\id \times \delta_i)^*  \cL^{(k)}_\Delta
\end{equation}
induced by \eqref{gsimplicial1}.

For $i=0$ we have an isomorphism of gerbes
\begin{equation}\label{gsimp2}
(\pr_{01}^k)^* \varsigma \colon (\del_0 \times \id)^* \cL^{(k-1)}_\Delta \to (\id \times \delta_0)^* \cL^{(k)}_\Delta
\end{equation}

Equip now $\maL^{(k)}_\Delta$ (for every $k$) with connection $(\nabla_\Delta, \omega_\Delta)$ where
\begin{equation}\label{gconnective_structure} \nabla_\Delta = (\pr^k_0)^* \nabla^\maL, \text{ and }
\end{equation}
\begin{multline}\label{gcurving}
\omega_\Delta := \sum_{i=0}^k t_i \omega_i - \sum \limits_{0<i <j \le k} (t_i dt_j -t_j dt_i) (\pr_{0,i, j}^k)^* c(\nu)=\\
\omega_0 - \sum \limits_{i=1}^{k} t_i (\pr_{0, i}^k)^* c(\varsigma) - \sum \limits_{0<i <j \le k} (t_i dt_j -t_j dt_i) (\pr_{0,i, j}^k)^* c(\nu)
\end{multline}
Equip finally the isomorphism \eqref{gsimp2} with the connection
\begin{equation}\label{gsimplicial connection}
\nabla^\varsigma_\Delta := \sum_{i=1}^k t_i \nabla_i= (\pr_{0 1}^k)^*\nabla^\varsigma +\sum_{i=2}^k t_ic((\pr_{01i}^k)^* \nu)
\end{equation}

\begin{lemma} With the choice of connections and curvings given in \eqref{gconnective_structure}, \eqref{gcurving} and \eqref{gsimplicial connection}, the maps  given by \eqref{gsimp1} and \eqref{gsimp2} are flat isomorphisms of gerbes with connections.
\end{lemma}

\begin{proof}
Recall that by the definition of flatness (cf. Lemma \ref{curvature1morphism}) we need to
verify that the curvatures   of the isomorphisms in \eqref{gsimp1} (for $i\ge 1$) and \eqref{gsimp2} (for $i=0$) vanish. The proof for $i \ge 1$ is straightforward; we verify the statement for $i=0$. It is clear that \eqref{gsimp2} is a morphism of
gerbes with connective structures, so we just need to verify its flatness. If we equip the left and right hand sides with curvings $\omega_1$ and $\omega_0$ respectively
and \eqref{gsimp2} with the connection $(\pr_{0 1}^k)^*\nabla^\varsigma$, the curvature of \eqref{gsimp2} is $(\pr_{0, i}^k)^* c(\varsigma)$. It follows that
if we use curvings given by \eqref{gcurving} and connection given by \eqref{gsimplicial connection} the curvature of \eqref{gsimp2} will be (we denote coordinates on $\Delta^{k-1}$ by $t_1$, $t_2$, \ldots, $t_k$ with  $t_1+t_2+\ldots t_k=1$)
\begin{multline}\label{longsum}
(\pr_{0, 1}^k)^* c(\varsigma)- \sum \limits_{i=1}^{k} t_i (\pr_{0, i}^k)^* c(\varsigma) - \sum \limits_{0<i <j \le k} (t_i dt_j -t_j dt_i) (\pr_{0,i, j}^k)^* c(\nu)
+\\ \sum \limits_{i=1}^{k} t_i (\pr_{1, i}^k)^* c(\varsigma) +  \sum \limits_{1<i <j \le k} (t_i dt_j -t_j dt_i) (\pr_{1,i, j}^k)^* c(\nu)+ d\left( \sum_{i=2}^k t_ic((\pr_{01i}^k)^* \nu) \right)=\\ \sum_{i=2}^k t_i\left((\pr_{0, 1}^k)^* c(\varsigma)+(\pr_{1, i}^k)^* c(\varsigma)-(\pr_{0, i}^k)^* c(\varsigma)\right)+\\
\sum \limits_{0<i <j \le k} (t_i dt_j -t_j dt_i) \left((\pr_{1,i, j}^k)^* c(\nu)- (\pr_{0,i, j}^k)^* c(\nu)\right)+ d\left( \sum_{i=2}^k t_ic((\pr_{01i}^k)^* \nu) \right).
\end{multline}
Using the identities
\begin{align*}
&(\pr_{0, 1}^k)^* c(\varsigma)+(\pr_{1, i}^k)^* c(\varsigma)-(\pr_{0, i}^k)^* c(\varsigma)=- d(\pr_{0,1, i}^k)^* c(\nu)\\
&(\pr_{1,i, j}^k)^* c(\nu)- (\pr_{0,i, j}^k)^* c(\nu)=-(\pr_{0,1, j}^k)^* c(\nu)+ (\pr_{0,1, i}^k)^* c(\nu)
\end{align*}
we rewrite the  sum \eqref{longsum} as
\begin{equation}\label{ressum}
 \sum \limits_{0<i <j \le k} (t_i dt_j -t_j dt_i) \left(-(\pr_{0,1, j}^k)^* c(\nu)+ (\pr_{0,1, i}^k)^* c(\nu)\right)+  \sum_{i=2}^k dt_i\,c((\pr_{01i}^k)^* \nu).
\end{equation}
But notice that
\begin{multline*}
 \sum \limits_{0<i <j \le k} (t_i dt_j -t_j dt_i) \left(-(\pr_{0,1, j}^k)^* c(\nu)+ (\pr_{0,1, i}^k)^* c(\nu)\right) =\\
  \sum \limits_{0<i,j \le k} (t_i dt_j -t_j dt_i)  (\pr_{0,1, i}^k)^* c(\nu)=-\sum \limits_{0<i \le k}dt_i\,(\pr_{0,1, i}^k)^* c(\nu).
\end{multline*}
Whence we conclude that the sum in the \eqref{ressum} (and thus the sum in \eqref{longsum}) is equal to $0$.
This verifies the flatness of the morphism \eqref{gsimp2}.
\end{proof}

\begin{lemma}\label{DS-Curvature} With the choice of connection as above
\begin{equation*}
\varOmega(\maL^{(k)}_\Delta) =  \sum_{i=0}^kt_i(\pr^k_{i})^*\varOmega(\cL)-
\sum_{0<i<j\le k} (t_idt_j-t_j dt_i) (\pr^k_{ij})^*c(\varsigma)-2 \sum_{0<i<j\le k} dt_idt_j (\pr^k_{0 i j})^* c(\nu).
\end{equation*}
\end{lemma}
\begin{proof}\ We use the identities
\begin{align*}
&  c((\pr_{0i}^k)^*\varsigma)=\omega_0-\omega_i\\
& d\,(\pr_{0,i, j}^k)^* c(\nu)=-(\pr_{0, i}^k)^* c(\varsigma)-(\pr_{i, j}^k)^* c(\varsigma)+(\pr_{0, j}^k)^* c(\varsigma)
\end{align*}
to compute
\begin{multline*}
\varOmega(\maL^{(k)}_\Delta) =\sum_{i=0}^kt_i(\pr^k_{i})^*\varOmega(\cL)  -\sum_{i=1}^kdt_i(\pr_{0, i}^k)^* c(\varsigma)- 2 \sum_{0<i<j\le k} dt_idt_j (\pr_{0, i, j}^k)^* c(\nu)
+\\ \sum_{0<i<j\le k} (t_idt_j-t_j dt_i) (-(\pr_{0, i}^k)^* c(\varsigma)-(\pr_{i, j}^k)^* c(\varsigma)+(\pr_{0, j}^k)^* c(\varsigma)).
\end{multline*}
Now, the statement follows since
\begin{multline*}
\sum_{0<i<j\le k} (t_idt_j-t_j dt_i) (-(\pr_{0, i}^k)^* c(\varsigma) +(\pr_{0, j}^k)^* c(\varsigma))=\\
\sum_{0<i,j\le k} (t_idt_j-t_j dt_i)  (\pr_{0, j}^k)^* c(\varsigma)= \sum_{0<j\le k} dt_j\,  (\pr_{0, j}^k)^* c(\varsigma).
\end{multline*}
\end{proof}

\begin{definition}
The collection of forms $\varOmega(\maL^{(k)}_\Delta)$ defines the element  $\varOmega(L) \in\Omega^3_{DS}(\cG)$, with vanishing $\Omega^{0,3}_{DS}(\cG)$ component.
\end{definition}

We therefore can form, as in Subsection \ref{defcurrent}, Equation \eqref{defdtsimp}, a twisted differential $d^t_{\varOmega(L)}$ and consider the twisted cohomology $H_{\varOmega (L)} ^\bullet (B\cG)$ of the classifying space $B\cG$.

\section{Projective $\cG$-families and the analytic index}\label{projective}

\subsection{Twisted vector bundles on groupoids}

Assume that   $L=(\maL, \varsigma, \nu)$ is a $\maG$-gerbe as before.
\begin{definition}\label{TwistedBundles} An $L$-twisted vector bundle $E=(\maE, \phi)$ is given by the following data:
\begin{enumerate}
\item An $\maL$-twisted vector bundle $\maE$ on $\maG^{(0)}$.
\item An isomorphism of twisted vector bundles $\phi \colon s^* \maE \to r^* \maE$ such that:
\begin{itemize}
\item The underlying
isomorphism of gerbes is $ \varsigma \colon s^*\maL  \rightarrow r^*\maL$.
\item   The $2$-morphism $\nu \colon  \varsigma_{01}^2 \circ  \varsigma_{12}^2 \to  \varsigma_{02}^2$ defines a $2$-morphism
$\nu \colon  \phi^2_{01} \circ  \phi^2_{12} \to  \phi^2_{02}$.
\end{itemize}
 \end{enumerate}
\end{definition}
We have used Definition \ref{morphisms for bundles} and denoted as before  $ \phi^2_{ij} := (\pi^2_{ij})^*  \phi$.
Assume that $E=(\maE, \phi)$ is an $L$-twisted vector bundle and let  $\nabla^{\cE}$ be a connection on $\cE$.
As for the gerbe, we denote by $\cE^{(k)} $  the pullback of $\cE$  to  $\cG^{(k)}$ via the map $\pr^k_0$.

Note also that there is an isomorphism $(\pr^k_{0i})^* \phi  \colon (\pr^k_i)^* \cE \to \cE^{(k)}$ and
similarly to \eqref{gsimplicial1},  \eqref{gsimplicial2}
we as well have canonical isomorphisms
\begin{equation}\label{simplicial1}
\del_i^*\cE^{(k-1)} \cong  \cE^{(k)}, \text{ for } i \ge 1,
\end{equation}
and  an isomorphism
\begin{equation}\label{simplicial2}
(\pr_{01}^k)^* \phi \colon \del_0^*\cE^{(k-1)} \to  \cE^{(k)}, \text{ for }i=0,
\end{equation}
with the underlying isomorphisms of gerbes \eqref{gsimplicial1},  \eqref{gsimplicial2} respectively.

The choice of the connection $\nabla^{\cE}$ induces connections $\nabla^{\cE^{(k)}}_i$ on $\cE^{(k)}$, $i=0,1,\ldots, k$ by
the following formula
\begin{equation*}
\nabla^{\cE^{(k)}}_i = ((\pr^k_{0i})^* \phi^{-1})^*(\pr^k_i)^* \nabla^{\cE}
\end{equation*}
Finally, we pull back   $\maE^{(k)}$   to $\cG^{(k)}\times \Delta^k$
via the projections on the first factor obtaining the bundle $\maE^{(k)}_\Delta$, etc.

Note that we have a canonical isomorphism of twisted bundles
\begin{equation}\label{simp1}
(\del_i \times \id)^*  \cE^{(k-1)}_\Delta \cong (\id \times \delta_i)^*  \cE^{(k)}_\Delta, \text{ for } i\ge 1
\end{equation}
with the underlying isomorphism of gerbes \eqref{gsimp1}.

For $i=0$ we have an isomorphism of twisted bundles
\begin{equation}\label{simp2}
(\pr_{01}^k)^* \phi \colon (\del_0 \times \id)^* \cE^{(k-1)}_\Delta \longrightarrow (\id \times \delta_0)^* \cE^{(k)}_\Delta
\end{equation}
with the underlying isomorphism of gerbes \eqref{gsimp2}.

On the bundle $\maE^{(k)}_\Delta$, we now consider the connection
\begin{equation}\label{simplicial connection}
(\nabla^{\cE})^{(k)} =\sum \limits_{i=0}^{k} t_i  \nabla^{\cE^{(k)}}_i +d_{dR}
\end{equation}
where $d_{dR}$ is the de Rham differential in $t$ variables. Direct calculation shows the following:

\begin{lemma}\label{lemma connection invariance} Equip $\maE^{(k)}_\Delta$ with the connection \eqref{simplicial connection}, and $(\pr_{01}^k)^* \varsigma $ with the connection
\eqref{gsimplicial connection}. Then \eqref{simp1} and \eqref{simp2} become isomorphisms of bundles with connections.
\end{lemma}

Let $(\theta^{\cE})^{(k)}$ be the curvature of $(\nabla^{\cE})^{(k)}$. Use the isomorphisms \eqref{simp1} and \eqref{simp2} to identify
$\End ((\del_i \times \id)^*  \cE^{(k-1)}_\Delta)$  with  $\End((\id \times \delta_i)^*  \cE^{(k)}_\Delta)$. Then Lemma \ref{lemma connection invariance} implies that   we have
\begin{equation*}
(\del_i \times \id)^* (\theta^{\cE})^{(k-1)} =
(\id \times \delta_i)^*  (\theta^{\cE})^{(k)}, \text{ for } i \ge 0.
\end{equation*}

\medskip

Now given a $\maG$-equivariant submersion $P \to \maG^{(0)}$ as before, one can similarly define a notion of transversally twisted $L$-bundle. Recall the \'etale groupoid $\tP=P\rtimes \cG$ with source and range maps also denoted $s$ and $r$ respectively. The submersion  $P \to \maG^{(0)}$ gives a submersion $P\rtimes \cG\to \cG$ which is the pull-back under the range map and which is also a groupoid morphism. Moreover, given a gerbe $\maL$ on $\cG^{(0)}$ and a transversely $\maL$-twisted bundle $\maE$ on $P\to \maG^{(0)}$, the bundle $r^*\maE$ is a transversally $r^*\maL$-twisted vector bundle over $P\rtimes \maG$ with respect to the second projection $\what{\pi}: P\rtimes \maG \to \maG$ and to the gerbe $r^*\maL$ on $\maG$.

\begin{definition} An transversally $L$-twisted vector bundle is given by the following data:
\begin{enumerate}
\item An $\maL$-transversely twisted bundle $\maE$ on $P\to \maG^{(0)}$.
\item An isomorphism of transversally twisted vector bundles $\phi \colon s^* \maE \to r^* \maE$ over the submersion $P\rtimes \cG\to \cG$ such that:
\begin{itemize}
\item The underlying
isomorphism of gerbes is $ \varsigma \colon s^*\maL  \rightarrow r^*\maL$.
\item   The $2$-morphism $\nu \colon  \varsigma_{01}^2 \circ  \varsigma_{12}^2 \to  \varsigma_{02}^2$ defines a $2$-morphism
$\nu \colon  \phi^2_{01} \circ  \phi^2_{12} \to  \phi^2_{02}$
where $ \phi^2_{ij}\colon = (\pi^2_{ij})^*  \phi$.
\end{itemize}
\end{enumerate}
\end{definition}

\subsection{Projective $\maG$-operators} $\;$

Assume as before that $\pi: P \to \maG^{(0)}$ is a $\maG$-equivariant submersion such that $\maG$ acts freely and properly  on $P$, and let $L=(\maL, \varsigma, \nu)$ be a $\maG$-gerbe. We assume for simplicity, and in view of our interest in the index theorem, that the action of $\maG$ on $P$ is cocompact.
A  $\maG$-invariant distance on the $\maG$-manifold $P$ will be any smooth  function $d=d_P$ on the closed submanifold $\{(p, p')\in P^2, \pi (p)= \pi (p')\}$ of $P^2$ such that
\begin{enumerate}
\item For any  fiber  $P_x:=\pi^{-1} (x)$ of $\pi$, the restriction $d_x$ of $d$ to $P_x \times P_x$  is a distance which defines the induced topology of $P_x$.
\item For any $\gamma\in \maG$ and any $(p,p')\in P_{r(\gamma)}^2$, we have $d_{s(\gamma)} ( p\gamma, p'\gamma) = d_{r(\gamma)} (p, p')$.
\end{enumerate}
It is easy to check using the properness of the action that such distance always exists.  We fix such $d$ and may use it in the definition of the proper support.

Let $\maE$ be a transversally $L$-twisted vector bundle on $P$.
Recall the vector space  $\Psi_\maL ( P | \maG^{(0)} ; \maE)$  of
 properly supported $\maL$-projective pseudo-differential operators on the submersion $P\to \maG^{(0)}$. Then,   $\Psi_\maL ( P | \maG^{(0)} ; \maE)$  defines a $\maG$-sheaf of algebras, and we describe the functorial isomorphisms now. First the similar space of (uniformly in $\maG$) properly supported pseudodifferential operators, which are projective for the gerbe $r^*\maL$, will suggestively be denoted $r^*\Psi_\maL (P | \maG^{(0)} ; \maE)$, so
$$
r^*\Psi_\maL (P | \maG^{(0)} ; \maE) := \Psi_{r^*\maL} (P\rtimes \maG | \maG ; r^*\maE) \text{ with again } r (p, \gamma)=p.
$$
We define similarly
$$
s^*\Psi_\maL (P | \maG^{(0)} ; \maE) := \Psi_{s^*\maL} (P\rtimes \maG | \maG ; s^*\maE) \text{ and } s (p, \gamma)=p\gamma.
$$
Notice that $s^*\maE$ is now a transversally twisted vector bundle on $P\rtimes \maG$ with respect to the gerbe $s^*\maL$ on $\maG$. The isomorphism $\varsigma: s^*\maL\to r^*\maL$ can then be pulled back to $P\rtimes \maG$ using the second projection $\what{\pi}:P\rtimes \maG \to \maG$ and yields an isomorphism of transversally twisted bundles
$$
\what{\pi}^*\varsigma: s^*\maE \longrightarrow r^*\maE.
$$
We get in this way the functorial  isomorphisms
of sheaves of algebras  by considering
$$
\Psi(\varsigma) \colon s^*\Psi_\maL (P | \maG^{(0)} ; \maE) \longrightarrow r^*\Psi_\maL (P | \maG^{(0)} ; \maE),
$$
 see Lemma \ref{Psi(s)}.

Given an operator $A \in \Psi_\maL ( P | \maG^{(0)} ; \maE)$, we define the operators $s^*A$ and $r^*A$ in $s^*\Psi_\maL (P | \maG^{(0)} ; \maE)$ and $r^*\Psi_\maL (P | \maG^{(0)} ; \maE)$ respectively, as follows. Recall that $s$ and $r$ are also local diffeomorphisms on $P\rtimes \cG$ such that we may take as domains  open sets  of the form $\pi^{-1} U$ for $U\subset \maG^{(0)}$. Therefore, we can use the following commutative diagrams and the usual pull-back of local vertical pseudodifferential operators under fiberwise preserving diffeomorphisms
\begin{equation*}
\begin{CD}
P\rtimes \maG
@>{\what{r}}>>
P\\
@V{\hat\pi}VV  @VV{\pi}V\\
\maG @>\ \ {r} \ \ >>\maG^{(0)}
\end{CD}
\; \text{ and } \;
\begin{CD}
P\rtimes \maG
@>{\what{s}}>>
P\\
@V{\hat\pi}VV  @VV{\pi}V\\
\maG @>\ \ {s} \ \ >>\maG^{(0)}
\end{CD}
\end{equation*}
Indeed, it is easy to check that we obtain in this way transversely twisted pseudodifferential operators $r^*A$ and $s^*A$ with respect to $r^*\maL$ and $s^*\maL$, which are, uniformly in $\maG$, properly supported.

\begin{definition} \label{def:4.4}
The vector space of $\maG$-invariant (fiberwise) pseudo-differential operators
$\Psi_\maL ( P | \maG^{(0)} ; \maE)^\maG $ is defined to be $\{ A \in \Psi_\maL ( P | \maG^{(0)} ; \maE)\ | \  \Psi(\varsigma)(s^* A) = r^* A \}$.
\end{definition}

\medskip We pause here to give an explicit example.

\subsection{An example of twisted $\maG$-family} \label{sub:ex} $\;$

Let $B$ be a smooth manifold on which the group $\Z$ acts. We denote the (action of) generator of $\Z$ by $T$.
Let $\maG$ be the corresponding action groupoid: $\maG^{(0)}= B$, $\maG^{(1)} = B \times \Z $ with $s(b, n) =b$, $r(b, n) = T^nb$, $(b, n ) \in B \times \Z$.  Let $\ell$ be a $\Z$ equivariant line bundle on $B$. We can construct  a $\maG$-gerbe $L$ as follows.
\begin{itemize}
\item The gerbe $\maL$ on the unit space $\maG^{(0)}$ is trivial.
\item The isomorphism $ \varsigma \colon s^*\maL  \rightarrow r^*\maL$ of gerbes on $\maG^{(1)}$ is given by $\varsigma|_{B\times n} = \ell^{\otimes n}$, $n \in \Z$.
\item A $2$-morphism  $\nu \colon  \varsigma_{01}^2 \circ  \varsigma_{12}^2 \to  \varsigma_{02}^2$ is given by the isomorphism
$(T^n)^*\ell^{\otimes m} \otimes \ell^{\otimes n} \to \ell^{\otimes(m+n)}$ on $B \times m \times n$, induced by the isomorphism $(T^n)^*\ell^{\otimes m} \to \ell^{\otimes m}$.

\end{itemize}

Note that we have canonical isomorphism $H^\bullet(\maG; \C) \cong H^\bullet(B_T; \C)$ where $B_T$ is the mapping torus of the map $T$: $B_T:= (B \times [0,1])/\sim$. Here $\sim$ is the equivalence relation given by $(b, 1) \sim (Tb, 0)$. Under this identification the de Rham class of $L$ can be described by the following $3$-form $\varOmega \in \Omega^3( B_T)$. Let $\phi \in C^\infty_c(0, 1)$ be such that $\int_0^1 \phi (t) dt =1$. Chose a connection $\nabla$ on $\ell$ and let $\theta \in \Omega^2(B)$ be its curvature. Then the form $\theta \phi(t) dt \in \Omega^3(B \times [0,1])$ descends to a well defined closed $3$ form $\varOmega$ on $B_T$. It is easy to see that the cohomology class of this form is, in general, nontrivial.

Let $\widehat{B}$ be the circle bundle of $\ell$, $\widehat{B} = \{\ell \setminus 0\}/\R^{>0}$, and let $p \colon \widehat{B} \to B$ be the natural projection. The $\Z$ action on $B$ naturally induces a $\Z$ action on $\widehat{B}$, since $\ell$ is equivariant. Choose an Hermitian metric on $\ell$. This choice induces canonically a trivialization of $p^* \ell$, and let $s \in \Gamma(\widehat{B}; p^* \ell)$ be the corresponding section. If we define the form $\alpha \in \Omega^1(\widehat{B})$ (depending on $s$) by
$(p^*\nabla)(s) = \alpha s$, then
\[
d \alpha = p^*\theta.
\]
It follows that with $\eta = \alpha \phi(t) dt \in \Omega^2(\widehat{B}_T)$  we have $d \eta = p^*\varOmega$. Here $\widehat{B}_T$ is the mapping torus of $\Z$ action on $\widehat{B}$. It follows that  the map
\[I_\eta \colon \xi \mapsto e^{- u \eta} \xi\] intertwines de Rham differential $ud$ with $d_\varOmega$ and induces an isomorphism between the twisted and untwisted cohomology of $\widehat{B}_T$.

Assume that $M$ is a compact spin manifold and $\widetilde{M} \to M$ is a $\Z$-covering. Let $\omega$ be an $S^1$-equivariant line bundle on $M \times S^1$.
denote by  $D$  a $\Z \times S^1$ equivariant Dirac-type operator on $\widetilde{M} \times S^1$, acting on sections of the spinor bundle twisted by $\omega$.
 Set $P:= \widetilde{M}\times \widehat{B}$; then $\pi \colon P \to B$ given by $\pi(m, b)=p(b)$, $m \in \widetilde{M}, b \in \widehat{B}$, is naturally a $\maG$-fibration. We will denote by $\mathcal{E}$ for the pull back of the spinor bundle tensored with $\omega$ to $P$. It is a $\Z$ equivariant bundle on $P$. It can be endowed with a structure of $L$- transversally twisted bundle. To do this we need to specify isomorphism of twisted vector bundles $\phi \colon s^* \maE \to r^* \maE$. In this example is given by isomorphism
\[
\pi^* \ell^{\otimes n} \otimes \maE \to (T^n)^* \maE
\]
given by the trivialization of $\pi^* \ell$ (induced by the trivialization of $p^* \ell$) and the isomorphism $(T^n)^* \maE \cong \maE$ (given by $\Z$-equivariance of $\maE$). The Dirac operator $D$ induces a constant family on the fibration $ (\widetilde{M} \times S^1) \times \widehat{B} \to \widehat{B}$ which descends, since $D$ is $S^1$-equivariant, to the family on $P \to B$ acting on the sections of $\maE$. It is easily seen to be $\maG$-invariant.

\bigskip

\subsection{Invariant operators and cross-products} $\;$

\medskip

For any c-soft $\maG$ sheaf of algebras $\cA$, one can define the convolution algebra $\cA \rtimes \maG$ as follows, see \cite{crainic}, 3.31 for the details. As a vector space it is the space of compactly supported sections $\Gamma_c(\maG; \cA)$.
The convolution product is given by
\[
fg (\gamma) = \sum_{\gamma_1 \gamma_2 = \gamma} f(\gamma_1) (\gamma_1 \cdot g(\gamma_2)).
\]
Now the sheaf of algebras $\Psi_\maL ( \pi ; \maE)$ associated with $\pi:P\to \cG^{(0)}$ and defined in Remark \ref{Sheafify} being a c-soft $\cG$-sheaf, we eventually obtain the well-defined cross-product algebra that we denote for simplicity by $\Psi_\maL ( P | \maG^{(0)} ; \maE) \rtimes \maG$.
As a space,\
$\Psi_\maL ( P | \maG^{(0)} ; \maE) \rtimes \maG$  can be identified with  the space of properly supported operators from  $\Psi_{r^*\maL} (P\rtimes \maG | \maG ; r^*\maE)$  such that the projection of the support in $\maG$ is compact, and
with the convolution product defined above.

Consider now $\phi \in C_c^\infty(P ; \mathbb{R}) $ such that for any $p\in P$,
\[\sum_{\{ \gamma \mid r(\gamma)=\pi(p)\}} \phi^2( p  \gamma) =1.\]
Let $\gamma \in \maG$. Define $\phi^\gamma \in C^\infty(\pi^{-1}(r(\gamma)))$ by $\phi^\gamma (p) := \phi(p \gamma), p \in r(\gamma)$.
Notice that $\phi^\gamma$ is a restriction to $\pi^{-1}(r(\gamma))$ of a smooth function
defined on $\pi^{-1}(U)$ for a sufficiently small neighborhood $U$ of $r(\gamma)$.
Indeed, let $V$ be a neighborhood of $\gamma$ in $\maG$ such that $r$ and $s$ restricted to
$V$ are diffeomorphisms. Then we can take $U=r(V)$. Notice that the action of $\gamma \colon
\pi^{-1}(r(\gamma)) \to \pi^{-1}(s(\gamma)) $ extends to a diffeomorphism, which we also denote by $\gamma \colon \pi^{-1}(r(V)) \to \pi^{-1}(s(V))$, and this defines the extension of $\phi^\gamma$ to $\pi^{-1}(U)$. Note also that $\phi^\gamma(p)$ is smooth
in $(\gamma, p) \in V \times \pi^{-1}(V)$, $\pi(p) =r(\gamma)$.
The straightforward verification  of the following lemma is left to the reader.

\begin{lemma}\label{mapI} The map
\begin{equation*}
I \colon \Psi_\maL ( P | \maG^{(0)} ; \maE)^\maG \rightarrow \Psi_\maL ( P | \maG^{(0)} ; \maE) \rtimes \maG\text{ defined by }I(A)(\gamma) =   \phi A \phi^\gamma,
\end{equation*}
is an algebra homomorphism.
\end{lemma}

\subsection{Analytic index of projective $\cG$-operators}

\begin{definition}
Let $m\geq 0$ be fixed. A fiberwise pseudodifferential operator
$$A \in \Psi_\maL^m (P | \maG^{(0)} ; \maE, \maE')
$$
is fiberwise elliptic if the principal symbol $\sigma_m(A) \in \Gamma(p^*\Hom(\maE, \maE'))$ is an isomorphism.
\end{definition}

We say that $Q \in \Psi_\maL^{-m} (P | \maG^{(0)} ; \maE', \maE)$ is a parametrix of $A\in \Psi_\maL^m (P | \maG^{(0)} ; \maE, \maE')$ if
$A\circ Q-1 \in \Psi_\maL^{-\infty} (P | \maG^{(0)} ; \maE', \maE')$ and $Q\circ A-1 \in \Psi_\maL^{-\infty} (P | \maG^{(0)} ; \maE, \maE)$. Recall that the action of $\maG$ on $P$ is free proper and co-compact.

\begin{lemma} Every $\cG$-invariant elliptic operator $A\in \Psi_\maL^m (P | \maG^{(0)} ; \maE, \maE')^\maG$ has a $\maG$-invariant parametrix $Q$ in $ \Psi_\maL^{-m} (P | \maG^{(0)} ; \maE', \maE)^\maG$.
\end{lemma}
\begin{proof} Consider the algebra of $\cG$-invariant transversally $L$-twisted fiberwise pseudo-differential operators
$\Psi_\maL^m (P | \maG^{(0)} ; \maE \oplus \maE', \maE \oplus \maE')^\maG\,,$
 and the element
$$\mathcal{A}= \begin{pmatrix} 0 &A^* \\ A &0\end{pmatrix}\,.$$
where $A^*$ is a formal adjoint of $A$, an element of $\Psi_\maL^m (P | \maG^{(0)} ; \maE', \maE)^\maG$, defined using some unitary structures on $\maL$ and $\maE\oplus \maE'$.


Consider a $\maG$-invariant distance  $d$
on $P$ as before. There exists a $\maG$-invariant function
$k: P \times P \rightarrow \mathbb{R}$ and $a,b \in \mathbb{R}^{+*}$
such that $k(x,y)=1$ for $d(x,y) < a$  and $k(x,y)=0$ for $d(x,y) >b$.
Consider the function of the real variable $f(u)= \frac{1-e^{-u^2}}{u}$.
One then checks that the operator with Schwartz kernel
$k(x,y) K_{f(\mathcal{A})}(x,y)$  defines a $\maG$-invariant parametrix for $\mathcal{A}$. One then immediately obtains a $\maG$-invariant parametrix (uniformly properly supported) for $A$.
\end{proof}

Let $D^+\in \Psi_\maL^m (P | \maG^{(0)} ; \maE, \maE')^\maG$ be elliptic and denote by $S^*(P | \maG^{(0)})$ the vertical cosphere bundle of the fibration $\pi$. Let $F\in \Psi_\maL^0 (P | \maG^{(0)} ; \maE, \maE')^\maG$ be such that $$\sigma_0(F)|_{S^*(P | \maG^{(0)})}=\sigma_m(D^+)|_{S^*(P | \maG^{(0)})}\, .$$
 Choose a parametrix
$G$ for $F$.
 Let $U_{D^+} \in \Psi_\maL^{0} (P | \maG^{(0)} ; \maE\oplus \maE')^\maG$ be an invertible operator such that $U_{D^+}-\begin{bmatrix} 0 &-G\\
F &0 \end{bmatrix} \in  \Psi_\maL^{-\infty} (P | \maG^{(0)} ; \maE\oplus \maE')^\maG$. An explicit construction of an example of such
operator is as follows.  Let $S_0=1-GF$, $S_1=1-FG$. Then set $U_{D^+}=\begin{bmatrix} S_0 &-(1+S_0)G\\F &S_1 \end{bmatrix}$. With such a choice the inverse is given by an explicit formula  $U_{D^+}^{-1}=\begin{bmatrix} S_0 &(1+S_0)G\\-F &S_1 \end{bmatrix}$.

\begin{definition}\label{defindex}
The index of $D^+$ is the class in the $K$-theory of the algebra $\Psi_\maL^{-\infty}(P | \maG^{(0)}; \maE\oplus \maE')^\maG$ defined by
\begin{equation*}
\ind(D^+) = \left[P_{D^+} - Q \right] \in K_0(\Psi_\maL^{0} (P | \maG^{(0)}, \maE\oplus \maE')^\maG, \Psi_\maL^{-\infty} (P | \maG^{(0)}, \maE\oplus \maE')^\maG) \end{equation*}
\begin{equation*} \cong K_0(\Psi_\maL^{-\infty} (P | \maG^{(0)}, \maE\oplus \maE')^\maG),
\end{equation*}
where $P_{D^+}$ and $Q$  are the idempotents given by $P_{D^+}= U_{D^+}  \begin{bmatrix} 1 &0 \\0 &0\end{bmatrix} U_{D^+}^{-1}$ and $Q= \begin{bmatrix}0 &0\\0 &1 \end{bmatrix}$.
\end{definition}


\begin{remark}
              One can  use the regular representation to complete the, twisted by $L$, convolution algebra  of  $\maG$ to a $C^*$-algebra
        $C^*_r(\maG, L)$ and define the analytic index in its $K$-theory.
        \end{remark}

We leave to the reader the standard $K$-theoretic proof that the index is well defined and is stable under the
homotopies of $D^+$ in the class of elliptic operators in $\Psi_\maL^m (P | \maG^{(0)} ; \maE, \maE')^\maG$.
Assume that the transversally $\maL$-twisted bundles $\maE$ and $\maE'$ are hermitian (and in particular $\maL$ is unitary) and that the fibers of $\pi$ are equipped with $\cG$-equivariant smoothly varying volume forms. In this situation for a projective family $D^+$ we can again define a formally adjoint projective family $(D^+)^*$ by forming the formal adjoints for each family $D_{\alpha}$.

\begin{lemma}  Then, identifying $K_0(\Psi_\maL^{-\infty} (P | \maG^{(0)}, \maE\oplus \maE')^\maG)$ with $K_0(\Psi_\maL^{-\infty} (P | \maG^{(0)}, \maE'\oplus \maE)^\maG)$ we have
\begin{equation*}\ind (D^+)^*=-\ind D^+.\end{equation*}
\end{lemma}
\begin{proof}
We have
$$
\ind (D^+)^* = \left[U_{(D^+)^*}  \begin{bmatrix} 0 &0 \\0 &1\end{bmatrix} U_{(D^+)^*}^{-1} - \begin{bmatrix}1 &0\\0 &0 \end{bmatrix} \right] \in K_0(\Psi_\maL^{-\infty} (P | \maG^{(0)}, \maE\oplus \maE')^\maG).
$$
By deforming in the usual way $D^+$ we may assume that $\sigma_m(D^+)|_{S^*(P | \maG^{(0)})}$ is an isometry. In this case we may choose $U_{(D^+)^*}=U_{D^{+}}^{-1}$, and the statement follows.
\end{proof}


\section{The twisted Connes map}\label{morphismPhi}

\subsection{Review of cyclic homology}
Let us first review some basic material on the cyclic theory.
The general reference is the book \cite{Loday}, here we review the basic definitions to fix the notations.

Let $A$ be a complex unitary algebra. Set $C_k(A)= A\otimes (A/\C1)^{\otimes k}$ and $C^k (A)$ for the dual vector space of $C_k(A)$. Let $u$ be a formal variable of degree $-2$.
The space of negative cyclic chains of degree $ l \in \Z$ is defined by
\begin{equation*}
\CC^-_{l}(A) = \left(C_{\bullet}(A)[[u]]\right)_l=\prod \limits_{-2n+k=l,\ n\ge 0} u^nC_k(A).
\end{equation*}
The boundary is given by $b+uB$ where $b$ and $B$ are the Hochschild and Connes boundaries of the cyclic complex.
The homology of this complex is denoted $\HC^-_{\bullet}(A)$.
When the algebra $A$ is $\Z_2$ graded they incorporate the relevant signs.
If $A$ is not necessarily unital denote by $A^+$ its unitalisation and set $\CC^-_{l}(A)=\CC^-_{l}(A^+)/\CC^-_{l}(\C)$.
If $I$ is an ideal in a unital algebra $A$ the relative cyclic complex  is defined by
$\CC^-_{\bullet}(A, I) = \Ker \left( \CC^-_{\bullet}(A) \to \CC^-_{\bullet}(A/I)\right) $. One has a natural morphism of complexes $\iota \colon \CC^-_{\bullet}(I) \to \CC^-_{\bullet}(A, I)$ induced by the homomorphism $I^+ \to A$.

Dually, we define the cyclic cohomology complex by
\begin{equation*}
\CC^{l}(A) = \left(C^{\bullet}(A)[[\u]]\right)_l=\bigoplus \limits_{2n+k=l,\ n\ge 0} (\u)^nC^k(A).
\end{equation*}
with coboundary given by $b+\u B$, where $\u$ has degree $+2$.

Recall that for an algebra $A$ we have the Chern character in cyclic homology   $\ch \colon K_0(A)\to HC_{0}^{-}(A)$.
It is  defined by the following formula. Let $P, Q \in M_n(A^+)$ be two idempotents in $n\times n$ matrices of the
algebra $A^+$, representing a class $[P-Q] \in K_0(A)$. Then

\begin{equation}\label{cyclicchern}
{\Ch\left([P-Q]\right) = \tr(P-Q)+ \sum_{n=1}^{\infty} (-u)^n\frac{(2n)!}{n!} \tr\left( \left(P-\frac{1}{2}\right)\otimes P^{\otimes (2n)}-  \left(Q-\frac{1}{2}\right)\otimes Q^{\otimes (2n)} \right)}
\end{equation}
We will use abusively the notation $\Ch\left([ P-Q ]\right)$ for the cyclic cycle defined above and $\ch\left([P-Q]\right)$ for
its well defined class in cyclic homology $HC_0^-(A)$.

We will also need to use the entire cyclic complex. For our purpose the
algebraic version from  \cite{ConnesBook}, IV.7.$\alpha$ Remark 7 b.
will be sufficient.
First recall that one has the periodic cyclic complex $\left(\CC^{per}_{\bullet}(A), b+uB \right) $ where
$\CC^{per}_{\bullet}(A) = C_{\bullet}(A)[u^{-1}, u]]$.
Assume we are given a  periodic chain $\alpha = \sum_{k\ge 0} \alpha_k u^k \in \CC^{per}_{m}(A) $,  $\alpha_k \in C_{2k+m}(A)$. Then $\alpha$ is called entire if
there exist a finite dimensional subspace $V\subset A$, $1 \in V$ and $C>0$ (depending on $\alpha$) such that $\alpha_k \in V\otimes (V/\C1)^{\otimes k}$ and $\|\alpha_k\| \le  C^k k! $.
Here the norms on $V\otimes (V/\C1)^{\otimes k}$ are induced by an arbitrary
norm on $V$.
We denote the entire cyclic complex of $A$ by $\CC^{entire}_{\bullet}(A)$.

Note that the chain $\Ch\left([P-Q]\right)$ defined in \eqref{cyclicchern} is an element in $\CC^{entire}_{0}(A)$.

\subsection{The Connes map} $\;$

\medskip

Assume first that $M \to B$ is a submersion.
Let $\maL$ be a gerbe with connection on $B$ and $\maE$ a transversally $\maL$-twisted bundle on $M$ with connection. Let $\varOmega(\cL)$
be the associated curvature.

Following Mathai and Stevenson \cite{MathaiStevenson}, we construct a  morphism of complexes
\begin{equation*}\Phi_{\nabla^{\maH}} \colon  \maC^\bullet_{\varOmega (\cL)}(B) \longrightarrow \CC^{\bullet}(\Psi_\maL^{-\infty} (M|B; \maE)) \end{equation*}
 as follows.
Denote by $\Delta^l $ the standard $l$-simplex
\begin{equation*}
\Delta^l := \{(\t_0, \cdots , \t_{l})\in \R^{l+1}\ | \  0\leq \t_i,  \sum_{i=0}^l \t_i=1\}.
\end{equation*}

%
%
Set for any $\xi \in \maC_\bullet(B)[\u]$,
$$
\Phi_{\nabla^\maH}(\xi) = (-1)^{N(\xi)}\sum_k (\u)^k \phi_k(\xi)\, ,
$$
where $N(\xi)=\frac{\deg \xi(\deg \xi+1)}{2}$ and the cochain
 $\phi_k(\xi) \in C^{l}\left(\Psi_\maL^{-\infty} (M|B; \maE)\right)$ for
$l= \deg \xi -2k$,  is defined by:
\begin{multline}\label{MSequation}
\;\; \phi_k(\xi) (A_0, \cdots, A_l) := \\ \int_{\Delta^l}\left \langle (\u)^{-k}\xi,  \Tr\left( A_0 e^{-\t_0{u \theta^{\maH}}{}} \pa^\maH (A_1) e^{-\t_1{u \theta^{\maH}}{}} \cdots e^{-\t_{l-1}{u \theta^{\maH}}{}} \pa^\maH (A_l)e^{-\t_l{u \theta^{\maH}}{}}\right)\right\rangle d\t_1\ldots d\t_l
\end{multline}
where   $A_1, \ldots, A_l \in \Psi_\maL^{-\infty}$, $A_0 \in \Psi_\maL^{-\infty}+ \C \cdot 1$, with the convention that $\Tr(1)=0$. Recall the twisted de Rham complexes reviewed in Appendix \ref{TwistedCohom}.


\begin{theorem}[cf. \cite{MathaiStevenson}] \label{MS}The map
\begin{equation}
\Phi_{\nabla^{\maH}} \colon  \maC^\bullet_{\varOmega (\cL)}(B) \longrightarrow \CC^{\bullet}(\Psi_\maL^{-\infty} (M|B; \maE))
 \end{equation} defined above is a morphism of complexes.
\end{theorem}

\begin{proof} This is an immediate consequence of the result of \cite{MathaiStevenson}[page 323-]  and an observation that our map is the transposed
of the map they construct.
\end{proof}
This morphism depends on the choice of our horizontal distribution $\maH$. However the results of \cite{MathaiStevenson} show that a different choice of $\maH$ leads to a chain homotopic morphism.
\begin{remark}
Note that the similar morphism in the non-twisted case was constructed
in the references \cite{BenameurHeitsch1, Gorokhovsky1, NestTsygan1, NestTsygan2}.
\end{remark}

We now  proceed to explain the construction of the $\Phi$-map in our general case. Our constructions are based on the previous results obtained  for discrete translation groupoids  in \cite{Angel}.

Assume again that we have an \'etale groupoid $\maG$, and a gerbe $L$ with connection on $\maG$. Assume that $\pi: P \to \maG^{(0)}$ is a $\maG$-equivariant
submersion and that $\maE$ is a  transversally $L$-twisted vector bundle with connection  on $P$.  Assume as before that the action of $\maG$ on $P$ is free, proper and for simplicity cocompact.
Recall that we also have for every $k \in \N$, a submersion $\cP^{(k)} \to \cG^{(k)}$, a gerbe $\maL^{(k)}$ on $\cG^{(k)}$ and a transversally  $\maL^{(k)}$-twisted vector
bundle $\cE^{(k)}$ on $\cP^{(k)}$ (recall that $\cP^{(k)}$ is defined in Section \ref{AppendixHomology}).
We thus obtain the following equality of sheaves
 $$
 \Psi_{\maL^{(k)}}(\pi^{(k)}; \maE^{(k)}) =(\pr_0^k)^* \Psi_{\maL}(\pi; \maE)\,
 $$
on every $\cG^{(k)}$.
See Subsection \ref{Projective families} for the definitions.
Using the constructions of Section \ref{projective} we obtain for every $k$ a submersion  $\cP^{(k)}\times \Delta^k \to \cG^{(k)}\times \Delta^k$, a gerbe $\cL^{(k)}_\Delta$
on  $\cG^{(k)}\times \Delta^k$ and a transversally $\cL^{(k)}_\Delta$-twisted  bundle $\cE^{(k)}_\Delta$ on $\cP^{(k)}\times \Delta^k$.

We can view  the algebra of operators $\Psi_{\maL^{(k)}}(\cP^{(k)}|\cG^{(k)}; \maE^{(k)})$   as  a subalgebra in
$\Psi_{\maL^{(k)}_\Delta}(\cP^{(k)}\times \Delta^k|\cG^{(k)}\times \Delta^k; \maE^{(k)}_\Delta)$ constant in the $\Delta^k$
direction, i.e. pulled back to $ \cG^{(k)}\times \Delta^k$ via the projection on the first factor.

Recall that the gerbe  $\cL^{(k)}_\Delta$ is equipped with connection given by \eqref{gconnective_structure}, \eqref{gcurving} and $\cE^{(k)}_\Delta$ is equipped with connection given in \eqref{simplicial connection}.
Since the action of $\maG$ on $P$ is proper, there exists a horizontal distribution $\maH$ on $P \to \maG^{(0)}$ which is $\maG$-invariant.
We therefore obtain a derivation (cf. \eqref{defpartial})
\begin{equation*}
\partial^{(k)} \colon \Psi_{\maL^{(k)}_\Delta}(\cP^{(k)}\times \Delta^k|\cG^{(k)}\times \Delta^k; \maE^{(k)}_\Delta)  \longrightarrow   \Omega^1(\cG^{(k)}, \Psi_{\maL^{(k)}_\Delta}(\cP^{(k)}\times \Delta^k|\cG^{(k)}\times \Delta^k; \maE^{(k)}_\Delta)),
\end{equation*}
and we get by \eqref{deftheta}, $\theta^{(k)} \in \Omega^2(\cG^{(k)}\times \Delta^k, \Psi_{\maL^{(k)}_\Delta}(\cP^{(k)}\times \Delta^k|\cG^{(k)}\times \Delta^k; \maE^{(k)}_\Delta))$,
such that $(\partial^{(k)})^2 = [\theta^{(k)}, \cdot]$.
We then have by straightforward verification the following identities
\begin{equation*}
(\del_i \times \id)^* \partial^{(k-1)} =
(\id \times \delta_i)^*  \partial^{(k)} \text{ for } i \ge 0.
\end{equation*}
\begin{equation*}
(\del_i \times \id)^* \theta^{(k-1)} =
(\id \times \delta_i)^*  \theta^{(k)} \text{ for } i \ge 0.
\end{equation*}
Note also that $\theta^{(k)}$ has components of $(2, 0)$ and $(1, 1)$ types only.

For a cocycle $c= \{c_k\}\; \in \maC_\Delta(\maG)^\bullet $,  we set
\begin{equation*}
\Phi^\cE(c) = (-1)^{N(c)} \sum  (\u)^j\phi_j(c)
\end{equation*}
where $\phi_j(c) = \sum_{k+l = \deg c -2j} \phi_j^{kl}(c_k)$ and  $\phi_j^{kl}(c_k) \in \CC^\bullet (\Psi^{-\infty}_\maL(P|\maG^{(0)}, \maE))$
is given by
\begin{multline}\label{philk}
\;\; \phi_j^{kl}(c_k)(A_0, \cdots, A_l) := \\     \int \limits_{\Delta^l}\left \langle (\u)^{-j} c_k,  \Tr\left( A_0 e^{-\t_0 u\theta^{(k)}} \pa^{(k)} (A_1) e^{-\t_1 u\theta^{(k)}} \cdots e^{-\t_{l-1} u\theta^{(k)}} \pa^{(k)} (A_l)e^{-\t_l u\theta^{(k)}}\right) \right\rangle d\t_1\ldots d\t_l
\end{multline}
Note that this expression is $0$ if $l > \dim \maG^{(0)}$ or $k > \dim  \,\maG^{(0)}+\deg \omega $, and therefore the sum in the formula for $\Phi^\cE$ is finite.

Recall that  the complex $\maC_\Delta( \maG)^\bullet$ is endowed with the differential $d^t_{\varOmega(L)}$ defined by \eqref{defdtsimp}
We then have the following simplicial analogue of  Theorem \ref{MS}
\begin{theorem}\label{PhiE}
The map $\Phi^\cE$  is
a morphism of complexes
\begin{equation*}
\Phi^\cE \colon (\maC_\Delta(\maG)^\bullet, d^t_{\varOmega(L)}) \longrightarrow C^\bullet(\maG, \CC^\bullet (\Psi^{-\infty}_\maL(P|\maG^{(0)}, \maE))).
\end{equation*}
\end{theorem}
Recall that $\CC^\bullet (\Psi^{-\infty}_\maL(P|\maG^{(0)}, \maE))$ defines a $\maG$-sheaf of algebras and the complex on the right hand side of the theorem is thus well defined, see also Appendix \ref{cross}. The proof is an immediate rephrasing of the classical one of Theorem \ref{MS} and is omitted. The above theorem is also an extension of the similar statement for discrete translation groupoids proved in \cite{Angel} and which corresponds to suspended foliations.

Using again the notations from Appendix \ref{cross} and composing the induced morphism on cohomology with the  map
\begin{equation*}
H^\bullet(C^\bullet(\maG, \CC^\bullet (\Psi_\maL^{-\infty}(P|\maG^{(0)}, \maE)))) \longrightarrow \HC^\bullet(\Psi^{-\infty}_\maL(P|\maG^{(0)}, \maE) \rtimes \maG )
 \,,
\end{equation*}
 which comes from \cite{bryni, crainic}, and then with the morphism
\[
I^* \colon \HC^\bullet(\Psi^{-\infty}_\maL(P|\maG^{(0)}, \maE) \rtimes \maG )
 \longrightarrow  \HC^\bullet (\Psi^{-\infty}_\maL(P|\maG^{(0)}, \maE)^\maG)\,
\]
induced by the map $I$ constructed in  Lemma \ref{mapI}, we eventually end up with the  map
\begin{equation*}
\Lambda^\cE \colon H^\bullet (\maC_\Delta(\maG)^\bullet, d^t_{\varOmega(L)}) \longrightarrow \HC^\bullet (\Psi^{-\infty}_\maL(P|\maG^{(0)}, \maE)^\maG).
\end{equation*}
Said differently, we end up with the morphism
\begin{equation}\label{eq:Lamb}
\Lambda^{\cE} \colon H_{\varOmega(L)}^\bullet (B\cG) \longrightarrow \HC^\bullet (\Psi^{-\infty}_\maL(P|\maG^{(0)}, \maE)^\maG).
\end{equation}


 Assume now that  the bundle $\cE$ is $\Z_2$ graded, $\cE= \cE^+\oplus \cE^-$ with the grading operator $\gamma$, $\gamma|_{\cE^\pm} = \pm1$.
In this case $\Psi^{-\infty}_\maL(\pi, \maE)$ is a sheaf of $\Z_2$ graded algebras  with the even part $ (\Psi^{-\infty}_\maL(\pi, \maE)) ^+ = \{A | A\gamma= \gamma A \}$.

We can use  the supertrace $\STr = \Tr \circ \gamma$ instead of  the trace in the formula \eqref{philk} to obtain  cyclic cochains on the algebras
$\Psi^{-\infty}_{\maL^{(k)}}(\cP^{(k)} |\cG^{(k)} ; \maE^{(k)})^+$ and the following generalization of Theorem \ref{PhiE}.
\begin{theorem}
We have a   morphism of complexes
\begin{equation*}
\Phi^{\cE, \gamma} \colon (\maC_\Delta(\maG)^\bullet, d^t_{\varOmega(L)}) \longrightarrow C^\bullet(\maG, \CC^\bullet (\Psi^{-\infty}_\maL((P|\maG^{(0)}, \maE))^+)),
\end{equation*}
which allows to  define  our  twisted  map on cohomology:
\begin{equation*}
\Lambda^{\cE, \gamma} \colon H_{\varOmega(L)}^\bullet (B\cG) \longrightarrow \HC^\bullet ((\Psi^{-\infty}_\maL(P|\maG^{(0)}, \maE))^+){{^\maG}}).
\end{equation*}
\end{theorem}



We are now in position to define the higher $c$-index associated with any element $[c]\in H^\bullet (\maC_\Delta(\maG)^\bullet, d^t_{\varOmega(L)})$. More precisely,

\begin{definition}
We  define the index pairing for a $\maG$-invariant $L$-projective elliptic  family $D$, with the groupoid cocycle $c$ as
\begin{equation*}
\ind_c (D) := \langle (u^*)^{-\frac{\deg c}{2}}\Lambda^{\cE} (c), \ch (\ind D) \rangle.
\end{equation*}
\end{definition}
This higher twisted $c$-index  $\ind_c( D)$ is a complex number which is  computed  by a local formula {\em \`a la Atiyah-Singer} in the next section.\\

\section{The higher twisted index formula}\label{Dirac}

The goal of this section is to give a superconnection construction of Connes' $\Phi$ map which will be used to prove  our higher  index theorem for  $L$-projective leafwise Dirac operators. Recall that we assume that $\cG$ acts freely, properly and co-compactly on $P$.

\subsection{Projective superconnections on fibrations}\label{Dirac operators}
We first review some material from \cite{BenaGoro} and consider the smooth  submersion $\pi: M\to B$. Assume that the fibers of $\pi$  are even dimensional complete Riemannian manifolds.
 Denote by $C(M|B)$ the Clifford algebra of the fiberwise cotangent bundle $T^*(M|B)=T^*M/ (\ker \pi_*)^{\perp}$.
Let $\maL$ be a unitary gerbe on $B$.
\begin{definition}\
\begin{itemize}
 \item  A twisted Clifford module is a transversally  $\maL$-twisted Hermitian $\Z_2$-graded vector bundle $\maE=\maE^+ \oplus \maE^-$ on $M$ together with the homomorphism $c \colon  C(M|B) \rightarrow \End (\maE)$ of bundles of unital $\Z_2$-graded $*$-algebras.
\item A Clifford  connection $\nabla^{\maE}$ on $\maE$ is an Hermitian  connection such that for any section $s$ of $\maE$ one has:
$$\nabla^{\maE} (c(a) s )=c(\nabla^{M|B} (a)) s + c(a) \nabla^{\maE} ( s ) \,.$$
\end{itemize}

\end{definition}
Clifford connections on Clifford modules always exist; the proof  is analogous to the
proof of existence of twisted connections in Lemma \ref{connection existence}.

Choose a horizontal distribution $\maH \subset TM$ so that
$
TM = \maH \oplus T (M|B).
$
This choice together with the Riemannian metric on the fibers of $\pi$ allows one to define a connection $\nabla^{M|B}$ on the fiberwise tangent bundle $T(M|B)$, see \cite{bgv} Section 10.1.
We denote by $R^{M|B}$ the curvature of this connection.
 Set
 $$
 \End_{C(M|B)}(\maE)= \{A \in \End(\maE)\ | \ [A, c(a)]=0 \text{ for every } a \in C(M|B) \}.
 $$
Let $\Gamma \in C(M|B)$ be the chirality operator defined locally by $\Gamma = i^{k / 2} e^1\ldots e^k$ where $k=\dim M-\dim B$ and $e^1, \ldots, e^k$ is the local orthonormal basis of $T^*(M|B)$. Define then the relative supertrace
$$
\str_{\maE/\maS} \colon  \End_{C(M|B)}(\maE) \longrightarrow C^{\infty} (M)\text{ by }\str_{\maE/\maS}(A)= 2^{-k/2} \str c(\Gamma) A.
$$
We fix from now on a Clifford connection $\nabla^{\maE}$ on $\maE$ and a descent datum $(U_{\alpha}, \maL_{\alpha \beta},
\mu_{\alpha \beta \gamma}, \maE_{\alpha}, \varphi_{\alpha \beta})$ for the transversally $\maL$-twisted Clifford
module $\maE$. The connection $\nabla^{\maE}$ defines a connection $(\nabla^{\maE}_{\alpha}, \nabla_{\alpha \beta}, \omega_{\alpha})$ (defined up to equivalence) on this descent datum. Each $\maE_{\alpha}$ is then a Clifford module on the submersion $\pi^{-1}U_{\alpha} \to U_{\alpha}$, and
each connection $\nabla^{\maE}_{\alpha}$ is a Clifford connection.

Recall  (see Proposition \ref{global curvature}) that one defines $\theta^{\maE} \in \Omega^2(M, \End(\maE))$
by setting $\theta^{\maE}|_{\pi^{-1}U_{\alpha}}= \theta^{\maE}_{\alpha} +\pi^* \omega_{\alpha}$.
Denote by  $c(R^{M|B})$  the action of the $2$-form with values in the
Clifford algebra obtained from $R^{M|B}$ via the Lie algebra isomorphism $\mathfrak{so}(T(M|B)) \to C^2(M|B)$.
Here $C^2(M|B) \subset C(M|B)$ is the subspace consisting of elements $\sum u_i v_i$, $u_i, v_i \in T^*(M|B)$ with $\sum \langle u_i, v_i \rangle=0$.
Define $\theta^{\maE/\maS}=\theta^{\maE}-c(R^{M|B})$.

The argument  in \cite{bgv}, Proposition 3.43,
shows that $\theta^{\maE/\maS} \in \Omega^2(M, \End_{C(M|B)}(\maE)) $. Over an open set $\pi^{-1}U_{\alpha}$ we have
$\theta^{\maE/\maS}|_{U_{\alpha}}=\theta^{\maE/\maS}_{\alpha} +\pi^*\omega_{\alpha}$ where $\theta^{\maE/\maS}_{\alpha} \in \End_{C(M|B)}(\maE_{\alpha})$ is defined via the equality $\theta^{\maE}_{\alpha}=\theta^{\maE/\maS}_{\alpha}+ c(R^{M|B})$. We can then define a differential form $\Ch_{\maL}(\maE/\maS)$ by
\begin{equation}\label{deftwistedchern}
\Ch_{\maL}(\maE/\maS) := \str_{\maE/\maS} e^{-u \theta^{\maE/\maS}}  \in \Omega^*(M)[u].
\end{equation}
The  proof of the following result  is standard and analogous to the proof of  Proposition \ref{twchern}.
\begin{lemma}
We have $ d_{\pi^*\varOmega}\, \Ch_{\maL}(\maE/\maS) = 0$ and the corresponding class is $H^*_{\pi^*\varOmega}(M)$ is independent of the
choice of Clifford connection $\nabla^{\maE}$.
\end{lemma}
We also introduce the fiberwise $\widehat{A}$-genus by $\widehat{A}(TM|B)= \widehat{A}\left(uR^{M|B}\right)\in \Omega^*(M)[u]$, where $ \widehat{A}(x)$ is the power series defined by
\begin{equation*}
\widehat{A}(x)=\left.\det\right.^{1/2}\left(\frac{x/2}{\sinh x/2}\right).
\end{equation*}

Using the above data,  we can define on each  submersion $\pi^{-1}U_{\alpha} \to U_{\alpha}$ a family of Dirac operators $D_{\alpha}$ acting on the sections of the bundle $\maE_{\alpha}$. Locally $D_{\alpha} = \sum_i c(e^i) \left(\nabla^{\maE}_{\alpha}\right)_{e_i}$ where $\{e_i\}$, $\{e^i\}$ are dual bases of $T(M|B)$ and $T^*(M|B)$ respectively. {The following is obvious from the definitions.

\begin{lemma} The collection $D=\{ D_{\alpha}\}$ defines an element in $\Psi_\maL^1(M|B; \maE)$.
\end{lemma}
With respect to the decomposition $\maE= \maE^+\oplus\maE^-$, the Dirac operator $D$ decomposes as $D= D^+\oplus D^-$ where
$$D^+ \in \Psi_\maL^1(M|B; \maE^+, \maE^-),\;D^- \in \Psi_\maL^1(M|B; \maE^-, \maE^+)\,.$$

Classical arguments show that $D^+$ is fiberwise elliptic and hence the analytical index $\ind (D^+)$ of $D^+$ is well defined in $K(\Psi_\maL^{-\infty}(M|B; \maE))$.


Recall that a twisted superconnection $\A$ on the descent datum $(U_{\alpha}, \maL_{\alpha \beta}, \mu_{\alpha \beta \gamma}, \maE_{\alpha}, \varphi_{\alpha \beta})$  is a collection
$(\A_\alpha)_{\alpha\in \Lambda}$ of superconnections on the vector bundles $\pi_*\maE_\alpha$ over the open sets $U_\alpha$  such that when $U_{\alpha\beta} \not = \emptyset$,
\begin{equation}\label{suco}
(\pi_*\varphi_{\alpha\beta})^*   \A_\beta   = \A_\alpha \otimes \id + \id\otimes \nabla_{\alpha\beta}.
\end{equation}

We say that $\A$ is a Bismut superconnection if each $\A_{\alpha}$ is. Specifically we have
\begin{equation*}
\A_{\alpha}= u^{-1/2} D_{\alpha}+\nabla_{\alpha}^{\maH} -\frac{1}{4}u^{1/2}c(T^{\maH})
\end{equation*}
where $T^{\maH}$ is defined in \eqref{th} and we refer to \cite{bgv} for the precise definition of the Bismut superconnection $\A_{\alpha}$.

  The facts that $\nabla^{\maH}$ is a connection
on the descent datum, that $D\in \Psi^1_{{\maL}} (M|B, \maE)$ and that $c(T^{\maH}) \in \Omega^2(B, \Psi^0(M|B, \maE))$,
imply that the conditions \eqref{suco} are satisfied.
From now on $\A$ will denote the  Bismut superconnection.

We will also consider as before the  rescaled Bismut superconnection  $\A_s$
defined by
\begin{equation*}
(\A_s)_{\alpha}= \A_{\alpha, s}: = s^{1/2}u^{-1/2}D_{\alpha}+\nabla_{\alpha}^{\maH} -\frac{1}{4}u^{1/2}s^{-1/2}c(T^{\maH})
\end{equation*}
where $s$ is either  a positive number or a multiple of the formal variable $u$. Denote by  $\theta^{\A_s}_{\alpha}$ the curvature of the rescaled superconnection
 $\A_{\alpha, s}$.
In particular we have forms $u\theta^{\A_{\alpha, u^{-1}}} \in \Omega^*(U_{\alpha}, \Psi(\pi^{-1}U_{\alpha}|U_{\alpha}, \maE_{\alpha}))[u^{1/2}]$.

\begin{proposition} \label{prop:theta}There exists a form $u\theta^{\A} \in \Omega^*(B, \Psi_{{\maL}}(M|B; \maE))[u^{1/2}]$ such that
\begin{equation*}\left.u\theta^{\A}\right|_{U_{\alpha}}= u(\theta^{\A_{\alpha, u^{-1}}} + \pi^* \omega_{\alpha}).
\end{equation*}
\end{proposition}

\begin{proof}
Recall that the curvature of $\nabla_{\alpha\beta}$ is equal to
 $\omega_\alpha- \omega_\beta$. Therefore from the equation \eqref{suco} we obtain
$$
\phi_{\alpha\beta} \left(u\theta^{\A_{\beta, u^{-1}}}\right)= u(\theta^{\A_{\alpha, u^{-1}}} + \pi^*\omega_\alpha - \pi^*\omega_\beta).
$$
The statement of the Proposition follows.
\end{proof}

Notice that $u\theta^{\A} = D^2 +$ forms of degree $>0$. Therefore we can
define
$$e^{u \theta^{\A}} \in \Omega^*(B, \Psi^{-\infty}_{{\maL}}(M|B; \maE))[u^{1/2}]$$
 by the usual Duhamel formula.

 Note that again the parity considerations as in the finite dimensional case show that the coefficients for the nonintegral
powers of $u$ are odd with respect to the grading and hence have a vanishing supertrace.

\subsection{Projective superconnections on foliations} $\;$

\medskip

Assume now that we are in the situation of Section \ref{morphismPhi} with the $\cG$-submersion $P\to \cG^{(0)}$. Assume that the twisted Clifford module  $\maE$ is $\Z_2$ graded and
$$D =\begin{bmatrix} 0 & D^- \\
D^+ & 0 \end{bmatrix}$$ is the Dirac type family
of operators.
We construct for every $k\in \N^*$ the Bismut superconnection $\A^{(k)}$ on the fibration $\cP^{(k)}\times \Delta^k \to \cG^{(k)}\times \Delta^k$.
Proposition \ref{prop:theta}  applied with $M=\cP^{(k)}\times \Delta^k$ and $B=\cG^{(k)}\times \Delta^k  $ gives the form $u \theta^{\A^{(k)}}$.

We
 similarly produce a collection of morphisms in the following way.
For a cochain $c \in (\maC_\Delta(\maG)^\bullet[u], d^t_{\varOmega(L)})$, $c= \{c_k\}$ we have
\begin{equation*}
\Phi^\A(c) =  (-1)^{N(c)} \sum_{j= -\infty}^m  (\u)^j\phi^\A_j(c)
\end{equation*}
where $m =\lfloor\frac{1}{2}\deg c\rfloor$, $N(c)=\frac{\deg c(\deg c+1)}{2}$, $\phi^\A_j(c) = \sum (\phi^\A_j)^{kl}(c_k)$ and $(\phi^\A_j)^{kl}$ is the element of  $C^k( \cG, C^l_{entire} (\Psi^{-\infty}_\maL(P|\maG^{(0)}))$ defined for $k+l = \deg c -2j$ by

\begin{multline}
(\phi^\A_j)^{kl}(c_k)(A_0, \cdots, A_l) := \\     \int \limits_{\Delta^l}\left \langle (\u)^{-j} c_k,  \STr \left( A_0 e^{-\t_0 u \theta^{\A^{(k)}}} [\A^{(k)}, A_1] e^{-\t_1 u \theta^{\A^{(k)}}} \cdots e^{-\t_{l-1} u \theta^{\A^{(k)}}} [\A^{(k)}, A_l] e^{ -\t_l u \theta^{\A^{(k)}}}\right)  \right\rangle d\t_1\ldots d\t_l
\end{multline}

Note that only integral powers of $u$ appear in this formula since the coefficients of half-integral powers of $u$ involve
the supertrace of odd operators. Notice that this formula extends  the usual untwisted one and yields an entire cyclic cocycle,  see for instance \cite{BenameurCarey}.

\medskip

Consider now the  sheaf of algebras $\cA \subset  \Psi^0_\maL(\pi, \maE)$ on $\cG^{(0)}$ given by:
\begin{equation*}
U \text { open subset of }   \cG^{(0)} \longrightarrow \cA(U)= \{ A \in \Psi^0_\maL(\pi, \maE)(U) \ | \ [D, A] \in \Psi^0_\maL(\pi, \maE) \}.
\end{equation*}
When no confusion arises we will also denote by $\cA$ the algebra of global sections, which is a subalgebra of
$\Psi^0_\maL(P|\maG^{(0)}, \maE)$. This is a $\cG$-subsheaf of $\Psi^0_\maL(\pi, \maE)$.
Consider now the complex $C^\bullet(\maG, \CC^\bullet_{entire} (\cA))$. Recall
that this is a total complex of the bicomplex given in the bidegree $l, k$ by
$\CC^l( \cA^{(k)})$, where $\cA^{(k)} = (\pr^k_0)^* \cA (\maG^{(k)}) \subset \Psi^0_\maL(\cP^{(k)}|\maG^{(k)}, \maE^{(k)})$. The cochains $(\Phi_\A^{(k)})^l(c)$, defined by the same formulas,
extend to define elements in $C^\bullet(\maG, \CC^\bullet_{entire} (\cA))$.

We therefore obtain a morphism:
\begin{equation*}
\Phi_\A \colon (\maC_\Delta(\maG)^\bullet, d^t_{\varOmega(L)}) \longrightarrow C^\bullet(\maG, \CC^\bullet_{entire} (\cA)).
\end{equation*}

For every $ 0 <s \le 1$ we also have the rescaled superconnection $\A_s$ and the corresponding
morphism $\Phi_{\A_s}$.
\medskip

We will now apply this construction  of the morphism $\Phi_\A$ to the groupoid $\cG \times \R^{+}$.  Let $s$ be the coordinate in $\R^{+}$ and $d_s$ the de Rham differential on $\R^{+}$.
Equip every submersion $\cP^{(k)}\times \Delta^k \to \cG^{(k)}\times \Delta^k$
  with the superconnection $d_s+\mathbb{A}_s$.
We obtain the morphism
\begin{equation*}
\Phi_{d_s+\A_s} \colon (\maC_\Delta^{vc}(\maG\times \R^{+})^\bullet, \u d_s^t+d^t_{\varOmega(L)}) \to C^\bullet(\maG, \CC^\bullet_{entire} (\cA)).
\end{equation*}
Here $\maC_\Delta^{vc}(\maG\times \R^{+})^\bullet$ denotes the currents with fiberwise compact support on (trivial) fibrations
$\maG^{(k)}\times \Delta^k \times \R^{+} \to \maG^{(k)}\times \Delta^k $, and $\cA$ is viewed as sections on
 $\cG \times \R^{+}$   constant in $\R^+$-direction.

For every $\tau \in \R^+$ let $i_\tau \colon \Omega^\bullet( \maG\times \R^{+}) \to \Omega^{\bullet-1}( \maG)$ be the {{map}} defined by
\[
i_\tau (\alpha + \beta ds) =\beta|_{s=\tau}.
\]
By transposition it induces a {{map}}
\[
i_\tau^t \colon \maC_\Delta(\maG)^\bullet \longrightarrow  \maC_\Delta^{vc}(\maG\times \R^{+})^{\bullet+1}
\]
Composing it with $\Phi_{d_s+\A_s}$ we obtain the {{transgression map}}

\begin{equation*}
\Psi_{\A_s}:= \Phi_{d_s+\A_s} \circ i_\tau^t \colon \maC_\Delta(\maG)^\bullet \longrightarrow C^{\bullet+1}(\maG, \CC^\bullet_{entire} (\cA)).
\end{equation*}
The following lemma is the exact twisted version of a classical relation, it can be verified by a direct computation:
\begin{lemma} \label{lem:trans}

\begin{equation*}
\u \frac{d}{ds} \Phi_{\A_s} = (b\pm \del+\u B) \circ \Psi_{\A_s} + \Psi_{\A_s} \circ d^t_{\varOmega(L)}
\end{equation*}
\end{lemma}

We now fix a cocycle $c \in (\maC_\Delta(\maG)^\bullet[u], d^t_{\varOmega(L)})$, and write
  $\phi_j(s)$ for $\phi^{\A_s}_j(c)$  and $\phi_j^{kl}(s)$ for $(\phi^{\A_s}_j)^{kl}(c)$. So $\phi_j^{pq}(s) \in C^p( \cG, C^q (\cA))$, $p+q = \deg c -2j$.
Write similarly
\begin{equation*}
\Psi_{\A_s}(c) = (-1)^{N(c)}\sum_{j= -\infty}^{m+1} (\u)^j\psi_j(s)
\end{equation*}
Then
\begin{equation*}
\frac{d}{ds}\phi_j(s) = (b\pm \del) \psi_{j+1}(s) + B\psi_{j}(s)
\end{equation*}
We also use the notations of Appendix \ref{JLO}, and we can state the following consequence of  Lemma \ref{lem:sub} and Lemma \ref{trace estimate}:

\begin{proposition}  \label{prop:estim} For a current $c$  there exist constants  $M$ and $C=C(j)$    such that
\begin{equation*}
|\phi_j^{kl}(s)(A_0, A_1, \ldots A_l)| \le C \mu(A_0) \nu_n(A_1) \ldots \nu_n(A_l) s^{-M -j/2}
\end{equation*}
and
\begin{equation*}
|\psi_j^{kl}(s)(A_0, A_1, \ldots A_l)| \le C \mu(A_0) \nu_n(A_1) \ldots \nu_n(A_l) s^{-M -(j+1)/2}
\end{equation*}
where $\nu_n(A) =\max \{\|(1+D^2)^i A (1+D^2)^{-i}\|\ |\ |i|\le n\}+ \max \{\|(1+D^2)^i [D, A] (1+D^2)^{-i}\|\ |\ |i|\le n\}$.
\end{proposition}
\begin{proof}   First consider the expression  $e^{-\t u \theta^{\A_s^{(k)}}}$.
By the Duhamel expansion formula, this  is a finite sum (polynomial in $u$) of the form $\sum u^i P_i$.
Each term of this expansion can be bounded by using Lemma \ref{lem:sub}. We obtain that for sufficiently large $r$  and for each $i$ the (fiberwise) norm of  $(1+D^2)^{-r} P_i$ is bounded uniformly in $s$, $\t$. If we assume also that $\t \ge \frac{1}{l}$ then we have:
\begin{equation} \label{termest}
\|(1+D^2)^{r} P_i\| \le C s^{-r-M},
\end{equation}
for some universal constant $M$. Consider now $\phi_j^{kl}(s)(A_0, A_1, \ldots , A_l)$. As in Lemma \ref{lem:sub} it is sufficient to bound the integrand assuming $\t_i \ge 1/l$ and then add the estimates. The expansion of the supertrace in powers of $u$ contains $u^j$ multiplied by
$s$ to the power $-j/2$ or higher and by the integral which can be bounded as in Lemmas \ref{lem:sub}, \ref{trace estimate} using the estimate \eqref{termest}.  The result follows immediately.

\end{proof}
The entire cocycle $\Phi_{\A_s}(c)$ can indeed be embodied by  a periodic one. Set  for $p>>1$  sufficiently large
\begin{equation}\label{deftpf}
T_p\Phi_{\A_s}(c):=  (-1)^{N(c)}\left( \sum_{j= -p}^m (\u)^j\phi_j(s) - (\u)^{-p}B \int_0^s \psi_{-p-1}(t) dt \right)\ .
\end{equation}
We use the notation $T_p\Phi_\A$ for  $T_p\Phi_{\A_1}$.

\begin{proposition}For every $0<s \le 1$ and sufficiently large $p$,  $T_p\Phi_{\A_s}(c)$ is a cocycle in the complex $C^{\bullet}(\maG, \CC^\bullet_{per} (\cA))$.
Its cohomology class does not depend on choices of $s$ and $p$.
\end{proposition}
\begin{proof}
We note that Proposition \ref{prop:estim} implies that for sufficiently large $p$ ($p > 2M$ in the notations of the proposition)
$\phi_{-p}(s) \to 0$ as $s \to 0$, and the integral $\int_0^s \psi_{-p-1}(t) dt$ converges, when evaluated on collections of operators $A_i$
with uniformly bounded norms $\mu$, $\nu_n$. It follows that
\begin{equation*}
 \phi_{-p}(s) = (b\pm \del) \int_0^s \psi_{-p}(t)dt + B \int_0^s\psi_{-p-1}(t)dt
\end{equation*}
We also have identities
\[(b\pm\del) \phi_j(s) + B\phi_{j-1}(s)=0.
\]

Together these identities easily imply that $T_p\Phi_{\A_s}(c)$ is a cocycle in the complex $C^{\bullet}(\maG, \CC^\bullet_{per} (\cA))$ for any $s$ and sufficiently large $p$.
We also have
\[
\u \frac{d}{ds}T_p\Phi_{\A_s}(c)= (-1)^{N(c)}(b\pm \del +\u B)\sum_{j= -p}^{m+1} (\u)^j\psi_j(s)
\]
Finally  for $q>p$, we have
\[\u \left( T_p\Phi_{\A_s}(c) - T_q\Phi_{\A_s}(c) \right)= (-1)^{N(c)}(b\pm \del +\u B)\sum_{j= -q}^{-p-1} (\u)^j\psi_j(s)
\]
\end{proof}
\medskip
Now, applying the map
\begin{equation*}
H^\bullet(C^\bullet(\maG, \CC^\bullet ({{\maA}})) \longrightarrow \HC^\bullet({{\maA}} \rtimes \maG )\,,
\end{equation*}
constructed in \cite{crainic}, see also  \cite{bryni}, together with the map
\begin{equation*}
\HC^\bullet({{\maA}} \rtimes \maG ) \longrightarrow \HC^\bullet ({{\maA}}^\maG)\,,
\end{equation*}
to the class of the above cocycle $[T_p\Phi_\A(c)]$,
we eventually end up with the  class $\Lambda_\A (c)$ {{which thus lives in the periodic cyclic cohomology}} $\HC^\bullet_{per} (\cA^\cG)$.
\begin{proposition}\label{transgress}
Let $\iota \colon \Psi^{-\infty}_\maL(P|\maG^{(0)}, \maE) \to \cA$ be the  inclusion. Then
\begin{equation*}
[\iota^*(T_p\Phi_{\A}(c))] =[\Phi^{\cE, \gamma} (c)] \text{ and }  \iota^* \Lambda_{\A} (c) =  \Lambda^{\cE, \gamma}(c)\,.
\end{equation*}
See Remark 6.3 for the definition of $\Lambda^{\cE, \gamma}(c)$.
\end{proposition}

\begin{proof} First notice that the second equality is an immediate consequence of the first. To prove the first notice that applying the construction of the morphism
$T_p\Phi_{\A}$ to the groupoid $\cG \times \R{{^+}}$, as in discussion preceding Lemma \ref{lem:trans}, one shows that the class of $T_p\Phi_{\A}(c)$ is independent of superconnection $\A$. Moreover, when restricted to the algebra of smoothing operators one can homotop $\A$ to the connection $\nabla^\cH$ and obtain that the cocycles
$\iota^*T_p\Phi_{\A}(c)$ and  $\Phi^{\cE, \gamma} (c)$ are cohomologous.
\end{proof}

We are now in position to compute our higher index in terms of the superconnection $\A$.

\begin{theorem} \label{superconnection computes index} For any cocycle $c \in (\maC_\Delta(\maG)^\bullet[\u], d^t_{\varOmega(L)})$, the higher twisted index is given by:
\begin{equation*}
\ind_c (D^+) \; = \;  \left \langle (u^*)^{-\frac{ {\rm deg}\, c }{2}}\Lambda_{\A} (c),  [1]  \right \rangle.
\end{equation*}

\end{theorem}
\begin{proof}
 Consider the twisted bundle $\widetilde{\maE}=\maE\oplus \maE$ with the grading given by $\Gamma=\begin{bmatrix} \gamma &0\\ 0 &-\gamma \end{bmatrix}$. The algebra  $\Psi_\maL (P|\cG^{(0)}; \widetilde \maE)$ of operators on $\widetilde{\maE}$
is naturally $\Z_2$ graded.
For an operator or (super) connection $K$ on $\maE$ set $\widetilde{K} = K \oplus K$; so for example $\widetilde{\A} =\A\oplus\A$, $\widetilde{\nabla}^{\maH}=\nabla^{\maH}\oplus \nabla^{\maH}$, etc. In particular $\widetilde{\cA}$ is a sheaf defined by
\begin{equation*}
\widetilde{\cA}(U)= \{ A \in \Psi^0_\maL(\pi, \widetilde {\maE})(U) \ | \ [\widetilde{D}, A] \in \Psi^0_\maL(\pi, \widetilde{\maE}) \}.
\end{equation*}

Let $F=D(1+D^2)^{-1/2}$. Then it is immediate that $F\in \Psi^0_\maL (P|\cG^{(0)}; \widetilde \maE)^\cG$ is odd with respect to  the grading $\Gamma$  and fiberwise elliptic.
Construct the invertible operator $U_D \in M_2 \left(\Psi_\maL^{-\infty} (P|\cG^{(0)}; \widetilde \maE)^\cG\right)$ by the same formula as before. Namely,
 choose a parametrix
$Q$ for $F$. Let $S_0=1-QF$, $S_1=1-FQ$. Then set $U_{D}=\begin{bmatrix} S_0 &-(1+S_0)Q\\F &S_1 \end{bmatrix}$. With such a choice the inverse is given by an explicit formula  $U_{D}^{-1}=\begin{bmatrix} S_0 &(1+S_0)Q\\-F &S_1 \end{bmatrix}$. Set $P_D=U_D^{-1}\begin{bmatrix} 1_{\maE} &0\\0 &0 \end{bmatrix}U_D$.

The choices in the constructions can be made so that  we have $(P_D)^\pm=P_{D^\pm}$, see  Definition~\ref{defindex}.
We then compute
\begin{multline*}
\ind_c (D^+) := \langle (\u)^{-\frac{ {\rm deg}\, c }{2}} \Lambda^\cE(c), \ch (\ind D^+) \rangle = \\ \frac{1}{2}\left(\langle (\u)^{-\frac{ {\rm deg}\, c }{2}} \Lambda^\cE(c), \Ch (\ind D^+) \rangle  - \langle (\u)^{-\frac{ {\rm deg}\, c }{2}} \Lambda^\cE(c), \Ch (\ind D^-) \rangle\right)
=\\ \frac{1}{2} \left \langle (\u)^{-\frac{ {\rm deg}\, c }{2}}
 \Lambda^{\widetilde{\cE}, \Gamma}(c), \left(\Ch P_D- \Ch \left(\begin{bmatrix} 0 &0\\0 &1_{\maE} \end{bmatrix} \right)\right) \right\rangle= \\
 \frac{1}{2} \left \langle (\u)^{-\frac{ {\rm deg}\, c }{2}}
 \Lambda_{\widetilde{\A}}(c), \left(\Ch P_D- \Ch \left(\begin{bmatrix} 0 &0\\0 &1_{\maE} \end{bmatrix} \right)\right) \right\rangle,
\end{multline*}
where the last equality is a consequence of  Proposition \ref{transgress}. Since $U_D$ is an even element of the algebra $\widetilde{\cA}$,
$\Ch(P_D) = \Ch \left(\begin{bmatrix} 1_{\maE} &0\\0 &0 \end{bmatrix} \right)$ in $\HC^-_0(\cA)$, and therefore
\begin{equation*}
\langle (\u)^{-\frac{ {\rm deg}\, c }{2}} \Lambda^\cE(c), \Ch (\ind D^+) \rangle = \frac{1}{2}\left \langle (u^*)^{-\frac{ {\rm deg}\, c }{2}}\Lambda_{\widetilde{\A}} (c), \left(\Ch \left(\begin{bmatrix} 1_{\maE} &0\\0 &0 \end{bmatrix} \right)- \Ch \left(\begin{bmatrix} 0 &0\\0 &1_{\maE} \end{bmatrix} \right)\right)   \right \rangle.
\end{equation*}

To complete the proof, we observe that
\begin{multline*}
\left \langle (u^*)^{-\frac{ {\rm deg}\, c }{2}}\Lambda_{\widetilde{\A}} (c),  \Ch  \left(\begin{bmatrix} 1_{\maE} &0\\0 &0 \end{bmatrix} \right)\right \rangle=\\
-\left \langle (u^*)^{-\frac{ {\rm deg}\, c }{2}}\Lambda_{\widetilde{\A}} (c),  \Ch \left(\begin{bmatrix} 0 &0\\0 &1_{\maE} \end{bmatrix} \right)    \right \rangle=
\left \langle (u^*)^{-\frac{ {\rm deg}\, c }{2}}\Lambda_{\A} (c),  [1]  \right \rangle
\end{multline*}
\end{proof}

Theorem \ref{superconnection computes index} allows by the standard method, see \cite{Gorokhovsky}, to deduce the higher twisted index theorem for foliations. For the sake of clarity we gather in the next paragraph  the main steps of the proof with the explicit  formula we get.

 \subsection{Proof of the higher index theorem}

The HKR-Connes morphism $ \chi \colon (\maC^\bullet(Z), d) \to (\CC^\bullet(C_c^\infty(Z), B)$ for an arbitrary manifold $Z$ is given by formula
\begin{equation*}
\chi (\xi) (a_0, a_1, \ldots, a_k) = \frac{1}{k!} \langle \xi, a_0da_1\ldots da_k \rangle
\end{equation*}
for a current $\xi$ and then extended $\u$-linearly.
It induces for our proper and free $\cG$-compact manifold $P$, a morphism,  denoted
$$
\chi_* \colon C^\bullet (\cG, \maC^\bullet(P)) \longrightarrow C^\bullet (\cG, \CC^\bullet(C_c^\infty (P)))\,.
$$
Composing the morphism $\chi$ with the integration quasi-isomorphism
$$
\int_\Delta \colon \maC_\Delta(\cP) \longrightarrow C^\bullet (\cG, \maC^\bullet(P))\,,
$$
we obtain a morphism
\begin{equation*}
\Xi:= \chi_* \circ \int_\Delta \colon \maC_\Delta(\cP) \longrightarrow C^\bullet (\cG, \CC^\bullet(C_c^\infty(P)))
\end{equation*}
Since the action of $\cG$ on $P$ is proper the inclusion
$$
\maC^\bullet(P/\cG)= \maC^\bullet(P)^\cG \to C^\bullet (\cG, \maC^\bullet(P))
$$
is a quasiisomorphism. Let $R \colon C^\bullet (\cG, \maC^\bullet(P)) \to \maC^\bullet(P/\cG)$ be a  (right) inverse quasiisomorphism constructed explicitely  using for instance Lemma \ref{proper}.
Composing the morphism induced by   $\chi_*$ on cohomology  with the composite map
\begin{equation*}
H^\bullet(C^\bullet (\cG, \CC^\bullet(C_c^\infty(P)))) \stackrel{\alpha_P}\longrightarrow \HC^\bullet(C_c^\infty(P) \rtimes \cG) \stackrel{\beta_P}\longrightarrow \HC^\bullet(C_c^\infty(P/\cG))
\end{equation*}
we obtain a map $\widehat{\chi}\colon H^\bullet(C^\bullet (\cG,  \maC^\bullet(P))) \longrightarrow \HC^\bullet(C^\infty(P/\cG))$ which is more precisely induced on cohomology by $\beta_P\circ \alpha_P\circ \chi_*$. We also have the similar maps:
\begin{equation*}
H^\bullet(C^\bullet (\cG, \CC^\bullet(\maA))) \stackrel{\alpha_\maA}\longrightarrow \HC^\bullet(\maA \rtimes \cG) \stackrel{\beta_\maA}\longrightarrow \HC^\bullet(\maA^\cG),
\end{equation*}
as well as the corresponding ones for the ideal $\Psi_\maL^{-\infty} (P\vert \maG^{(0)}; \maE)$ that we denote by $\alpha_{\Psi^{-\infty}}$ and $\beta_{\Psi^{-\infty}}$. Recall that $i:\Psi_\maL^{-\infty} (P\vert \maG^{(0)}; \maE) \hookrightarrow \maA$ and set $I: C_c^\infty (P)  \hookrightarrow \maA$, then we have
$$
i^*\circ (\beta_\maA\circ \alpha_\maA)= (\beta_{\Psi^{-\infty}}\circ \alpha_{\Psi^{-\infty}})\circ i^* \text{ and } I^*\circ (\beta_\maA\circ \alpha_\maA)= (\beta_{P}\circ \alpha_{P})\circ I^*.
$$

\begin{lemma} \label{integration} Let $\varrho$ be a cocycle in  $C^\bullet (\cG,  \maC^\bullet(P))$. Then
 \begin{equation*} \langle \widehat{\chi}([\varrho]), [1] \rangle = \int_{P/\cG} R(\varrho):=\langle R(\varrho), 1 \rangle .\end{equation*}
\end{lemma}
\begin{proof}
Since both sides of the equality depend only on the cohomology class $[\varrho]$, we can assume that $\varrho$ is $\cG$-invariant, in which case the result is clear.
\end{proof}
\medskip

Let now $\widehat{A}(P|\cG^{(0)})\, \Ch_L (\maE/S) \in \Omega_\bullet^\Delta(\cP)$ denote the form whose component
in $\Omega^*(\cP^{(k)} \times \Delta^k)$ is  $\widehat{A} (\cP^{(k)}|\cG^{(k)}) \Ch_{\maL^{(k)}} (\cE^{(k)}_\Delta/S)$.
This is a closed form in $(\Omega_\bullet^\Delta(\cP), d^t_{\varOmega(L)})$. Therefore if $c$ is a cocycle in $(\maC^\bullet_\Delta(\cP)[(\u)^{-1}, \u], d^t_{\varOmega(L)})$ then the pairing already described yields  a cocycle
$\pi^*c\cap (\, \widehat{A}(P|\cG^{(0)}) \Ch_L(\maE/S)\,)$  in the untwisted complex $(\maC^\bullet_\Delta(\cP)[(\u)^{-1}, \u], d).$ Hence, we can integrate it  and eventually define, when $c$ is an even cocycle, the higher twisted characteristic number  (with the appropriate normalization constant):
\begin{multline*}
\left<\pi^*c ,   \widehat{A}(P|\cG^{(0)}) \Ch_L(\maE/S) \right> :=\\  (2 \pi i )^{-l/2} \left(\int_{P/\cG} \circ R\circ  \int_\Delta \right) \left((\u)^{-\frac{ {\rm deg}\, c +l }{2}}(\pi^*c)\cap \widehat{A}(P|\cG^{(0)})\, \Ch_L(\maE/S)\right).
\end{multline*}
The classical short time limit of the heat operator associated with our Bismut superconnection yields:
\begin{proposition}\label{short time limit} Let $c$ be a cocycle in $(\maC^\bullet_\Delta(\cG)[(\u)^{-1}, \u], d^t_{\varOmega(L)})$.
Then
\begin{equation*}[I^* T_p\Phi_{\A}(c) ] =(2 \pi i \u)^{-l/2}   \left[{\Xi}\left((\pi^*c)\cap \widehat{A}(P|\cG^{(0)}) \Ch_L(\maE/S)\right)\right].\end{equation*}
where $l = \dim P - \dim \cG^{(0)}$ is the fiberwise dimension and the operator $T_p$ is defined in \eqref{deftpf}.
\end{proposition}

We are now in position to state our main theorem:

\begin{theorem} [The higher twisted index theorem for foliations] \label{thm:main}\
 Let $c$ be an even-degree cocycle in $(\maC^\bullet_\Delta(\cG)[(\u)^{-1}, \u], d^t_{\varOmega(L)})$. Then the higher twisted index is given by:
$$
\ind_c (D^+) =\left<\pi^*c ,   \widehat{A}(P|\cG^{(0)}) \Ch_L(\maE/S) \right>.
$$
\end{theorem}

\begin{proof}
We simply combine  the results of Theorem \ref{superconnection computes index}, Lemma \ref{integration} and Proposition \ref{short time limit}. We can rewrite the expression $\left<\pi^*c ,   \widehat{A}(P|\cG^{(0)}) \Ch_L(\maE/S) \right> $ as
\begin{multline*}
(2 \pi i )^{-l/2}  \left< \left(R\circ  \chi_*^{-1} \circ \Xi \right) \left((\u)^{-\frac{ {\rm deg}\, c +l }{2}}(\pi^*c)\cap \widehat{A}(P|\cG^{(0)})\, \Ch_L(\maE/S)\right), 1\right>\\=  \left< (\u)^{-\frac{ {\rm deg}\, c}{2}}(R\circ  \chi_*^{-1}) (I^* T_p\Phi_{\A}(c)), 1\right>.
\end{multline*}
Applying Proposition \ref{short time limit} to $\varrho:= \chi_*^{-1} (I^* T_p\Phi_{\A}(c)))$ we deduce
\begin{multline*}
\left<\pi^*c ,   \widehat{A}(P|\cG^{(0)}) \Ch_L(\maE/S) \right> =   \left<(\u)^{-\frac{ {\rm deg}\, c}{2}}(\beta_P\circ \alpha_P) (I^*T_p\phi_\A (c)), [1]\right>\\= \left< (\u)^{-\frac{ {\rm deg}\, c}{2}} (\beta_\maA\circ \alpha_\maA) (T_p\phi_\A (c)), [1]\right> =  \left<  (\u)^{-\frac{ {\rm deg}\, c}{2}} \Lambda_\A (c), [1]\right>.
\end{multline*}
Therefore, the proof is complete since it is reduced to  Theorem \ref{superconnection computes index}.
\end{proof}

\bigskip

\appendix

\section{Twisted cohomology}\label{TwistedCohom}

For a smooth manifold $M$ let $\varOmega \in \Omega^3(M)$ be a closed  $3$-form. Denote
by $u$ a formal variable of degree $-2$. The twisted de Rham complex is defined as the
complex $\Omega^*(M)[u]$ with the differential $d_{\varOmega}=ud +u^2\varOmega\wedge \cdot$. Its homology is denoted by $H_\Omega^\bullet (M)$.
Note that if $\varOmega'=\varOmega+d\eta $ is cohomologous to $\varOmega$ then the complexes $\left(\Omega^*(M)[u], d_{\varOmega} \right)
$  and $\left(\Omega^*(M)[u], d_{\varOmega'} \right)
$ are isomorphic via the isomorphism
\begin{equation}\label{defI}
I_{\eta} \colon \xi \mapsto e^{-u\eta}\wedge \xi
\end{equation}
We have more precisely  \cite{BenaGoro}:

\begin{lemma}\label{Imap}
The map induced by $I_{\eta}$ on cohomology depends only on the class of $\eta$ in
$\Omega^2(M)/\, d\Omega^1(M)$.
\end{lemma}

The identity $d_{\varOmega+\varOmega'} (\xi \wedge \eta) =  d_{\varOmega}\xi \wedge \eta + (-1)^{|\xi|} \xi \wedge d_{\varOmega'}  \eta$ implies that the product of forms induces the product $H^{\bullet}_{\varOmega}(M) \otimes H^{\bullet}_{\varOmega'}(M) \to H^{\bullet}_{\varOmega+\varOmega'}(M)$ and in particular endows
 $H^{\bullet}_{\varOmega}(M)$ with the structure of   $H^{\bullet}(M)[u]$-module.

We will also need to consider the following situation. Let $\pi \colon M \to B$ be an oriented fibration with compact fibers.
Then we have an integration along the fibers map $ \int_{M|B}\colon \Omega^*(M) \to \Omega^{*-k}(B)$ such that  $ \int_{M|B} d \eta = d \int_{M|B} \eta$, where $k=\dim M - \dim B$. Let $\varOmega \in \Omega^3(B)$ be a closed form. Then $\int_{M|B} \pi^*\varOmega \wedge \eta = \varOmega\wedge \int_{M|B} \eta$ and therefore
\begin{equation*}
\int_{M|B} d_{\pi^*\varOmega} \eta = d_{\varOmega} \int_{M|B} \eta.
\end{equation*}
Hence we obtain a chain map $\int_{M|B} \colon \left(\Omega^*(M){{[u]}}, d_{\pi^* \varOmega} \right)_{\bullet} \longrightarrow
\left(\Omega^*(B){{[u]}}, d_{\varOmega} \right)_{\bullet-k}$.

Let $\maL$ be a gerbe. The choice of connection defines a closed  $3$-form $\varOmega$.
We can therefore consider the complex $\left(\Omega^*(M)[u], d_{\varOmega} \right)$.
A different choice of connection leads to a different complex $\left(\Omega^*(M)[u], d_{\varOmega'} \right)$.
Lemma \ref{defomega} together with Lemma \ref{Imap}  implies that there exists however a canonical isomorphism
$I$ of the homologies of these complexes. We thus denote this (canonical isomorphism class of the) homology space by $H^{\bullet}_{\maL}(M)$.

\subsubsection{} \label{Currents} In the sequel it will be convenient to also work with regular currents and twisted de Rham homology. For a manifold $M$,  $\maC_l(M)$, the space of smooth currents of degree $l$, is given by $\maC_l(M) =\Omega^{\dim M -l}(M, or)$,
where $or$ is the orientation bundle. It
is a subspace of the vector space of distributions on $\Omega^{l}(M)$ with the pairing defined by the formula
\begin{equation*}
\langle\xi, \omega  \rangle=\int \xi \wedge \omega.
\end{equation*}

The differential $d^t \colon \maC_l(M) \to \maC_{l-1}(M)$ is  the transposed
of the de Rham differential up to a sign:
\begin{equation*}
\langle d^t \xi, \omega \rangle = (-1)^{\deg \xi+1} \langle\xi, d \omega \rangle, \quad \xi \in \maC_*(M), \omega \in \Omega^*(M).
\end{equation*}
We also have an operation of cap product $\maC_l(M) \otimes \Omega^m(M) \to \maC_{l-m}(M)$ defined by
\begin{equation*}
\langle  \xi \cap \alpha, \omega \rangle =   \langle\xi, \alpha \wedge \omega \rangle, \xi \in \maC_*(M), \alpha, \omega \in \Omega^*(M).
\end{equation*}
With these conventions we have:
\begin{equation*}
\xi \cap (\alpha \wedge \beta) = (\xi \cap \alpha) \cap \beta, \ \langle  \xi \cap \alpha, 1 \rangle= \langle\xi, \alpha   \rangle
\end{equation*}
and
\begin{equation}\label{dtLeibniz}
d^t(\xi \cap \alpha) = (d^t\xi) \cap \alpha+ (-1)^{\deg \xi} \xi \cap d \alpha
\end{equation}

Introduce the formal variable
$\u $ of degree $2$.
Form now the complexes   $\maC^\bullet(M) = \maC_*(M)[\u]$ with the differential of degree $1$ given by $\u d^t$.  If $\varOmega$ is a form of degree $3$ we can also introduce the differential $d_{\varOmega}^t$ by the formula
\begin{equation}\label{defdt}
d_{\varOmega}^t \xi = \u d^t \xi  - (\u)^2 \varOmega\cap \xi =  \u d^t \xi + (-1)^{\deg \xi+1}(\u)^2\xi \cap \varOmega.
\end{equation}

  We define a complex-valued pairing between  $\maC^\bullet(M)[(\u)^{-1}]$  and $\Omega^*(M)[u]$ by
\begin{equation*}
\langle (\u)^k \xi, u^l \omega  \rangle = \begin{cases} \langle \xi, \omega \rangle &\text{ if } k+l =0\\
0 &\text{ otherwise}  \end{cases}
\end{equation*}
The cap-product extends by
\begin{equation*}
((\u)^k \xi) \cap (u^l \omega) =(\u)^{k+l} (\xi \cap \omega).
\end{equation*}

We then have
\begin{equation*}
(\u d^t)(\xi \cap \alpha) = (d_{\varOmega}^t\xi) \cap \alpha+ (-1)^{\deg \xi} \xi \cap d_{\varOmega} \alpha.
\end{equation*}

\section{Review of the map $\Phi$}\label{cross}
 If $\cA$ is a sheaf of algebras over $\cG$, one can form a $\cG$-sheaf
$CC^\bullet\left(\cA\right)$, defined as the sheafification of the presheaf $U
\mapsto CC^*\left(\Gamma_c\left(U;\cA\right)\right)$. From now on assume that $\cA$ is a c-soft sheaf and we refer to such an $\cA$ as a $\cG-$algebra.
We can define a complex $C^\bullet\left(\cG,
CC^\bullet\left(\cA\right)\right)$. On the other hand one can form a
cross-product algebra $\cA\rtimes \cG$.

Assume $\lambda \in \Gamma(\cG^{(0)}; CC^\bullet\left(\cA\right))$ is $\cG$-invariant. Then one can construct a cyclic cocycle $c_\lambda$ on $\cA\rtimes \cG$ as follows.
Let $Z^{(k)}=
 \{(\gamma_1, \cdots, \gamma_k)\in \maG^{(k)} \mid  r(\gamma_1) = s(\gamma_k) \in \cG^{(0)}\}
$ be the subset of $\cG^{(k)}$. Let $\pi \colon Z^{(k)} \to \cG^{(0)}$ be the map given by
\[
\pi(\gamma_1, \cdots, \gamma_k)= r(\gamma_1).
\]
Note that we have the following sequence of maps of vector spaces:
\[
\Gamma_c(\cG^k; (r^*\cA)^{\boxtimes k}) \to    \Gamma_c(Z^{(k)}; (r^*\cA)^{\boxtimes k}) \to \Gamma_c(Z^{(k)}; (\pi^*\cA)^{\otimes k})\to
\Gamma_c(\cG^{(0)}; \cA^{\otimes k})
\]
Here the first map is induced by the inclusion $Z^{(k)} \hookrightarrow \cG^k$. The second map is induced by the isomorphism of sheaves $(r^*\cA)^{\boxtimes k}) \cong (\pi^*\cA)^{\otimes k}$ on $Z^{(k)}$ given by $\cG$ action on $\cA$. Finally the third map is induced by the map $\pi$. Denote the composition of these maps by $M_k$. Then we can define the cyclic cocycle $c_\lambda$ by
\[
c_\lambda(f_0, f_1, \ldots, f_k) = \lambda( M_{k+1}(f_0\otimes f_1\otimes\ldots\otimes f_k)).
\]

One can extend this construction to  produce a
morphism of complexes
\begin{equation}
\Phi \colon C^\bullet\left(\cG, CC^\bullet\left(\cA\right)\right) \to CC^\bullet\left(\cA \rtimes
\cG\right)
\end{equation}
This map has been constructed in \cite{crainic}, based on the
constructions in \cite{bryni}, compare also \cite{ni90}, \cite{ft87},
\cite{gj93}. An important particular case is when $A$ is the sheaf
$C^{\infty}(B)$ of smooth functions on $B$. In this case there is a
canonical morphism of complexes of sheaves $\iota:
\Omega_*\left(B\right) \to CC^*\left( C_0^{\infty}\left(B\right)\right)$
defined on the current $c$ of degree $m$ by
\begin{equation}\label{iota}
\iota(c) (a_0, a_1, \ldots, a_m) =\frac{1}{m!}\left \langle c,
a_0da_1 \ldots da_m\right \rangle.
\end{equation}
Composing $\iota$ with the map in \eqref{phi} one obtains (up to chain homotopy) Connes'
map \cite{ConnesTransverse, ConnesBook}, also denoted by $\Phi$:
\begin{equation}
\Phi: C^*\left(\cG, \Omega_*\left( B\right)\right) \to
CC^*\left(C_0^{\infty}(\cG)\right).
\end{equation}

We will use the following properties of the map $\Phi$, compare
\cite{crainic}:
\begin{itemize}
\item If $\maA$, $\maB$ two $\cG$-algebras and $f \colon \maA \to \maB$ is a
$\cG$-homomorphism, we get a natural induced homomorphism,  which we
also denote by $f$, from $\maA \rtimes \cG$ to $\maB \rtimes \cG$. Also the
induced map $f^* :CC^*\left(\maB\right) \to CC^*\left(\maA\right)$ is
$\cG$-equivariant, and hence defines a map, also denoted by $f^*$,
from $C^*\left(\cG, CC^*\left(\maB\right)\right)$ to $C^*\left(\cG,
CC^*\left(\maA\right)\right)$. Then the following diagram commutes:
\begin{equation*}
\begin{CD}
 C^\bullet \left(\cG, CC^\bullet\left(\maB\right)\right)  @>\ \ \Phi \ \ >>
CC^\bullet\left(\maB\rtimes \cG\right)\\
@V f^*   VV    @VVf^* V \\
 C^\bullet \left(\cG, CC^\bullet\left(\cA\right)\right)  @>\ \ \  \Phi \ \ \ >>  CC^\bullet\left(\maA\rtimes \cG\right)
\end{CD}
\end{equation*}

\item If $\lambda \in \Gamma(\cG^{(0)}; CC^\bullet\left(\cA\right))$ is $\cG$-invariant then
\[
\Phi(\lambda)=c_\lambda
\]
\end{itemize}

\section{Some JLO estimates for complete manifolds}\label{JLO}
Let $ \| \cdot \|$ denote the operator norm of $B L^2(N, E)$.
\begin{lemma} \label{lem:sub} Let $N$ be a complete Riemannian manifold and $D$  a Dirac type operator on $N$ acting on a Clifford module $E$.
Assume that the $V_i$ ($0\leq i \leq l$) are pseudodifferential operators on $P$  such that
\begin{equation*}
\|V_i\|_{v_i, t} = \|(1+D^2)^{-t/2} V_i (1+D^2)^{(-v_i+t)/2}\|<\infty, \,  v_i=\ord V_i
\end{equation*}
for every $t \in \mathbb{R}$ and every $i$.
Let $q=\max \left \{\Sigma_{i=0}^l \ord V_i, 0\right\} $. Then for $0 <s \le 1$
\begin{equation*}
\left \| \int \limits_{\Delta^l} V_0 e^{-\sigma_0 s^2 D^2} V_1 e^{-\sigma_1 s^2D^2} \ldots V_l e^{-\sigma_l s^2 D^2} d\sigma_1\ldots d\sigma_l\right\| \le C(q) \|V_0\|_{v_0, t_0}\ldots \|V_l\|_{v_l, t_l}  s^{-q}
\end{equation*}
for some $t_i$.

\end{lemma}
\begin{proof}

Using barycentric subdivision, decompose the integral into the sum of integrals over $l+1$ nonintersecting domains $U_i \subset \Delta^l $ such that
$\sigma_i \ge 1/l$ on $U_i$. Notice that we have the following identity:
\begin{multline*}
V_0 e^{-\sigma_0 s^2 D^2} V_1 e^{-\sigma_1 s^2D^2} \ldots V_l e^{-\sigma_l s^2 D^2} = \\
(V_0(1+D^2)^{-v_0/2})e^{-\sigma_0 s^2 D^2}( (1+D^2)^{v_0/2}V_1(1+D^2)^{-(v_0/2+v_1/2)})e^{-\sigma_1 s^2 D^2} \\
\ldots(1+D^2)^{(v_0/2+\ldots +v_{i-1}/2)} V_i (1+D^2)^{-(v_0/2+\ldots +v_i/2)} ((1+D^2)^{v_0/2+\ldots+v_l/2} e^{-\sigma_i s^2 D^2}) \ldots\\
((1+D^2)^{-(v_{l-1}/2+v_l/2)}V_{l-1}(1+D^2)^{v_l/2})e^{-\sigma_{l-1} s^2 D^2}( (1+D^2)^{-v_l/2}V_le^{-\sigma_l s^2 D^2})e^{-\sigma_l s^2 D^2}
\end{multline*}
On $U_i$, using the estimate
\begin{equation*}
\| (1+D^2)^{v_0/2+\ldots +v_l/2} e^{-\sigma_i s^2 D^2} \| \le c( q)l^{q/2} s^{-q} \text { when }  \sigma \ge \frac{1}{l}, \ 0<s \le 1.
\end{equation*}
we bound the integrand, and obtain the inequality
\begin{equation*}
\left \| \int_{U_i} V_0 e^{-\sigma_0 s^2 D^2} V_1 e^{-\sigma_1 s^2D^2} \ldots V_l e^{-\sigma_l s^2 D^2} d\sigma_1\ldots d\sigma_l\right\| \le c(q) \|V_0\|_{v_0, t_0}\ldots \|V_l\|_{v_l, t_l} \frac{l^{q/2} s^{-q}}{l!}
\end{equation*}
The inequality follows since $\frac{l^{q/2} }{l!}$ can be bounded uniformly in $l$.
\end{proof}

Choose any $p$ such that    $(1/2 \dim N +1/2) > p > 1/2 \dim N$ and let $Q$ be a properly supported parametrix for $(1+D^2)^p$, with
$
1-Q(1+D^2)^p =R
$
a properly supported smoothing operator. For a zero order pseudodifferential operator $V$ with  compactly supported Schwartz kernel
let $\mu(V) = \|V Q\|_1+ \|VR\|_1$ where $\| \cdot \|_1$ denotes the
trace norm.
\begin{lemma}\label{trace estimate}
 Let $V_i$ be as in Lemma \ref{lem:sub}. Assume also that $V_0$ is a zero order pseudodifferential operator  with  compactly supported Schwartz kernel. Then
\begin{multline*}
  \left \| \int \limits_{\Delta^l} V_0 e^{-\sigma_0 s^2 D^2} V_1 e^{-\sigma_1 s^2D^2} \ldots V_l e^{-\sigma_l s^2 D^2} d\sigma_1\ldots d\sigma_l \right\|_1  \le \\ C'(q) \mu(V_0) \|V_1\|_{v_1, t_1}\ldots \|V_l\|_{v_l, t_l}   s^{-q - \dim N - 1}
\end{multline*}
for some $t_i$.
\end{lemma}

\begin{proof}
Notice that
\[ \int \limits_{\Delta^l} V_0 e^{-\sigma_0 s^2 D^2} V_1 e^{-\sigma_1 s^2D^2} \ldots V_l e^{-\sigma_l s^2 D^2} d\sigma_1\ldots d\sigma_l =  V_0Q  \cdot (1+D^2)^{p} P + V_0 R \cdot P\]
where
\[P =\int \limits_{\Delta^l}  (1+D^2)^{p} e^{-\sigma_0 s^2 D^2} V_1 e^{-\sigma_1 s^2D^2} \ldots V_l e^{-\sigma_l s^2 D^2} d\sigma_1\ldots d\sigma_l
\]
Now $(1+D^2)^{p} P$ and $P$ are bounded operators, whose norms can be estimated by Lemma \ref{lem:sub}.

\end{proof}


\end{document}